\documentclass{amsart}

\usepackage{kotex}
\usepackage{amsmath}
\usepackage{mathrsfs}
\usepackage{amsfonts}
\usepackage{amscd}
\usepackage{slashed}
\usepackage{stmaryrd}

\usepackage{amssymb}
\usepackage{graphicx}
\usepackage{caption}
\usepackage{subcaption}
\usepackage{dsfont}
\usepackage{color}
\usepackage{hyperref}
\usepackage{bbm}
\usepackage{a4wide}
\usepackage{sseq}
\usepackage{tikz-cd}
\usepackage[abs]{overpic} 

\newtheorem{theorem}{Theorem}[section]
\newtheorem{lemma}[theorem]{Lemma}
\newtheorem{proposition}[theorem]{Proposition}

\newtheorem{example}[theorem]{Example}

\theoremstyle{definition}
\newtheorem*{definition*}{Definition}
\newtheorem{definition}[theorem]{Definition}

\theoremstyle{remark}
\newtheorem*{remark*}{Remark}
\newtheorem{remark}[theorem]{Remark}

\theoremstyle{question}
\newtheorem*{question*}{Question}

\numberwithin{equation}{section}

\newcommand{\myproof}[2]{Proof of {#1} {#2}}

\begin{document}

\title[Deformation quantization of coisotropic A-branes]{Deformation quantization of coisotropic A-branes}

\author{Naichung Conan Leung}
\address{The Institute of Mathematical Sciences and Department of Mathematics, The Chinese University of Hong Kong, Shatin, Hong Kong}
\email{leung@ims.cuhk.edu.hk}

\author{YuTung Yau}
\address{Yau Mathematical Sciences Center, Tsinghua University, Beijing, China}
\email{ytyau@tsinghua.edu.cn}

\author{Zhenyu Ye}
\address{The Institute of Mathematical Sciences and Department of Mathematics, The Chinese University of Hong Kong, Shatin, Hong Kong}
\email{zyye@math.cuhk.edu.hk}

\thanks{}

\maketitle


\begin{abstract}
	Motivated by the Gukov--Witten proposal that the endomorphism algebra of a coisotropic A-brane should be described by a deformation quantization of its transverse holomorphic symplectic geometry at the perturbative level, we study deformation quantization on transverse holomorphic symplectic manifolds $(M,\Omega)$. Extending the approach of Nest--Tsygan to this setting, we show that, under an admissibility condition, equivalence classes of deformation quantizations are in bijection with
	\[
	\hbar H^2\bigl(\operatorname{F}_1\Omega^\bullet(M,\mathbb C),d\bigr)
	[[\hbar]],
	\]
	where $\operatorname{F}_1\Omega^\bullet(M,\mathbb C)$ is the differential ideal generated by transverse $(1,0)$-forms.
	
	We further illustrate this framework through examples arising from mapping tori, $T$-duality, and coisotropic lifts in cotangent bundles. These examples clarify the roles of admissibility and the transverse Atiyah class, as well as the distinction between deformation quantizations and global star products. Finally, for a coisotropic A-brane $M$ induced from a transverse complex manifold $Q$, we use transverse jets and the Fock representation to construct a perturbative level action $End(M)\curvearrowright Hom(M,Q)$ under suitable compatibility assumptions, providing
	a geometric realization of the brane quantization philosophy.

\end{abstract}

\section{Introduction}
Deformation quantization originated from quantizing classical observables on a symplectic phase space. Intriguingly, it also appears in mirror symmetry through the Gukov--Witten \emph{brane quantization} proposal \cite{GukWit2009}, which relates quantization to endomorphism algebras of \emph{coisotropic A-branes}. Introduced by Kapustin--Orlov \cite{KapOrl2003}, such branes are expected to belong to an enlarged Fukaya category required for homological mirror symmetry \cite{Kon1995}. 

Although a general mathematical construction of this enlarged category
is still lacking, the Gukov--Witten proposal provides a physical prediction
for the endomorphism algebra of a space-filling coisotropic A-brane and its module. Several
works have implemented aspects of this proposal in the study of coisotropic A-branes\cite{AENLV2014, ChaLeuLiSueYau2026, GukKorNawPeiSab2023, Kuw2022, Qiu2023}. More generally, coisotropic A-branes need not be space-filling. A simple
model is a brane supported on $L\times Y\subset X_0\times Y$, where
$L$ is Lagrangian inside a symplectic space $X_0$ and $Y$ is holomorphic symplectic, equipped
with suitable brane data. 

Non-space-filling coisotropic branes arise in various settings, such as topological Chern--Simons theory couple with matter \cite{AgaCosMcnVaf2017}, SYZ transformation \cite{ChaLeuZha2018} and brane reduction\cite{qcr}.
Nevertheless, a general mathematical framework for the endomorphism algebra of a coisotropic A-brane and its module remains to be developed.

The present work aims to fill this gap by developing a systematic theory of deformation quantization and its module for the sheaf of transverse holomorphic functions on a transverse holomorphic symplectic manifold $(M,\Omega)$, and establishing corresponding existence and classification results.  

In the setting of coisotropic A-branes which naturally carries transverse holomorphic symplectic structure, the endomorphism algebra of a brane $\mathcal{B}=(M,\Omega)$ is modeled by one such deformation quantization. Our formulation of deformation quantization part is sheaf-theoretic, following the viewpoint of Nest--Tsygan \cite{NesTsy2001}. The existence and classification results are then established using Fedosov's approach \cite{Fed1994}, which has been extended to the holomorphic and algebraic symplectic settings in \cite{BezKal2004,NesTsy2001,Vai2002,Yek2008}. Within this framework, we give an explicit construction of modules over such deformation quantizations in the transverse complex manifold setting. We expect the resulting module action to realize the composition of brane morphisms $\operatorname{End}(\mathcal B)\curvearrowright\operatorname{Hom}(\mathcal L,\mathcal B)$ when $\text{supp}(\mathcal{L})\subset \text{supp}(\mathcal{B})$.

\subsection{The main results}
\quad\par
A transverse holomorphic symplectic manifold \((M,\Omega)\) carries the distribution \(A:=\ker\Omega\subset T_{\mathbb C}M\) which determines its transverse holomorphic structure. We denote by \(\mathcal O\) the sheaf of smooth complex-valued functions annihilated by $A$, regarded as the structure sheaf of $(M, \Omega)$. In Definition \ref{dqofo}, a \emph{deformation quantization} of $(M, \Omega)$ is a sheaf \(\mathcal A^\hbar_M\) of flat, \(\hbar\)-adically complete \(\mathbb C[[\hbar]]\)-algebras whose classical limit is $\mathcal{O}$ and whose semiclassical commutator recovers the Poisson bracket induced by $\Omega$. We further require \(\mathcal A_M^\hbar\) to be locally equivalent to a transverse holomorphic star product.

To study existence and classification of such deformation quantization, we use Fedosov's method on the Weyl bundle over $M$,
\[
W = \widehat{\operatorname{Sym}} (B^*)[[\hbar]],
\qquad
B = T_\mathbb{C} M / A.
\]
A Fedosov connection is a flat connection on $W$ compatible with its fibrewise Weyl algebra structure (see Definition \ref{Definition: Fedosov connection} for its precise meaning). Its sheaf of flat sections defines a deformation quantization of $(M, \Omega)$, which we call a \emph{Fedosov quantization}. Our first main result shows that every deformation quantization arises in this way:

\begin{theorem}[$=$ Theorem \ref{dqtofedosov}]
	Let \((M,\Omega)\) be a transverse holomorphic symplectic manifold. Then every deformation quantization of \((M,\Omega)\) is equivalent to a Fedosov quantization.
\end{theorem}

Associated with each Fedosov connection is a closed formal $2$-form in $-\Omega + \hbar \, \Omega^2(M,\mathbb C)[[\hbar]]$, whose de Rham cohomology class depends only on the equivalence class of the resulting Fedosov quantization. This defines the \emph{non-commutative period map}
\[
\operatorname{per}: \operatorname{DQ}(M, \Omega) \to -[\Omega] + \hbar H^2(M, \mathbb{C})[[\hbar]],
\]
where $\operatorname{DQ}(M, \Omega)$ denotes the set of equivalence classes of deformation quantizations of $(M, \Omega)$.

Our existence and classification results are obtained by studying the image and fibres of this period map under a natural admissibility condition. We call \((M,\Omega)\) \emph{admissible} if the following natural restriction maps are surjective:
\[
H^i(M,\mathbb C)\longrightarrow 
H^i(M, \mathcal{O}),
\qquad i=1,2.
\]
Let \(F_*\Omega^*(M,\mathbb C)\) denote the Leray--Cartan filtration associated with the transverse holomorphic structure; in particular, \(F_1\Omega^*(M,\mathbb C)\) is the differential ideal generated by transverse \((1,0)\)-forms. We summarize the main results as follows; see Theorems~\ref{Theorem 3.18} and~\ref{Theorem: classification} for the precise statements.

\begin{theorem}
	Let \((M,\Omega)\) be a transverse holomorphic symplectic manifold.
	\begin{enumerate}
		\item \textbf{Existence.} If the restriction map $H^2(M,\mathbb C) \to 
		H^2(M, \mathcal{O})$ is surjective, then there exists a deformation quantization of \((M,\Omega)\).
		\item \textbf{Classification.} If $(M, \Omega)$ is admissible, then the noncommutative period map is injective and its image is in bijection with
		\[
		\hbar H^2\bigl(
		\operatorname F_1\Omega^*(M,\mathbb C),d
		\bigr)[[\hbar]].
		\]
	\end{enumerate}
\end{theorem}

Thus, for admissible transverse holomorphic symplectic manifolds, the noncommutative period map provides the complete classification of deformation quantizations.

For our final main result, we consider a coisotropic A-brane $\mathcal B$ and a Lagrangian A-brane $\mathcal L$ in a symplectic manifold $X$, such that $C:=\operatorname{supp}(\mathcal B) \supset Q:=\operatorname{supp}(\mathcal L)$.
Suppose that $Q$ is a transverse holomorphic Lagrangian submanifold
of $C$ and a section of a transverse holomorphic Lagrangian fibration
$\pi:C\to Q$; see Definition~\ref{transverse holomorohic lagrangian} and \ref{Definition5.11}. In this setting, we construct the action
\[
\operatorname{End}(\mathcal B)\curvearrowright
\operatorname{Hom}(\mathcal L,\mathcal B)
\]
expected from composition of brane morphisms. We expect that the techniques developed in \cite{BGKP2016} can be
adapted to extend our construction, under suitable conditions, to cases where no such transverse complex Lagrangian fibration assumption exists.

Our construction arises from transverse complex geometry. Let $(Q,A')$ be a
transverse complex manifold with characteristic distribution
$A'_{\mathbb R}$. Then $M:=(TQ/A'_{\mathbb R})^*$ carries a natural
transverse holomorphic symplectic structure $\Omega$. Write
$j:Q\hookrightarrow M$ for the zero section. Via the natural inclusion
$M\subset T^*Q$, we regard $M$ as the support of the coisotropic A-brane
$\mathcal B$ and the zero section $Q\subset T^*Q$ as the support of the
Lagrangian A-brane $\mathcal L$. The following is a simplified statement
of our result; see Theorem~\ref{Theorem: lift of flat connections} for
the precise statement.

\begin{theorem}
	Let $(Q,A')$ be a transverse complex manifold, and let $(M,\Omega)$ be
	the associated transverse holomorphic symplectic manifold constructed
	above. Write $j:Q\hookrightarrow M$ for the zero section. Then there
	exist a deformation quantization $(\mathcal A_M^\hbar,\psi)$ of
	$(M,\Omega)$ and a left $\mathcal A_M^\hbar$-module structure on
	$j_*\mathcal O_{A'}[[\hbar]]$.
\end{theorem}

\subsection{Relations to mirror symmetry through coisotropic A-branes}
\quad\par
A coisotropic A-brane $\mathcal B$ in a symplectic manifold $(X, \omega_X)$ is a coisotropic submanifold $M \subset X$ equipped with a Hermitian line bundle $L$ with a unitary connection $\nabla^L$ of curvature $-2\pi\sqrt{-1} F_L$ such that $\Omega = F_L + \sqrt{-1} \omega_X \vert_M$ defines a transverse holomorphic symplectic structure on $M$. Conversely, every transverse holomorphic symplectic manifold can be realized as a coisotropic A-brane in an ambient symplectic manifold, possibly in the presence of a nontrivial B-field; see Proposition \ref{Proposition: THS as coisotropic A-branes}.

While a mathematical construction of the endomorphism algebra
\[
\operatorname{Hom}(\mathcal B,\mathcal B)
\]
in the conjectural enlarged Fukaya category remains elusive in general, its classical limit is better understood. The physical analyses of Kapustin--Orlov \cite{KapOrl2003} and Kapustin--Li \cite{KapLi2005} identify the corresponding BRST cohomology with the sheaf cohomology $H^*(M,\mathcal O)$.

Beyond the classical limit, considerably more is understood when $\mathcal B$ is space-filling, so that $M = X$. In this case, the endomorphism algebra is expected to be governed by holomorphic deformation quantization of $(X, \Omega)$ \cite{GukWit2009}. Moreover, worldsheet instanton corrections are absent in this space-filling setting \cite{Her2012, KapWit2007}, so the relevant quantum corrections are expected to be perturbative and local on $X$. This locality already points naturally toward a sheaf of noncommutative algebras deforming $\mathcal{O}$, a perspective implicit in the discussion of Kapustin--Witten \cite{KapWit2007} and developed further in \cite{GaiWit2022}.

For a general, possibly non-space-filling, coisotropic A-brane, the full endomorphism algebra is expected to involve both perturbative and non-perturbative effects. The perturbative corrections should deform the local algebra encoded by the transverse holomorphic geometry, whereas non-perturbative contributions are expected to introduce genuinely global phenomena and higher $A_\infty$-operations. From this viewpoint, the cohomological description $H^*(M,\mathcal O)$ should be regarded as the shadow of a richer cochain-level structure.

We view the present work as developing the perturbative part of this picture. We propose that the appropriate perturbative model for the endomorphism algebra of a general coisotropic A-brane $\mathcal{B}$ is a deformation quantization of its transverse holomorphic symplectic support $(M, \Omega)$. As a future direction, we plan to apply our framework to a perturbative study of the ``quantization commutes with reduction'' principle for endomorphism algebras of coisotropic A-branes \cite{qcr}.

\subsection{Relation to existing frameworks of deformation quantization}
\quad\par

Our approach occupies an intermediate position among existing frameworks for deformation quantization. We work with sheaves of \(\mathbb C[[\hbar]]\)-algebras on \((M,\Omega)\), without requiring a globally defined star product, while retaining enough geometric structure for an explicit Fedosov-type construction and a classification via the noncommutative period map.

The closest precedents are the constructions of Nest--Tsygan \cite{NesTsy2001} and Bezrukavnikov--Kaledin \cite{BezKal2004}, which deform the structure sheaf of a holomorphic symplectic space. Our setting is more general in a different direction: $M$ need not be a complex manifold, but only a transverse complex manifold, allowing foliations with non-closed or even dense leaves.

Vaisman \cite{Vai2002} develops a Fedosov-type construction for symplectic ringed spaces, a framework that encompasses transverse holomorphic symplectic manifolds. His construction, however, produces global star products and requires the vanishing of the Atiyah class. This condition fails in many natural examples, including K3 surfaces, and is therefore too restrictive for our purposes.

Bandiera--Sti\'enon--Xu \cite{BanStiXu2021} prove a formality theorem for arbitrary Lie pairs. Their framework does not directly formulate deformation quantization, nor does it provide corresponding existence and classification results in the form needed here. Although such structures should in principle be recoverable from formality, making this passage explicit is substantially more involved and less geometrically transparent than the Fedosov-type approach pursued here.

More general categorical formulations are provided by deformation-quantization algebroids and related stack-theoretic approaches \cite{BreRoj2022, KasSch2012, Kon2001}, in which local quantizations may glue only up to higher descent data. The present work remains at the sheaf level. While less general than these categorical frameworks, this level of structure is sufficient not only to parametrize deformation quantizations of a coisotropic A-brane $\mathcal{B}$, but also to construct geometrically a distinguished quantization of its support $(M, \Omega)$ through its action on a proposed morphism space
\[
\operatorname{Hom}(\mathcal L,\mathcal B),
\]
where $\mathcal L\subset\mathcal B$ is a Lagrangian A-brane satisfying suitable compatibility conditions; see Theorem~\ref{Theorem: lift of flat connections}. This provides an explicit realization of the brane-quantization philosophy.

\subsection{Organization of this paper and notations}
\quad\par
The paper is organized as follows. Section~\ref{Section: Geometry of transverse holomorphic structures} reviews the geometry of transverse holomorphic structures \cite{Gom1980}, providing the necessary geometric preliminaries. Section~\ref{Section3} introduces transverse holomorphic symplectic manifolds and formulates their deformation quantization. Section~\ref{section4} develops the Fedosov quantization method and establishes the existence and classification theorems for the deformation quantization of admissible transverse holomorphic symplectic manifolds. Finally, in Section~\ref{Section: Examples}, we present several explicit examples illustrating our deformation quantization framework and its relation to coisotropic A-branes. Building on one of these examples, we rigorously construct a perturbative module action in Subsection~\ref{section5.6}.

Throughout this paper, unless otherwise specified, let $(M, \Omega)$ be a transverse holomorphic symplectic manifold and $A$ be its underlying transverse holomorphic structure. We define the quotient bundle $B := T_{\mathbb C}M/A$. We adopt the following basic notation and conventions:
\begin{itemize}
	\item $C^{\infty}(M)$: the algebra of $\mathbb{C}$-valued smooth functions on $M$.
	\item $C_M^\infty$: the sheaf of $\mathbb{C}$-valued smooth functions on $M$.
	\item $\mathcal{O}$: the sheaf of $A$-invariant smooth complex-valued functions on $M$, defined on any open subset $U \subset M$ by 
	\[ \mathcal{O}(U) := \{ f \in C^{\infty}(U) \mid \mathcal{L}_a f = 0, \, \forall a \in \Gamma(U, A) \}. \]
	\item $d_A^{\mathrm{Bott}}$: the Bott connection on $B$, or the induced connection on $B^*$, $\widehat{\operatorname{Sym}} (B^*)$, etc.
	\item $\operatorname{Sym}(Q)$ (resp. $\widehat{\operatorname{Sym}}(Q)$): the symmetric algebra bundle (resp. completed symmetric algebra bundle) associated with a given complex vector bundle $Q$ over $M$.
	\item $JM$: the jet bundle of the trivial complex line bundle over $M$.
	\item $J_A M$: the transverse jet bundle with respect to $A$.
	\item $\nabla^{\operatorname{G}}$: the Grothendieck connection on the jet bundle.
	\item $\mathcal{J}$: the sheaf of sections of $J_A M$, viewed as an $\mathcal{O}$-module.
	
	\item $W$: the Weyl bundle $\widehat{\operatorname{Sym}}(B^*)[[\hbar]]$.
	\item $*$: the fibrewise Moyal product on $W$.
	\color{black}
\end{itemize}

\section{Geometry of transverse holomorphic structures}
\label{Section: Geometry of transverse holomorphic structures}
Recall the following notion, whose systematic study goes back to the work of G\'omez-Mont \cite{Gom1980}.

\begin{definition}
	A \emph{transverse holomorphic structure} on a smooth manifold $M$ is an involutive complex distribution $A \subset T_{\mathbb C}M$ such that $A + \overline{A} = T_{\mathbb C}M$.
\end{definition}

Throughout this section, $M$ always denotes a smooth manifold, $A$ is a transverse holomorphic structure on $M$, and $B$ denotes the quotient bundle $T_{\mathbb C}M/A$. We refer to the pair $(M, A)$ as a \emph{transverse complex manifold}, and the complex rank of $B$ as its \emph{transverse complex dimension}. Below, we review the geometry of $(M, A)$. The material is essentially due to G\'omez-Mont \cite{Gom1980}; our aim is simply to present it in a form and notation consistent with the exposition of this paper.

In Subsection \ref{Subsection: Nirenberg integrable distributions}, we recall the basic properties of $A$. In Subsection \ref{Subsection: Atiyah class}, we recall the notion of the Atiyah class of $(M, A)$. Subsection \ref{Subsection: Leray--Cartan filtration} is devoted to the Leray--Cartan filtration induced by $A$. In Section \ref{Subsection: transverse jet bundles}, we discuss the formal geometry for $(M, A)$.

\subsection{Transverse holomorphic structures}
\quad\par
\label{Subsection: Nirenberg integrable distributions}
Associated with $A$ are several natural sheaves and differential complexes that encode its transverse and leafwise geometry. Let $\mathcal{O}$ denote the sheaf of smooth complex-valued functions on $M$ constant along $A$, which we call the \emph{structure sheaf} of $(M, A)$. More precisely, for any open subset $U \subset M$,
\begin{equation*}
	\mathcal{O}(U) := \{ f \in \mathcal{C}^\infty(U, \mathbb{C}) \mid \mathcal{L}_a f = 0 \text{ for all } a \in \Gamma(U, A) \}.
\end{equation*}
The quotient bundle $B = T_{\mathbb C}M/A$ admits a canonical flat $A$-connection $d_A^{\operatorname{Bott}}$, called the \emph{Bott connection}, and defined as follows: for all $a \in \Gamma(M, A)$ and $v \in \Gamma(M, T_{\mathbb C}M)$,
\begin{equation}
	\label{Equation: Bott connection}
	(d_A^{\operatorname{Bott}} \mathfrak{q}(v))(a) = \mathfrak{q}( [a, v] ),
\end{equation}
where $\mathfrak{q}: T_{\mathbb C}M \to B$ denotes the natural projection. It yields the sheaf $\mathcal{B}$ of $d_A^{\operatorname{Bott}}$-closed sections of $B$, which is naturally an $\mathcal{O}$-module.\par

Since $A$ is a Nirenberg integrable complex distribution, the Nirenberg Complex Frobenius Theorem \cite{Nir1957} implies that every point of $M$ admits a neighborhood $U$ diffeomorphic to a product $U_1 \times U_2$, where $U_1 \subset \mathbb{R}^m$ and $U_2 \subset \mathbb{C}^{m'}$ are open subsets, such that
$$\left\{ \frac{\partial}{\partial x^1}, \dots, \frac{\partial}{\partial x^m}, \frac{\partial}{\partial \bar{z}^1}, \dots, \frac{\partial}{\partial \bar{z}^{m'}} \right\}$$
forms a frame of $A$ over $U$. Here, $(x^1, ..., x^m)$ and $(z^1, ..., z^{m'})$ denote the standard real and complex coordinates on $\mathbb{R}^m$ and on $\mathbb{C}^{m'}$, respectively. We shall refer to
$$ (x^1, \dots, x^m, z^1, \dots, z^{m'})$$
as \emph{adapted coordinates} on $U$. Since coordinate vector fields commute, \eqref{Equation: Bott connection} immediately shows that $\left\{ \mathfrak{q} \left( \tfrac{\partial}{\partial z^1} \right), ..., \mathfrak{q} \left( \tfrac{\partial}{\partial z^{m'}} \right) \right\}$ is a $d_A^{\operatorname{Bott}}$-closed frame of $B$ over $U$, hence $\mathcal{B}$ is a locally free $\mathcal{O}$-module.

For any (possibly infinite rank) complex vector bundle $E$ over $M$, we define:
\begin{equation*}
	\Omega_A^{i,j}(M, E) = \Gamma(M, \textstyle\bigwedge^j A^* \otimes \bigwedge^i B^* \otimes E)
\end{equation*}
When $E = M \times \mathbb{C}$ is the trivial line bundle, we simply write $\Omega_A^{i, j}(M)$. Let $d_A: \Omega_A^{0, *}(M) \to \Omega_A^{0, *}(M)$ denote the leafwise de Rham differential of $A$. Indeed, the leafwise de Rham cohomology, i.e. the cohomology of the complex $(\Omega_A^{0, *}(M), d_A)$, denoted by $H_A^*(M)$, is naturally isomorphic to the sheaf cohomology $H^*(M, \mathcal{O})$.

\subsection{Atiyah class}
\quad\par
\label{Subsection: Atiyah class}
By the \emph{Atiyah class} of the transverse complex manifold $(M,A)$, we mean the Atiyah class of the Lie pair $(T_{\mathbb C}M,A)$ (see, for example, \cite[Section 2]{VoglaireXu2015}). We briefly recall its construction.

Let $\nabla$ be a $T_{\mathbb C}M$-connection on $B$ extending the Bott connection $d_A^{\operatorname{Bott}}$, i.e. $\nabla_a b = (d_A^{\operatorname{Bott}} b)(a)$ for all $a \in \Gamma(M, A)$ and all $b \in \Gamma(M, B)$. Then it defines a $d_A^{\operatorname{Bott}}$-closed $1$-form
\[
R \in \Omega_A^{0, 1}(M, B^* \otimes \operatorname{End}(B))
\]
given by
\[
R(a)(\mathfrak{q}(X)) := \nabla_a\nabla_X - \nabla_X\nabla_a - \nabla_{[a, X]} \in \Gamma(M, \operatorname{End}(B)),
\]
for all $a \in \Gamma(M, A)$ and $X \in \Gamma(M, T_{\mathbb C}M)$. We call $R$ the \emph{Atiyah cocycle} of $\nabla$; if $R = 0$, we say that $\nabla$ is $A$-\emph{compatible}. The cohomology class $[R]$ is independent of the choice of $\nabla$, and we call it the \emph{Atiyah class} of the transverse complex manifold $(M,A)$. In the special case where $(M,A)$ is a complex manifold, i.e. when $A\cap\overline{A}=0$, this definition agrees with the classical Atiyah class.

\subsection{Leray--Cartan filtration}
\quad\par
\label{Subsection: Leray--Cartan filtration}
We continue to use the notation introduced in the previous subsection. Associated with $A$ is a decreasing filtration on the de Rham complex $(\Omega^*(M, \mathbb{C}), d)$:
$$ \dots \subset F_k \Omega^*(M, \mathbb{C}) \subset \dots \subset F_1 \Omega^*(M, \mathbb{C}) \subset F_0 \Omega^*(M, \mathbb{C}) = \Omega^*(M, \mathbb{C}),$$
known as the \emph{Leray--Cartan filtration}. The $k$-th filtered piece is defined by
\[
F_k\Omega^i(M, \mathbb{C}) = \{ \alpha \in \Omega^i(M, \mathbb{C}) \mid \iota_{a_1} \cdots \iota_{a_{i-k+1}} \alpha = 0, \ \forall a_1, \dots, a_{i-k+1} \in \Gamma(M, A) \}.
\]
Equivalently, an $i$-form belongs to $F_k\Omega^i(M, \mathbb{C})$ if and only if it is annihilated by the contraction with any $i-k+1$ vector fields tangent to $A$. Recall that the dual bundle $B^*$ may be canonically identified with the annihilator of $A$ in $T^*M \otimes \mathbb{C}$. Consequently, $\Omega_A^{*, 0}(M)$ admits a canonical embedding into $\Omega^*(M, \mathbb{C})$ and one has
$$F_k \Omega^*(M, \mathbb{C}) = \Omega_A^{k, 0}(M) \wedge \Omega^*(M, \mathbb{C}).$$
It follows that the associated graded complex satisfies $\text{Gr}_k \Omega^*(M, \mathbb{C}) = \Omega_A^{k, *}(M)$ and hence there is a short exact sequence of cochain complexes
$$ 0 \longrightarrow (F_{k+1}\Omega^*(M, \mathbb{C}), d) \longrightarrow (F_k \Omega^*(M, \mathbb{C}), d) \longrightarrow (\Omega_A^{k, *-k}(M), d_A^{\operatorname{Bott}}) \longrightarrow 0.$$
Here, by abuse of notations, $d_A^{\operatorname{Bott}}$ denotes the differential on $\Omega_A^{k, *}(M)$ induced by the Bott connection on $B$. For $k = 0$, the resulting long exact sequence in cohomology takes the form
\begin{align*}
	\cdots \longrightarrow H^1(F_1 \Omega^*(M, \mathbb{C}), d) \longrightarrow & H^1(M, \mathbb{C}) \longrightarrow H_A^1(M)\\
	\longrightarrow H^2(F_1 \Omega^*(M, \mathbb{C}), d) \longrightarrow & H^2(M, \mathbb{C}) \longrightarrow H_A^2(M) \longrightarrow \cdots
\end{align*}

\subsection{Formal geometry for transverse holomorphic structures}
\quad\par
\label{Subsection: transverse jet bundles}

We review the construction of transverse holomorphic jet bundles and transverse holomorphic (poly)differential operators associated with $(M, A)$. Although closely related ideas appear in \cite{Gom1980} in the study of the deformation theory of transverse holomorphic structures, these notions are not introduced there in an explicit form. Our presentation primarily follows \cite{WanYau2025}, where the constructions are developed in the more general setting of Nirenberg-integrable complex distributions. An alternative algebraic description of these objects can be found in \cite{BanStiXu2021}.

\begin{remark}
	In \cite{WanYau2025}, these objects are referred to simply as \emph{transverse jet bundles} and \emph{transverse differential operators}. We include the adjective \emph{holomorphic} throughout to emphasize that we are working in the setting of a transverse holomorphic structure.
\end{remark}

We begin by recalling that, for each $r \in \mathbb{N}$, there is a canonical short exact sequence:
\begin{equation}
	\label{Equation: short exact sequence of jet bundles}
	0 \longrightarrow \mathrm{Sym}^r (T^*M \otimes {\mathbb{C}}) \longrightarrow J^r M \longrightarrow J^{r-1}M \longrightarrow 0,
\end{equation}
where $J^rM$ denotes the $r$-th order jet bundle of the trivial complex line bundle $M \times \mathbb{C} \to M$ (We adopt the convention that $J^{-1}M$ is the zero vector bundle). The (infinite) jet bundle $JM$ of $M \times \mathbb{C} \to M$ is defined as the projective limit $\varprojlim_r J^rM$.

The transverse holomorphic jet bundle was originally introduced in \cite{WanYau2025} by a coordinate-free construction as a subbundle of $JM$. For our purpose, however, an equivalent local description is more convenient. Let $(x^1, ..., x^m, z^1, ..., z^{m'})$ be adapted coordinates on an open subset $U \subset M$. Denote by $y^i$, $u^j$, $\overline{u}^k$ the local sections of $\widehat{\operatorname{Sym}} (T^*M \otimes \mathbb{C})$ corresponding to $dx^i$, $dz^j$, $d\overline{z}^k$, respectively. For each $r \in \mathbb{N}$, there is a natural vector bundle isomorphism:
\begin{equation}
	\label{Equation: local trivialization of finite jet bundle}
	J^rM \vert_U \cong \operatorname{Sym}^{\leq r} (T^*M \otimes \mathbb{C}) \vert_U, \quad s \mapsto \sum_{\lvert I \rvert + \lvert J \rvert + \lvert K \rvert \leq r} \frac{1}{I! \cdot J! \cdot K!} \frac{\partial^{\lvert I \rvert + \lvert J \rvert + \lvert K \rvert} f}{\partial x^I \partial z^J \partial \overline{z}^K}(p) y^I u^J \overline{u}^K,
\end{equation}
where $\operatorname{Sym}^{\leq r} (T^*M \otimes \mathbb{C}) := \bigoplus_{i=0}^r \operatorname{Sym}^i (T^*M \otimes \mathbb{C})$, $f$ is a local function whose $r$th order jet at $p \in U$ is $s$ and $I \in \mathbb{N}^m$ and $J, K \in \mathbb{N}^{m'}$ range over multi-indices satisfying $\lvert I \rvert + \lvert J \rvert + \lvert K \rvert \leq r$. Passing to the projective limit gives the local identification
\begin{equation}
	\label{Equation: local trivialization of infinite jet bundle}
	JM \vert_U \cong \widehat{\operatorname{Sym}} (T^*M \otimes \mathbb{C}) \vert_U.
\end{equation}

Using \eqref{Equation: local trivialization of finite jet bundle}, we define a subbundle $J_A^rM$ of $J^rM$ by requiring that, under this local identification, $J_A^rM \vert_U \cong \operatorname{Sym}^{\leq r} (B^*) \vert_U$. The short exact sequence \eqref{Equation: short exact sequence of jet bundles} then restricts to:
\[
0 \longrightarrow \mathrm{Sym}^r (B^*) \longrightarrow J^r_A M \longrightarrow J^{r-1}_A M \longrightarrow 0.
\]
The \emph{transverse holomorphic jet bundle} $J_AM$ is defined as the projective limit $\varprojlim_r J_A^rM$.

The bundle $JM$ carries a canonical flat connection, the \emph{Grothendieck connection} $\nabla^{\operatorname{G}}$. A smooth section $\sigma \in \Gamma(M, JM)$ is $\nabla^{\operatorname{G}}$-flat if and only if $\sigma = j^\infty f$ for some $f \in C^\infty(M)$, where $j^\infty$ denotes the prolongation map. By \cite[Proposition 2.6]{WanYau2025}, $\nabla^{\operatorname{G}}$ preserves the subbundle $J_AM$. Under the local identification \eqref{Equation: local trivialization of infinite jet bundle}, it is given by
\begin{equation}
	\label{Equation: local formula of Grothendieck connection}
	\nabla^{\operatorname{G}} = dx^{i} \wedge \frac{\partial}{\partial x^i} + d\overline{z}^j \wedge \frac{\partial}{\partial \overline{z}^j} + dz^k \wedge \left( \frac{\partial}{\partial z^k} - \frac{\partial}{\partial u^k} \right).
\end{equation}
Here, the coordinate vector fields act on the coefficient functions $a_J$ in the expansion $s = \sum_J a_J u^J$ of a local smooth section $s$ of $J_AM$. If $f \in \mathcal{O}(M)$, then locally
\begin{equation}
	\label{Equation: local formula of prolongation}
	j^\infty f
	=
	\sum_I
	\frac{1}{I!} \frac{\partial^{\lvert I \rvert} f}{\partial z^I} u^I.
\end{equation}
Therefore, $j^\infty$ sends sections of $\mathcal{O}$ to sections of $J_AM$. Its image is characterized as follows.

\begin{proposition}
	The prolongation map $j^\infty$ restricts to an isomorphism from $\mathcal{O}$ to the sheaf of $\nabla^{\operatorname{G}}$-flat sections of $J_A M$.
\end{proposition}
\begin{proof}
	This is a consequence of \cite[Proposition 3.17]{WanYau2025}. Alternatively, one may give a direct proof using the local formulas \eqref{Equation: local formula of Grothendieck connection} and \eqref{Equation: local formula of prolongation}; we omit the details.
\end{proof}

We now introduce transverse holomorphic polydifferential operators, which play a central role in the construction of deformation quantizations. Although they can be defined intrinsically, for our purposes it is more convenient to work with their local description in adapted coordinates.

\begin{definition}
	A \emph{transverse holomorphic polydifferential operator} is a $\mathbb{C}$-multilinear morphism of sheaves
	\[
	D: \underbrace{\mathcal{O} \times \cdots \times \mathcal{O}}_{k \text{ factors}} \to \mathcal{O}
	\]
	for some $k > 0$ such that there exists $r \in \mathbb{N}$ for which $D$ is locally of the form
	\begin{equation}
		\label{Equation: local form of differential operator}
		D(f_1, \ldots, f_k) = \sum_{\lvert I_1 \rvert, \ldots, \lvert I_k \rvert \leq r} a_{I_1, \ldots, I_k} \frac{\partial^{\lvert I_1 \rvert} f_1}{\partial z^{I_1}} \cdots \frac{\partial^{\lvert I_k \rvert} f_k}{\partial z^{I_k}}
	\end{equation}
	in any system of adapted coordinates $(x^1, \ldots, x^m, z^1, \ldots, z^{m'})$, where the coefficients $a_{I_1, \ldots, I_k}$ are local sections of $\mathcal{O}$ independent of $f_1, \ldots, f_k$. When $k = 1$ (resp. $k = 2$), such a morphism is called a \emph{transverse holomorphic differential} (resp. \emph{bidifferential}) \emph{operator}.
\end{definition}

\begin{remark}
	In \cite{WanYau2025}, the bundle $D_A(M)$ of transverse holomorphic differential operators is constructed together with a natural flat $A$-connection $d_A^{D_A(M)}$. A transverse holomorphic differential operator is defined there as a $d_A^{D_A(M)}$-flat section of $D_A(M)$. By \cite[Proposition 2.9]{WanYau2025}, every such section canonically determines a sheaf morphism $\mathcal{O} \to \mathcal{O}$ and this correspondence is injective. The present definition is simply this equivalent sheaf-theoretic description, which admits an immediate extension to the polydifferential setting.
\end{remark}

For a transverse holomorphic polydifferential operator $D$, we define its \emph{prolongation} to be the $C^\infty(M)$-multilinear map
\[
j^\infty(D): \underbrace{\Gamma(M, J_AM) \times \cdots \times \Gamma(M, J_AM)}_{k \text{ factors}} \to \Gamma(M, J_AM),
\]
which, under the local identification \eqref{Equation: local trivialization of infinite jet bundle}, is given by
\[
j^\infty(D)(s_1, \ldots, s_k)
:=
\sum_{\lvert I_1 \rvert, \ldots, \lvert I_k \rvert \leq r} j^\infty(a_{I_1, \ldots, I_k}) \frac{\partial^{\lvert I_1 \rvert} s_1}{\partial u^{I_1}} \cdots \frac{\partial^{\lvert I_k \rvert} s_k}{\partial u^{I_k}},
\]
where $D$ is expressed locally as in \eqref{Equation: local form of differential operator}. The prolongation $j^\infty(D)$ is compatible with the prolongation of functions: for all open subset $U \subset M$ and all $f_1, \ldots, f_k \in \mathcal{O}(U)$,
\[
j^\infty (D(f_1, \ldots, f_k)) = j^\infty(D)(j^\infty f_1, \ldots, j^\infty f_k).
\]
Moreover, $j^\infty (D)$ is compatible with the Grothendieck connection in the sense that
\[
\nabla^{\operatorname{G}}\circ j^\infty(D)
=
j^\infty(D)\circ\nabla^{\operatorname{G}},
\]
where $\nabla^{\operatorname{G}}$ acts on the $k$-fold product
$\Gamma(M,J_AM)\times\cdots\times\Gamma(M,J_AM)$ via the Leibniz rule.

\section{Deformation quantization on transverse holomorphic symplectic manifolds}\label{Section3}

In Subsection~\ref{Subsection: C-presymplectic forms}, we recall the concept of a transverse holomorphic symplectic manifold $(M, \Omega)$, which serves as the fundamental geometric setting for our study, and prove the corresponding Darboux theorem. Subsection~\ref{Subsection: symplectic connections} introduces symplectic connections on $(M, \Omega)$. In Subsection~\ref{3.3}, we investigate star products on $\mathcal{O}[[\hbar]]$ defined by transverse holomorphic bidifferential operators. In particular, we examine the relationship between the existence of global star products and the Atiyah class. Finally, Subsection~\ref{3.1} formalizes the notion of a deformation quantization for $(M, \Omega)$.

\subsection{Transverse holomorphic symplectic manifolds}
\quad\par
\label{Subsection: C-presymplectic forms}
We now introduce the primary geometric setting of this paper.

\begin{definition}[\cite{KapOrl2003, qcr}]
	A \emph{transverse holomorphic symplectic form} on a smooth manifold $M$ is a closed complex $2$-form $\Omega \in \Omega^2(M, \mathbb{C})$ of constant rank, such that the kernel $A$ of the map 
	\begin{equation}
		\label{Equation: contraction of omega}
		T_{\mathbb C}M \to T^*M \otimes \mathbb{C}, \quad v \mapsto \iota_v \Omega,
	\end{equation}
	is a transverse holomorphic structure on $M$.
\end{definition}

In what follows, let $(M, \Omega)$ be a \emph{transverse holomorphic symplectic manifold}. This means that $\Omega$ is a transverse holomorphic symplectic form on $M$. Let $A$ denote the underlying transverse holomorphic structure. We call the real distribution $A_\mathbb{R} \subset TM$ for which $A_\mathbb{R} \otimes \mathbb{C} = A \cap \overline{A}$ the \emph{characteristic distribution} of $(M, \Omega)$.

Because $\Omega$ descends to a non-degenerate $2$-form $\underline{\Omega}$ on the quotient bundle $B := T_{\mathbb C}M / A$, the transverse complex dimension of $(M, A)$ is even; we denote it by $2n$.

\begin{remark}
	The triple $(T_{\mathbb C}M, A, \underline{\Omega})$ is an example of a \emph{symplectic Lie pair} in the sense of \cite[Definition 3.3]{VoglaireXu2015}.
\end{remark}

The nondegeneracy of $\underline{\Omega}$ induces a Poisson bracket on $\mathcal{O}(U)$, for every open subset $U \subset M$, by
\begin{equation}
	\label{Equation: Poisson bracket}
	\{f, g\}:=\underline{\Omega}^{-1}(df, dg), \quad f, g \in \mathcal{O}(U).
\end{equation}

The Darboux theorem in the holomorphic symplectic setting is well known (cf. \cite{GauJoo2015}), whose proof is essentially the holomorphic version of Moser's argument. We now state the Darboux theorem for transverse holomorphic symplectic manifolds.

\begin{theorem}
	\label{Theorem: Darboux}
	Let $(M, \Omega)$ be a transverse holomorphic symplectic manifold of transverse complex dimension $2n$. Then every point of $M$ admits a neighborhood $U$ equipped with adapted coordinates $(x^1, \dots, x^m, z^1, \dots, z^{2n})$ such that $\Omega|_U = \sum_{j=1}^n dz^j \wedge dz^{j+n}$, which we call \emph{transverse holomorphic Darboux coordinates}.
\end{theorem}
\begin{proof}
	By choosing adapted coordinates as in
	Subsection \ref{Subsection: Nirenberg integrable distributions},
	it suffices to consider the product case $M = U_1 \times U_2 \subset \mathbb{R}^m \times \mathbb{C}^{2n}$.
	The closedness of $\Omega$ and the identities
	$\iota_V\Omega=0$ for $V\in\Gamma(M, A)$ imply that
	$\Omega=\pi^*\Omega'$ for a holomorphic symplectic form
	$\Omega'$ on $U_2$, where $\pi$ is the projection onto the
	second factor. Applying the holomorphic Darboux theorem to
	$\Omega'$ and retaining the coordinates on $U_1$ yields the
	desired adapted coordinates, after shrinking the neighborhood
	if necessary.
\end{proof}

\subsection{Symplectic connections}
\quad\par
\label{Subsection: symplectic connections}
Recall that Fedosov's approach to deformation quantization relies on the choice of a symplectic connection, which is used to construct the Fedosov connection. In the setting of this paper, we likewise need to formulate the notion of a \emph{symplectic connection} for a transverse holomorphic symplectic manifold $(M, \Omega)$. At first sight, there are two natural candidates for this notion. The first is a symplectic connection on the associated symplectic Lie pair $(T_{\mathbb C}M, A, \underline{\Omega})$ in the sense of Voglaire--Xu \cite{VoglaireXu2015}. The second is a symplectic connection on the symplectic ringed space $(M, \mathcal{O}, \mathcal{B}, \Omega)$ in the sense of Vaisman \cite{Vai2002}.

These two notions are not equivalent. Indeed, the latter is substantially more restrictive, as its existence is equivalent to the vanishing of the Atiyah class of the Lie pair $(T_{\mathbb C}M, A)$. Since the former notion is sufficient for our purposes, we adopt it throughout this paper. In this section, we recall its basic properties; for further details, see \cite[Propositions 4.1--4.2]{VoglaireXu2015}.

We begin by recalling the notion of \emph{torsion} for a connection on the bundle $B = T_{\mathbb C}M / A$.

\begin{definition}
	Let $\nabla$ be a $TM\otimes\mathbb{C}$-connection on $B$.
	\begin{itemize}
		\item The \emph{torsion} of $\nabla$ is the smooth section $T^\nabla \in \Omega^2(M, B)$ given by
		\begin{equation*}
			T^\nabla(V, V') := \nabla_V \mathfrak{q}(V') - \nabla_{V'} \mathfrak{q}(V) - \mathfrak{q}([V, V']),
		\end{equation*}
		for all $V, V' \in \Gamma(M, TM\otimes\mathbb{C})$.
		\item The $TM\otimes\mathbb{C}$-connection $\nabla$ is said to be \emph{torsion-free} if $T^\nabla = 0$.
	\end{itemize}
\end{definition}

A torsion-free $TM\otimes\mathbb{C}$-connection on $B$ necessarily extends the Bott-connection $d_A^{\mathrm{Bott}}$; see, for instance, \cite[Lemma 1.4]{LiaStiXu2019}.

The following definition is adapted from the notion of symplectic connection introduced in \cite{VoglaireXu2015}.

\begin{definition}
	A \emph{symplectic connection} on the transverse holomorphic symplectic manifold $(M, \Omega)$ is a torsion-free $(TM\otimes\mathbb{C})$-connection $\nabla$ on $B$ such that $\nabla \Omega = 0$.
\end{definition}

We now recall the existence and classification of symplectic connections.

\begin{proposition}[\cite{VoglaireXu2015}, Proposition 4.2]
	\label{Proposition 4.4}
	There exists a symplectic connection on the transverse holomorphic symplectic manifold $(M, \Omega)$. Moreover, the space of symplectic connections on $(M, \Omega)$  is an affine space modelled by $\operatorname{Sym}^3 B^*$.
\end{proposition}

\subsection{Star product}\label{3.3}
\quad\par
A first attempt at defining deformation quantization of $(M, \Omega)$, closely paralleling the smooth symplectic setting, is the following notion.

\begin{definition}
	\label{Definition: star product}
	A \emph{star product} on $\mathcal{O}[[\hbar]]$ is a $\mathbb{C}[[\hbar]]$-bilinear sheaf morphism
	\begin{equation*}
		\star: \mathcal{O}[[\hbar]] \times \mathcal{O}[[\hbar]] \longrightarrow \mathcal{O}[[\hbar]]
	\end{equation*}
	satisfying the following conditions:
	\begin{enumerate}
		\item $(\mathcal{O}[[\hbar]], \star)$ is a sheaf of associative, unital $\mathbb{C}[[\hbar]]$-algebras with unit $1$;
		\item There exists a sequence $\{ C_r \}_{r=1}^\infty$ of transverse holomorphic bidifferential operators such that, for any open subset $U \subset M$ and $f, g \in \mathcal{O}(U)$,
		\begin{equation*}
			f \star g = f \cdot g + \sum_{r=1}^\infty \hbar^r C_r(f, g),
		\end{equation*}
		where $f \cdot g$ denotes the standard commutative pointwise product in $\mathcal{O}$.
	\end{enumerate}
\end{definition}

For a star product $\star$ on $\mathcal{O}[[\hbar]]$, we denote its commutator by $[ \,\cdot\,,\,\cdot\,]_\star$, i.e. $[f, g]_\star := f \star g - g \star f$.

\begin{remark}
	\label{Remark: Poisson bracket induced by star product}
	The associativity of $\star$ immediately implies that
	\[
	\{f,g\} := C_1(f,g) - C_1(g,f) \quad \text{for } f, g \in \mathcal{O}(U),
	\]
	defines a Poisson bracket on $\mathcal{O}(U)$. Thus, $\mathcal{O}$ acquires the structure of a sheaf of Poisson algebras.
\end{remark}

\begin{example}
	\label{Example: Moyal product}
	Denote by $(x^1, ..., x^m)$ the standard real coordinates on $\mathbb{R}^m$ and by $(z^1, ..., z^{2n})$ the standard complex coordinates on $\mathbb{C}^{2n}$. Let $M$ be an open subset of $\mathbb{R}^m \times \mathbb{C}^{2n}$ and $\Omega$ be a transverse holomorphic symplectic form on $M$ of the form $\Omega = \tfrac{1}{2} \Omega_{ij} dz^i \wedge dz^j$, where $(\Omega_{ij})$ is a $2n \times 2n$ skew-symmetric matrix whose entries $\Omega_{ij}$ are all constants. Let $(\Omega^{ij})$ be the inverse matrix of $(\Omega_{ij})$. In this case, there is a star product $\star_{\operatorname{M}}$ on $\mathcal{O}[[\hbar]]$ defined by
	\begin{equation}
		f \star_{\operatorname{M}} g = \sum_{r=0}^\infty \frac{1}{r!} \left( \frac{\hbar}{2} \right)^r \Omega^{i_1j_1} \cdots \Omega^{i_rj_r} \frac{\partial^r f}{\partial z^{i_1} \cdots \partial z^{i_r}} \frac{\partial^r g}{\partial z^{j_1} \cdots \partial z^{j_r}},
	\end{equation}
	for all open subset $U \subset M$ and all $f, g \in \mathcal{O}(U)[[\hbar]]$. We call this star product the \emph{Moyal product}.
\end{example}

\begin{definition}\label{Definition: gauge equivalence}
	Let $\star$ and $\star'$ be two star products on $\mathcal{O}[[\hbar]]$. A \emph{gauge equivalence} between $\star$ and $\star'$ is an isomorphism of sheaves of $\mathbb{C}[[\hbar]]$-algebras
	\[
	T : (\mathcal{O}[[\hbar]], \star) \to (\mathcal{O}[[\hbar]], \star')
	\]
	of the form
	\begin{equation*}
		T = \operatorname{Id} + \sum_{r=1}^{\infty} \hbar^r T_r,
	\end{equation*}
	where each $T_r: \mathcal{O} \to \mathcal{O}$ is a transverse holomorphic differential operator.
\end{definition}

\begin{proposition}
	\label{Proposition: gauge equivalent Poisson bracket}
	Suppose $\star$ and $\star'$ are gauge equivalent star products on $\mathcal{O}[[\hbar]]$. Then their associated Poisson brackets, as stated in Remark \ref{Remark: Poisson bracket induced by star product}, coincide.
\end{proposition}
\begin{proof}
	Let $\{ \,\cdot\,, \,\cdot\, \}$ and $\{ \,\cdot\,, \,\cdot\, \}'$ denote the Poisson brackets associated with $\star$ and $\star'$ respectively. Choose a gauge equivalence $T = \mathrm{id} + \hbar T_1 + \hbar^2 T_2 + \cdots$ between $\star$ and $\star'$. Then
	\begin{equation}
		T([f, g]_\star) = [T(f), T(g) ]_{\star'},
	\end{equation}
	for any open subset $U \subset M$ and $f, g \in \mathcal{O}(U)$. The left hand side takes the form
	\[
	T([f, g]_\star) = T\big(\hbar \{f, g\} + O(\hbar^2)\big) = \hbar \{f, g\} + O(\hbar^2),
	\]
	while the right hand side takes the form
	\begin{align*}
		[T(f), T(g) ]_{\star'} &= [f + \hbar T_1(f), g + \hbar T_1(g) ]_{\star'} + O(\hbar^2) \\
		&= [f, g ]_{\star'} + \hbar\big(T_1(f)g + fT_1(g) - T_1(g)f - gT_1(f)\big) + O(\hbar^2).
	\end{align*}
	The cross terms involving $T_1$ cancel out perfectly, leaving $[f, g ]_{\star'} = \hbar\{f, g\}' + O(\hbar^2)$. Equating the first-order terms in $\hbar$ from both sides yields $\{f, g\} = \{f, g\}'$.
\end{proof}

Although Definition \ref{Definition: star product} directly parallels
the notion of deformation quantization in smooth symplectic geometry,
the existence of a global star product on $\mathcal O[[\hbar]]$
is a restrictive requirement. By \cite[Theorem 2.1]{Vai2002},
the vanishing of the Atiyah class of $(M,\Omega)$ guarantees the
existence of such a star product. The following partial converse
is a transverse holomorphic analogue of
\cite[Proposition 3.1]{GutRaw2003}, where Gutt--Rawnsley show that
a star product satisfying order bounds on its first two coefficients
determines a symplectic connection. We include a coordinate-free proof in the present setting for completeness.

\begin{theorem}
	\label{Theorem: vanishing of Atiyah class}
	Let $(M, \Omega)$ be a transverse holomorphic symplectic manifold. Suppose there exists a star product $\star$ on $\mathcal{O}[[\hbar]]$ defined by a sequence of holomorphic bidifferential operators $\{C_r\}_{r=1}^\infty$ (as in Definition \ref{Definition: star product}), such that $C_1$ and $C_2$ are of order $1$ and $2$ in each argument, respectively. Then the Atiyah class of $(M, \Omega)$ vanishes.
\end{theorem}
\begin{proof}
	It suffices to construct a $\mathcal B$-connection in the sense
	of \cite[Section 1]{Vai2002}, namely a $\mathbb C$-linear sheaf morphism $D:\mathcal B \to 
	\mathcal B^*\otimes_{\mathcal O}\mathcal B$ satisfying $D(fb)=df\otimes b+fDb$ for sections $f$ of $\mathcal{O}$ and $b$ of $\mathcal{B}$. Since $\mathcal B$ is a locally free $\mathcal O$-module, such a
	$\mathcal B$-connection $D$ extends, just as in the holomorphic setting, by the smooth
	Leibniz rule to a $TM\otimes\mathbb C$-connection $\nabla$ on $B$ extending $d_A^{\operatorname{Bott}}$. The resulting connection $\nabla$ is $A$-compatible because $D$ preserves $d_A^{\operatorname{Bott}}$-flat sections.
	
	Transverse holomorphic Darboux coordinates provide local $d_A^{\operatorname{Bott}}$-closed
	frames of $B$ consisting of sections of the form $X_f := \{\,\cdot\, f\}$ for sections $f$ of $\mathcal{O}$. Thus it remains to prescribe $D_{X_f}X_g$ and verify
	compatibility with the connection Leibniz rules. For sections $f, g$ of $\mathcal{O}$ and for $r=1,2$, write
	\begin{equation*}
		C_r^+(f,g):=\tfrac{1}{2}(C_r(f,g)+C_r(g,f)).
	\end{equation*}
	The order bound on $C_1$ and unitality of $\star$ imply that $C_1^+$ defines a global $d_A^{\operatorname{Bott}}$-flat section of $\operatorname{Sym}^2 B$. Put $Y_f:=C_1^+(\,\cdot\,,f)$ and $C_{2,f}^+:=C_2^+(\,\cdot\,,f)$. For any local sections $f, g, h$ of $\mathcal{O}$, prescribe
	\begin{equation}
		\label{Equation: star product determines connection}
		(D_{X_f}X_g)(h) := X_fX_gh + \tfrac{1}{2} (Y_fY_g + Y_gY_f - Y_{C_1^+(f, g)})(h) + (d_{\operatorname{Hoch}}C_{2, h}^+)(f, g),
	\end{equation}
	where $(d_{\operatorname{Hoch}}C_{2, h}^+)(f, g) := fC_2^+(g, h) - C_2^+(fg, h) + C_2^+(f, h)g$.
	
	We first check that the right-hand side is a derivation in $h$. Symmetry of $C_2^+$ gives
	\begin{equation*}
		(d_{\operatorname{Hoch}} C_{2, h}^+)(f, g) = (d_{\operatorname{Hoch}} C_2^+)(f,g,h) + (d_{\operatorname{Hoch}} C_{2, f}^+)(g, h).
	\end{equation*}
	Associativity of $\star$ through order $\hbar^2$ implies
	\begin{equation*}
		(d_{\operatorname{Hoch}} C_2^+)(f,g,h)= X_hX_gf + Y_hY_gf - X_fX_gh - Y_fY_gh.
	\end{equation*}
	Since $X_hX_gf = -X_{\{f, g\}} h$ and $Y_hY_gf=Y_{C_1^+(f,g)}(h)$, substitution yields
	\[
	(D_{X_f}X_g)(h)
	=-X_{\{f, g\}}h
	+\tfrac{1}{2} ([Y_g,Y_f]+Y_{C_1^+(f,g)})(h)
	+(d_{\operatorname{Hoch}}C_{2,f}^+)(g,h).
	\]
	Every term in this expression is a
	derivation in $h$.
	
	For fixed $h$, the last two terms in
	\eqref{Equation: star product determines connection}
	are symmetric biderivations in $f,g$: this follows from
	the biderivation property of $C_1^+$ for the second term,
	and from the order bound on $C_2$ and unitality for the last.
	Nondegeneracy of $\Omega$ therefore defines an intrinsic section $\alpha(h)$ of $\mathcal{B}^* \otimes_\mathcal{O} \mathcal{B}^*$ satisfying $(D_{X_f}X_g)(h) = X_fX_gh + \alpha(h)(X_f, X_g)$. Setting $(D_bb')(h):=b(b'(h))+\alpha(h)(b,b')$ tensoriality gives the connection Leibniz rules, while the preceding verification on local Hamiltonian frames shows that this expression is a derivation in $h$. Thus $D$ is the required $\mathcal B$-connection.
\end{proof}

\subsection{Deformation quantization}\label{3.1}
\quad\par
To accommodate a broader class of transverse holomorphic symplectic manifolds, we instead adopt the approach of Nest--Tsygan \cite{NesTsy2001} and formulate deformation quantization in terms of a locally trivial sheaf of algebras rather than a globally trivial one. To this end, consider a pair
\begin{equation}
	\label{Equation: formal deformation}
	(\mathcal{A}_M^\hbar, \psi)
\end{equation}
consisting of a sheaf of flat, $\hbar$-adically complete $\mathbb{C}[[\hbar]]$-algebras $\mathcal{A}_M^\hbar$ on $M$ and a morphism of sheaves of $\mathbb{C}$-algebras $\psi: \mathcal{A}_M^\hbar \to \mathcal{O}$ such that the following forms a short exact sequence of $\mathbb{C}$-linear sheaves:
\begin{equation*}
	\begin{tikzcd}
		0 \ar[r] & \mathcal{A}_M^\hbar \ar[r, "\hbar (\,\cdot\,)"] & \mathcal{A}_M^\hbar \ar[r, "\psi"] & \mathcal{O} \ar[r] & 0,
	\end{tikzcd}
\end{equation*}
where $\hbar (\,\cdot\,)$ denotes multiplication by $\hbar$.

\begin{definition}
	A \emph{local differential trivialization} of $(\mathcal{A}_M^\hbar, \psi)$ is a triple $(U, \star, \tau)$, where $U \subset M$ is an open subset, $\star$ is a star product on $\mathcal{O} \vert_U[[\hbar]]$ (in the sense of Definition \ref{Definition: star product}) and
	\[
	\tau: \mathcal{A}_M^\hbar \vert_U \to (\mathcal{O} \vert_U[[\hbar]], \star)
	\]
	is an isomorphism of sheaves of $\mathbb{C}[[\hbar]]$-algebras such that for all open subsets $V \subset U$ and $f \in \mathcal{O}(V)$,
	\[
	(\psi \circ \tau^{-1})(f) = f.
	\]
\end{definition}

\begin{definition}
	A \emph{differential atlas} on $(\mathcal{A}_M^\hbar, \psi)$ is a collection $\{ (U_i, \star_i, \tau_i) \}_{i \in I}$ of local differential trivializations of $(\mathcal{A}_M^\hbar, \psi)$ such that $\{ U_i \}_{i \in I}$ forms an open cover of $M$ and for all $i, j \in I$ with $U_{ij} := U_i \cap U_j \neq \emptyset$, the composition
	\begin{equation*}
		\tau_j \circ \tau_i^{-1}: (\mathcal{O}|_{U_{ij}}[[\hbar]], \star_i \vert_{U_{ij}}) \to (\mathcal{O}|_{U_{ij}}[[\hbar]], \star_j \vert_{U_{ij}})
	\end{equation*}
	is a gauge equivalence between the two star products. A \emph{differential structure} on $(\mathcal{A}_M^\hbar, \psi)$ is a maximal differential atlas on it.
\end{definition}

Consequently, when $(\mathcal{A}_M^\hbar, \psi)$ admits a differential structure $\{ (U_i, \star_i, \tau_i) \}_{i \in I}$, for each $i \in I$ one obtains a Poisson bracket $\{\,\cdot\,,\,\cdot\,\}_i$ induced by $\star_i$ as stated in Remark \ref{Remark: Poisson bracket induced by star product}. By Proposition \ref{Proposition: gauge equivalent Poisson bracket}, these Poisson brackets $\{\,\cdot\,,\,\cdot\,\}_i$ glue to a global one on $\mathcal{O}$. In defining a deformation quantization of $(M, \Omega)$, we require that this recovers the Poisson bracket $\{\,\cdot\,,\,\cdot\,\}$ on $\mathcal{O}$ induced by the transverse holomorphic symplectic form $\Omega$ on $M$, which is given by \eqref{Equation: Poisson bracket}.

\begin{definition}\label{dqofo}
	A \emph{deformation quantization} of the transverse holomorphic symplectic manifold $(M, \Omega)$ is a pair $(\mathcal{A}_M^\hbar, \psi)$ of the form \eqref{Equation: formal deformation} together with a differential structure on it such that for all open subset $U \subset M$ and all $\widetilde{f}, \widetilde{g} \in \mathcal{A}_M^\hbar(U)$,
	\begin{equation}
		\label{Equation: semiclassical to Poisson}
		\{f, g\} = \psi \left( \frac{1}{\hbar} [\widetilde{f}, \widetilde{g}] \right),
	\end{equation}
	where $f = \psi(\widetilde{f})$, $g = \psi(\widetilde{g})$, and $[\widetilde{f}, \widetilde{g}] := \widetilde{f} \widetilde{g} - \widetilde{g} \widetilde{f}$ denotes the commutator in $\mathcal{A}_M^\hbar$.
\end{definition}

When referring to a deformation quantization of $(M, \Omega)$, we often simply write $(\mathcal{A}_M^\hbar, \psi)$, without explicitly expressing the differential structure on it. We say a local differential trivialization of $(\mathcal{A}_M^\hbar, \psi)$ \emph{adapted} if it is contained in the differential structure of $(\mathcal{A}_M^\hbar, \psi)$.

\begin{definition}
	\label{Definition: equivalence of deformation quantization}
	Two deformation quantizations $(\mathcal{A}_M^\hbar, \psi)$ and $((\mathcal{A}')_M^\hbar, \psi')$ of $(M, \Omega)$ are said to be \emph{equivalent} if there exists an isomorphism of sheaves of $\mathbb{C}[[\hbar]]$-algebras
	\begin{equation*}
		\Phi: (\mathcal{A}_M^\hbar, \psi) \to ((\mathcal{A}')_M^\hbar, \psi')
	\end{equation*}
	such that $\psi = \psi' \circ \Phi$, and that for all adapted local differential trivializations $(U, \star, \tau)$ of $(\mathcal{A}_M^\hbar, \psi)$ and $(U', \star', \tau')$ of $((\mathcal{A}')_M^\hbar, \psi')$ with $U \cap U' \neq \emptyset$,
	\begin{equation*}
		\tau' \circ \Phi \circ \tau^{-1}: (\mathcal{O} \vert_{U \cap U'}[[\hbar]], \star \vert_{U \cap U'}) \to (\mathcal{O} \vert_{U \cap U'}[[\hbar]], \star' \vert_{U \cap U'})
	\end{equation*}
	is a gauge equivalence between the star products $\star \vert_{U \cap U'}$ and $\star' \vert_{U \cap U'}$.
\end{definition}

We denote by $\operatorname{DQ}(M, \Omega)$ the set of equivalence classes of deformation quantizations of $(M, \Omega)$.

\begin{proposition}\label{localstar}
	Let $(\mathcal A_M^\hbar,\psi)$ be a deformation quantization of $(M, \Omega)$. Then for any point $p \in M$, there exists an adapted local differential trivialization $(U, \star, \tau)$ of $(\mathcal A_M^\hbar,\psi)$ such that $U$ is a neighbourhood of $p$ with transverse holomorphic Darboux coordinates as in Theorem \ref{Theorem: Darboux} and $\star$ is the Moyal product on $\mathcal{O} \vert_U[[\hbar]]$ with respect to these coordinates.
\end{proposition}
\begin{proof}
	Fix a point $p \in M$. By Theorem \ref{Theorem: Darboux}, we can  choose transverse holomorphic Darboux coordinates $(x^1, \dots, x^m, z^1, \dots, z^{2n})$ on a neighbourhood $U \cong U_1 \times U_2 \subset \mathbb{R}^m \times \mathbb{C}^{2n}$ as in Subsection \ref{Subsection: Nirenberg integrable distributions}, together with an adapted local differential trivialization $(U, \star', \tau')$ of $(\mathcal{A}_M^\hbar, \psi)$.
	
	Observe that the pullback along the natural projection $\pi: U \to U_2$ defines an isomorphism of sheaves of $\mathbb{C}[[\hbar]]$-algebras
	\[
	\pi^*: (\mathcal{O}_{U_2}[[\hbar]], \widetilde{\star}) \to (\mathcal{O} \vert_U[[\hbar]], \star),
	\]
	where $\mathcal{O}_{U_2}$ denotes the sheaf of holomorphic functions on $U_2$, $\widetilde{\star}$ denotes the holomorphic Moyal product on $\mathcal{O}_{U_2}[[\hbar]]$ associated with the holomorphic symplectic form $\sum_{i=1}^{n} dz^i \wedge dz^{i+n}$ on $U_2$, and $\star$ denotes the Moyal product on $\mathcal{O} \vert_U[[\hbar]]$ with respect to the above transverse holomorphic Darboux coordinates. On the other hand, transporting $\star'$ through $\pi^*$ yields a (holomorphic) star product $\widetilde{\star}'$ on $\mathcal O_{U_2}[[\hbar]]$.
	
	The local classification of star products on holomorphic symplectic manifolds \cite[Theorem~4.1]{NesTsy2001} guarantees the existence of a (holomorphic) gauge equivalence $\widetilde{\gamma}: (\mathcal{O}_{U_2} [[\hbar]], \widetilde{\star}') \to (\mathcal{O}_{U_2}[[\hbar]], \widetilde{\star})$. Conjugating by $\pi^*$, we obtain a gauge equivalence in the sense of Definition \ref{Definition: gauge equivalence}:
	$$\gamma := \pi^* \circ \widetilde{\gamma} \circ (\pi^*)^{-1}: (\mathcal{O} \vert_U [[\hbar]], \star') \to (\mathcal{O} \vert_U [[\hbar]], \star).$$
	Consequently, $(U, \star, \gamma \circ \tau')$ forms an adapted local differential trivialization of $(\mathcal{A}_M^\hbar, \psi)$.
\end{proof}

\section{Main results: existence and classification}\label{section4}
This section presents the main results of this paper. Let $(M, \Omega)$ be a transverse holomorphic symplectic manifold, and we follow the relevant notations in the previous section.

We begin in Section \ref{4.1} by introducing the Weyl bundle $W$ for $(M, \Omega)$ and gauge transformations of $W$. Section \ref{4.2} develops the basic algebraic structures on the space of $W$-valued differential forms. In Section \ref{4.3}, we introduce Fedosov connections. In Section \ref{4.4}, we reduce the classification of deformation quantizations of $(M, \Omega)$ to the classification of Fedosov connections on $W$. Finally, in Section \ref{4.5} we use the non-commutative period map to prove our classification theorem for deformation quantizations of $(M, \Omega)$ under the admissibility assumption.

\subsection{The Weyl bundle}\label{4.1}
\quad\par
In this section, we construct an algebra bundle $(W, *)$, called the \emph{Weyl bundle}, associated with the transverse holomorphic symplectic manifold $(M, \Omega)$. Its underlying (infinite rank) vector bundle is defined by
\begin{equation*}
	W := \widehat{\operatorname{Sym}}(B^*)[[\hbar]],
\end{equation*}
where $\widehat{\operatorname{Sym}}(B^*)$ is the bundle of completed symmetric algebras on $B^*$. Recall that we write local foliated coordinates on $M$ generally as $(x^1, ..., x^m, z^1, ..., z^{2n})$. In particular, $(dz^1, ..., dz^{2n})$ is a local $d_A^{\operatorname{Bott}}$-closed frame of $B^*$. To distinguish elements of $\widehat{\operatorname{Sym}}(B^*)$ from differential forms, we write $u^i$ in place of the local $1$-form $dz^i$ when it is regarded as a local section of $\widehat{\operatorname{Sym}}(B^*)$. A local smooth section of $W$ is a formal power series of the form
\begin{equation*}
	\sum_{r, l \geq 0} \hbar^r a_{r, i_1, \dots, i_l} u^{i_1} \cdots u^{i_l},
\end{equation*}
where $a_{r, i_1, \dots, i_l}$ are local smooth complex-valued functions on $M$. Here, we suppress the notation for symmetric products in the expression $u^{i_1} \cdots u^{i_l}$.

We equip $W$ with the structure of an algebra bundle via a $\mathcal{C}^\infty(M, \mathbb{C})[[\hbar]]$-bilinear map
\begin{equation*}
	*: \Gamma(M, W) \times \Gamma(M, W) \to \Gamma(M, W),
\end{equation*}
called the \emph{fibrewise Moyal product}. It is defined by
\begin{equation}
	\label{Equation: fibrewise Moyal product}
	f * g = \sum_{r=0}^\infty \frac{1}{r!} \left( \frac{\hbar}{2} \right)^r \Omega^{i_1j_1} \cdots \Omega^{i_rj_r} \frac{\partial^r f}{\partial u^{i_1} \cdots \partial u^{i_r}} \frac{\partial^r g}{\partial u^{j_1} \cdots \partial u^{j_r}}.
\end{equation}
Here, writing $\Omega = \frac{1}{2} \Omega_{ij} dz^i \wedge dz^j$, we denote by $(\Omega^{ij})$ the inverse matrix of $(\Omega_{ij})$.

We assign weights on $W$ by setting the weight of $u^i$ to be $1$ and the weight of $\hbar$ to be $2$. For $r \in \mathbb{N}$, we denote by $W_{(r)}$ the subbundle of $W$ of weight at least $r$. This defines a decreasing filtration of $W$, which we call the \emph{weight filtration}:
\begin{equation*}
	W_{(0)} \supset W_{(1)} \supset \cdots \supset W_{(r)} \supset \cdots.
\end{equation*}

The weight filtration on the Weyl bundle naturally induces a filtration on the sheaf of derivations that are continuous with respect to the filtration topology. More precisely, for each $r \in \mathbb{N}$, let $\operatorname{Der}_{\geq r}$ denote the sheaf of $\mathcal{C}_M^\infty$-modules defined by assigning to every open subset $U \subset M$ the $\mathcal{C}^\infty(U)$-module consisting of all $\mathcal{C}^\infty(U)$-linear maps
$$D: \Gamma(U, W) \to \Gamma(U, W)$$
satisfying:
\begin{itemize}
	\item for all $a, b \in \Gamma(U, W)$, $D(a * b) = (Da) * b + a * (Db)$.
	\item for all $r' \in \mathbb{N}$ and $a \in \Gamma(U, W_{(r')})$, $Da \in \Gamma(U, W_{(r'+r)})$.
\end{itemize}
These sheaves form a decreasing filtration: for each open subset $U \subset M$,
$$\operatorname{Der}_{\geq 0}(U) \supset \operatorname{Der}_{\geq 1}(U) \supset \cdots \supset \operatorname{Der}_{\geq r}(U) \supset \cdots.$$
For each $r \in \mathbb{N}$, there is a short exact sequence
\begin{center}
	\begin{tikzcd}
		0 \ar[r] & \bigoplus_{l = \lceil \frac{r}{2}+1 \rceil}^\infty \hbar^l \cdot \mathcal{C}^\infty(U) \ar[r, "\iota"] & \Gamma(U, W_{(r+2)}) \ar[r, "\operatorname{ad}"] & \operatorname{Der}_{\geq r}(U) \ar[r] & 0
	\end{tikzcd}
\end{center}
where $\iota$ is the inclusion map, $\operatorname{ad}$ is given by $s \mapsto \tfrac{1}{\hbar} [s, \,\cdot\,]_*$ and $[\,\cdot\,,\,\cdot\,]_*$ denotes the commutator of $*$. On the other hand, the natural short exact sequence
\begin{center}
	\begin{tikzcd}
		0 \ar[r] & \Gamma(U, W_{(r+3)}) \ar[r] & \Gamma(U, W_{(r+2)}) \ar[r] & \Gamma(U, \operatorname{Gr}_{r+2} W) \ar[r] & 0
	\end{tikzcd}
\end{center}
induces a short exact sequence
\begin{center}
	\begin{tikzcd}
		0 \ar[r] & \operatorname{Der}_{\geq (r+1)}(U) \ar[r] & \operatorname{Der}_{\geq r}(U) \ar[r] & \mathcal{A}^{(r)}(U) \ar[r] & 0
	\end{tikzcd}
\end{center}
where $\operatorname{Der}_{\geq (r+1)}(U) \to \operatorname{Der}_{\geq r}(U)$ is the inclusion map,
\begin{equation}
	\label{Equation: graded piece of derivations}
	\mathcal{A}^{(r)}(U) := \bigoplus_{l = 0}^{\lceil \frac{r}{2}+1 \rceil - 1} \hbar^l \cdot \Gamma(U, \operatorname{Sym}^{r+2-2l} (B^*)),
\end{equation}
and the morphism $\operatorname{Der}_{\geq r}(U) \to \mathcal{A}^{(r)}(U)$ assigns $D \in \operatorname{Der}_{\geq r}(U)$ to the unique element $a \in \mathcal{A}^{(r)}(U)$ satisfying $D - \tfrac{1}{\hbar}[a, \,\cdot\,]_* \in \operatorname{Der}_{\geq (r+1)}(U)$.

\begin{definition}
	A \emph{gauge transformation} of the Weyl bundle $W$ is a $\mathcal{C}^\infty(M)$-algebra automorphism
	\[
	g: (\Gamma(M, W), *) \to (\Gamma(M, W), *)
	\]
	such that for all $r \in \mathbb{N}$ and $a \in \Gamma(M, W_{(r)})$, $g(a) - a \in \Gamma(M, W_{(r+1)})$.
\end{definition}

Gauge transformations of the Weyl bundle forms a sheaf of groups, and we denote it by $\mathcal{G}$.

Let $r \geq 1$ and $U \subset M$ be open. For every $D \in \operatorname{Der}_{\geq r}(U)$, since it strictly increases the weight of a section $a \in \Gamma(U, W)$, the series
\begin{equation*}
	\exp \left( D \right)(a) := \sum_{r=0}^\infty \frac{1}{r!} D^r (a)
\end{equation*}
converges in the filtration topology of $\Gamma(U, W)$. Thus $\exp \left( D \right)$ defines an automorphism of the $\mathcal{C}^\infty(U)$-module $\Gamma(U, W)$. Moreover, if $r' \in \mathbb{N}$ and $a \in \Gamma(M, W_{(r')})$, then
$$\exp(D) a - a \in \Gamma(M, W_{(r'+r)}) \subset \Gamma(M, W_{(r'+1)}),$$
so $\exp(D)$ is a gauge transformation of $W \vert_U$.

If $D, D' \in \operatorname{Der}_{\geq r}(U)$, then the Baker--Campbell--Hausdorff formula gives
\[
\exp(D) \exp(D') = \exp(D + D' + R(D, D')),
\]
for some $R(D, D') \in \operatorname{Der}_{\geq (r+1)}(U)$. It then follows from $\exp(0) = id$ and $(\exp(D))^{-1} = \exp(-D)$ that the subset $\mathcal{G}^{(r)}(U) := \{\exp(D): D \in \mathrm{Der}_{\ge r}(U)\}$ forms a subgroup of $\mathcal{G}(U)$. These subgroups define a sheaf of groups $\mathcal{G}^{(r)}$.

Suppose $g \in \mathcal{G}(U)$. Since $g - \mathrm{Id}$ raises the weight by at least one, the series
\[
\log g := \sum_{r=1}^\infty (-1)^{r-1} \frac{(g - \mathrm{id})^r}{r}
\]
converges in the filtration topology and $\log g \in \operatorname{Der}_{\geq 1}(U)$. Moreover, $\exp(\log g) = g$. Consequently, $\mathcal{G}$ admits a descending filtration:
\begin{equation*}
	\mathcal{G} = \mathcal{G}^{(1)} \supset \mathcal{G}^{(2)} \supset \dots \supset \mathcal{G}^{(r)} \supset \dots
\end{equation*}

\subsection{Basic structures on Weyl bundle valued differential forms}\label{4.2}
\quad\par

We extend $*$ to a $\mathcal{C}^\infty(M)[[\hbar]]$-bilinear product
\begin{equation*}
	\Omega^*(M, W) \times \Omega^*(M, W) \to \Omega^*(M, W)
\end{equation*}
by combining the fibrewise Moyal product on $W$ with the wedge product of differential forms on $M$. We also endow $\Omega^*(M, W)$ with the decreasing filtration
\begin{equation}
	\label{Equation: Weight filration on Weyl valued forms}
	\operatorname{F}_k \Omega^*(M, W) = \sum_{i=0}^k \Omega_A^{i, 0}(M) \wedge \Omega^{*-i}(M, W_{(k-i)}),
\end{equation}
where the sum is taken in the ordinary sense rather than as a direct sum. This filtration combines the Leray--Cartan filtration on $\Omega^*(M, \mathbb{C})$ with the weight filtration on $W$. 
Equipped with this filtration, $(\Omega^*(M, W), *)$ forms a filtered algebra.

The product $*$ therefore induces a product on the associated graded algebra of the filtration \eqref{Equation: Weight filration on Weyl valued forms}, which we denote by the same symbol by a slight abuse of notation. Since $W := \widehat{\operatorname{Sym}}(B^*)[[\hbar]]$ is naturally graded. Hence, the underlying vector space of the associated graded algebra is $\Omega_A^{*, *}(M, W)$. More explicitly,
\[
\operatorname{Gr}_k \Omega^i(M, W) = \bigoplus_{j=0}^{\min \{i, k\}} \Omega_A^{j, i-j}(M, \operatorname{Gr}_{k-j} W), \quad \text{for } i, k \geq 0.
\]
The Bott connection on $B$ induces a flat $A$-connection on $\bigwedge B^* \otimes W$, which we continue to denote by $d_A^{\operatorname{Bott}}$. Since the coefficients $\omega^{ij}$ are local sections of $\mathcal{O}$, the operator $d_A^{\operatorname{Bott}}$ is a graded derivation of $*$. More precisely, for any homogeneous $\alpha, \beta \in \Omega_A^{*, *}(M, W)$,
\begin{equation*}
	d_A^{\operatorname{Bott}} (\alpha * \beta) = (d_A^{\operatorname{Bott}} \alpha) * \beta + (-1)^{\lvert \alpha \rvert} \alpha * (d_A^{\operatorname{Bott}} \beta),
\end{equation*}
where $\lvert \alpha \rvert$ is the total degree of $\alpha$.

Define a $\mathcal{C}^\infty(M)[[\hbar]]$-linear operator
\begin{equation*}
	\delta: \Omega^*(M, W) \to \Omega^*(M, W)
\end{equation*}
by $\delta := dz^i \wedge \tfrac{\partial}{\partial u^i}$, and define a $\mathcal{C}^\infty(M)[[\hbar]]$-linear operator
\begin{equation}
	\label{Equatoin: pi zero star}
	\pi_{0, *}: \Omega^*(M, W) \to \Omega_A^{0, *}(M)[[\hbar]]
\end{equation}
as the composition of the canonical projection $\Omega^*(M, W) \to \Omega^*(M, \mathbb{C})[[\hbar]]$ with the restriction map $\Omega^*(M, \mathbb{C})[[\hbar]] \to \Omega_A^{0, *}(M)[[\hbar]]$. Both operators descend naturally to $\mathcal{C}^\infty(M)[[\hbar]]$-linear maps $\Omega_A^{*, *}(M, W) \to \Omega_A^{*, *}(M, W)$ and $\Omega_A^{*, *}(M, W) \to \Omega_A^{0, *}(M)[[\hbar]]$, which, by abuse of notation, will again be denoted by $\delta$ and $\pi_{0, *}$, respectively.

Observe that $(\Omega_A^{*, *}(M, W), d_A^{\operatorname{Bott}}, \delta)$ is a bicomplex. We shall require the following lemma concerning the associated total complex $( \operatorname{tot} ( \Omega_A^{*, *}(M, W) ), d_A^{\operatorname{Bott}} - \delta )$.

\begin{lemma}
	\label{Lemma 3.19}
	The inclusion map $\iota: \Omega_A^{0, *}(M)[[\hbar]] \to \Omega_A^{*, *}(M, W)$ together with the map $\pi_{0, *}$ form a deformation retract 
	\begin{equation*}
		(\Omega_A^{0, *}(M)[[\hbar]], d_A) \hookrightarrow ( \operatorname{tot} ( \Omega_A^{*, *}(M, W) ), d_A^{\operatorname{Bott}} - \delta ).
	\end{equation*}
\end{lemma}

We begin with constructing a $\mathcal{C}^\infty(M)[[\hbar]]$-linear map
\begin{equation*}
	\delta^{-1}: \Omega_A^{*, *}(M, W) \to \Omega_A^{*, *}(M, W),
\end{equation*}
which will serve as a homotopy operator.
It is defined as follows. For $a \in \Omega_A^{p, *}(M, \operatorname{Sym}^l B^*)$,
\begin{equation}
	\label{equation: homotopy operator}
	\delta^{-1} a := \frac{1}{l+p} u^i \iota_{s^i} a \quad \text{for } l + p > 0,
\end{equation}
where $s^i = \mathfrak{q}\left(\tfrac{\partial}{\partial z^i} \right)$, and $\delta^{-1} a := 0$ for $l + p = 0$.

\begin{remark}
	Via a splitting $T_{\mathbb C}M \cong A \oplus B$ of smooth vector bundles, which induces an isomorphism $\Omega^*(M, W) \cong \Omega_A^{*, *}(M, W)$, we can identify $\delta^{-1}$ as a map $\Omega^*(M, W) \to \Omega^*(M, W)$.
\end{remark}

\begin{proof}[\myproof{Lemma}{\ref{Lemma 3.19}}]
	Clearly, $\delta \circ \iota = 0$, $\pi_{0, *} \circ \delta = 0$, and $\pi_{0, *} \circ \iota = \operatorname{Id}$ on $\Omega_A^{0, *}(M)[[\hbar]]$. Together with the fact that $d_A^{\operatorname{Bott}}$ coincides with the leafwise de Rham differential $d_A$ when restricted on $\Omega_A^{0, *}(M)[[\hbar]]$, we see that
	\begin{align*}
		\iota: & \,(\Omega_A^{0, *}(M)[[\hbar]], d_A) \to ( \operatorname{tot} ( \Omega_A^{*, *}(M, W) ), d_A^{\operatorname{Bott}} - \delta ),\\
		\pi_{0, *}: & \,( \operatorname{tot} ( \Omega_A^{*, *}(M, W) ), d_A^{\operatorname{Bott}} - \delta ) \to (\Omega_A^{0, *}(M)[[\hbar]], d_A)
	\end{align*}
	are cochain maps. By \cite[Section 1.3]{ChaLia2024}, the equality
	\begin{equation}
		\label{Equation: formal Hodge decomposition}
		\operatorname{Id} - \iota \circ \pi_{0, *} = \delta \circ \delta^{-1} + \delta^{-1} \circ \delta
	\end{equation}
	holds on $\Omega_A^{*, *}(M, W)$. It follows easily from \eqref{equation: homotopy operator} that $d_A^{\operatorname{Bott}} \circ \delta^{-1} + \delta^{-1} \circ d_A^{\operatorname{Bott}} = 0$. Hence, the equality
	\begin{equation}
		\label{Equation: refined formal Hodge decomposition}
		\operatorname{Id} - \iota \circ \pi_{0, *} = (d_A^{\operatorname{Bott}} - \delta) \circ (-\delta^{-1}) + (-\delta^{-1}) \circ (d_A^{\operatorname{Bott}} - \delta)
	\end{equation}
	holds on $\Omega_A^{*, *}(M, W)$, concluding the proof.
\end{proof}

\begin{remark}
	\label{Remark: application of lemma}
	We will apply Lemma \ref{Lemma 3.19} in the following situation. Let $k \geq 2$, $i \geq 1$ and let 
	\begin{equation*}
		\alpha \in \operatorname{Gr}_k \Omega^i(M, W)
	\end{equation*}
	be $(d_A^{\operatorname{Bott}} - \delta)$-closed. It then follows from \eqref{Equation: refined formal Hodge decomposition} that
	\begin{equation*}
		\alpha + (d_A^{\operatorname{Bott}} - \delta)\beta = \tau,
	\end{equation*}
	where $\beta := \delta^{-1} \alpha \in \operatorname{Gr}_k \Omega^{i-1}(M, W)$ and $\tau := \pi_{0, *}(\alpha) \in \Omega_A^{0, i}(M)[[\hbar]]$ is of weight $k$. By the cochain map property of $\pi_{0, *}$ from Lemma \ref{Lemma 3.19}, we see that $\tau$ is $d_A$-closed.
\end{remark}

\subsection{Fedosov connections}\label{4.3}
\quad\par
Consider a symplectic connection $\nabla$ on $(M, \Omega)$. It induces a $T_{\mathbb C}M$-connection $\nabla^W$ on $W = \widehat{\operatorname{Sym}} B^* [[\hbar]]$. 
Since $\nabla \omega = 0$, $\nabla^W$ defines a $*$-graded derivation; namely,
\begin{equation*}
	\nabla^W(\alpha * \beta) = (\nabla^W \alpha) * \beta + (-1)^{\lvert \alpha \rvert} \alpha * (\nabla^W \beta).
\end{equation*}
Moreover, there exists a unique element $R^\nabla \in \Omega^2(M, \operatorname{Sym}^2 B^*)$ such that $\pi_{0, *}(R^\nabla) = 0$ and $(\nabla^W)^2 = \tfrac{1}{\hbar}[R^\nabla, \,\cdot\,]_*$, where $[\,\cdot\,, \,\cdot\,]_*$ is the graded commutator of $*$. Consequently, 
\begin{equation*}
	(\Omega^*(M, W), R^\nabla, \nabla^W, \tfrac{1}{\hbar}[\,\cdot\,, \,\cdot\,]_*)
\end{equation*}
forms a curved differential graded Lie algebra. Given an element $\gamma\in\Omega^1(M,W)$, we may twist this structure to obtain a new curved differential graded Lie algebra
\begin{equation*}
	(\Omega^*(M, W), R^{\nabla, \gamma}, D^{\nabla, \gamma}, \tfrac{1}{\hbar}[\,\cdot\,, \,\cdot\,]_*),
\end{equation*}
where
\begin{equation*}
	D^{\nabla, \gamma} := \nabla^W + \tfrac{1}{\hbar}[\gamma, \,\cdot\,]_* \quad \text{and} \quad R^{\nabla, \gamma} := R^\nabla + \nabla^W \gamma + \tfrac{1}{2\hbar}[\gamma, \gamma]_*.
\end{equation*}
Note that there is a unique element $\tilde{\Omega} \in \Omega_A^{1, 0}(M, B^*) \subset \Omega^1(M, B^*)$ such that $\delta = \tfrac{1}{\hbar}[\tilde{\Omega}, \,\cdot\,]_*$. Indeed, we are only interested in $W$-valued $1$-forms $\gamma$ lying in the affine space $-\tilde{\Omega} + \operatorname{F}_3 \Omega^1(M, W)$.

\begin{remark}
	\label{Remark: choice of symplectic connections}
	Let $\nabla'$ be another symplectic connection on $(M, \Omega)$ and let $\gamma' \in \Omega^1(M, W)$. The corresponding twisted connection $D^{\nabla',\gamma'}$ can always be expressed in terms of the given symplectic connection $\nabla$ and a suitable element $\gamma\in\Omega^1(M,W)$.
	
	Since both $\nabla$ and $\nabla'$ extend the Bott connection, their difference $\widetilde{\Gamma} := \nabla' - \nabla$ defines an element of $\Omega_A^{1, 0}(M, \operatorname{End}(B))$. Because both connections preserve $\Omega$, the identification $B\cong B^*$ induced by $\Omega$ allows us to view $\widetilde{\Gamma}$ as an element
	$$\Gamma \in \Omega_A^{1, 0}(M, \operatorname{Sym}^2 B^*) \in \operatorname{F}_3\Omega^1(M, W).$$
	It follows that $(\nabla')^W = \nabla^W + \tfrac{1}{\hbar}[\Gamma, \,\cdot\,]_*$. Setting $\gamma := \gamma' + \Gamma$, we obtain
	\begin{equation*}
		D^{\nabla', \gamma'} = D^{\nabla, \gamma} \quad \text{and} \quad R^{\nabla', \gamma'} = R^{\nabla, \gamma}.
	\end{equation*}
	In particular, if $\gamma' \in -\tilde{\Omega} + \operatorname{F}_3 \Omega^1(M, W)$, then $\gamma \in -\tilde{\Omega} + \operatorname{F}_3 \Omega^1(M, W)$ as well.
\end{remark}

In view of the preceding remark, we henceforth fix a symplectic connection $\nabla$ on $(M, \Omega)$. Accordingly, whenever no confusion can arise, we suppress the superscript $\nabla$ and write $D^\gamma$, $R^\gamma$, etc., in place of $D^{\nabla,\gamma}$, $R^{\nabla,\gamma}$, and so on.

\begin{definition}
	\label{Definition: Fedosov connection}
	A \emph{Fedosov connection} on $W$ is a flat $T_{\mathbb C}M$-connection on $W$ of the form
	\begin{equation*}
		D^\gamma = \nabla^W + \tfrac{1}{\hbar}[\gamma, \, \cdot\,]_*
	\end{equation*}
	for some $\gamma \in -\tilde{\Omega} + \operatorname{F}_3 \Omega^1(M, W)$.
\end{definition}

\begin{remark}
	The flatness of $D^\gamma$ means $(D^\gamma)^2 = 0$. Since $(D^\gamma)^2 = \frac{1}{\hbar}[R^\gamma, \,\cdot\,]_*$, this is equivalent to requiring the curvature form $R^\gamma$ to take values in the center of $[\,\cdot\,,\,\cdot\,]_*$, i.e. $R^\gamma \in \Omega^2(M, \mathbb{C})[[\hbar]]$. 
\end{remark}

The action of gauge transformations $g$ of $W$ on $\Gamma(M, W)$ induces, by conjugation, an action on the space of Fedosov connections $D^\gamma$ on $W$, which we denote by $g \cdot D^\gamma$. We say that two Fedosov connections $D^\gamma$ and $D^{\gamma'}$ are \emph{gauge equivalent} if there exists a gauge transformation $g$ of $W$ such that $D^{\gamma'} = g \cdot D^\gamma$. We briefly recall some basic properties of this action. For each $s \in \Gamma(M, W_{(3)})$ and each $\gamma \in -\tilde{\Omega} + \operatorname{F}_3 \Omega^1(M, W)$, we define
\begin{equation}
	\label{Equation: lifted gauge transform}
	s \bullet \gamma := \gamma - \frac{\exp(\frac{1}{\hbar} [s,\,\cdot\,]_*) - \operatorname{Id}}{\frac{1}{\hbar} [s,\,\cdot\,]_*} (D^\gamma s).
\end{equation}
Here the quotient is understood as the series
\begin{equation*}
	\frac{\exp(\frac{1}{\hbar} [s,\,\cdot\,]_*) - \operatorname{Id}}{\frac{1}{\hbar} [s,\,\cdot\,]_*} = \sum_{r=0}^\infty \frac{1}{(r+1)!} \left( \frac{1}{\hbar} \right)^r \underbrace{[s, \dots, [s, \,\cdot\,]_* \cdots ]_*}_{r\text{-fold iterated commutator}},
\end{equation*}
which is a well defined operator on $\Gamma(M, W)$ because $\tfrac{1}{\hbar}[s, \,\cdot\,]_* \in \operatorname{Der}_{\geq 1}(M)$. Also, since $D^\gamma$ preserves the filtration on $\Omega^*(M, W)$, we have
\[
s \bullet \gamma \in -\tilde{\Omega} + \operatorname{F}_3 \Omega^1(M, W).
\]

Recall that any gauge transformation of $W$ is of the form $\exp(\tfrac{1}{\hbar}[s,\,\cdot\,]_*)$, where $s \in \Gamma(M, W_{(3)})$. The following proposition is standard, but we include a proof for completeness.

\begin{proposition}
	\label{Proposition: gauge transformation of Fedosov connection}
	For all $s \in \Gamma(M, W_{(3)})$ and all $\gamma \in -\tilde{\Omega} + \operatorname{F}_3 \Omega^1(M, W)$,
	\[
	D^{s \bullet \gamma} = g \cdot D^\gamma \quad \text{and} \quad R^{s \bullet \gamma} = g(R^\gamma),
	\]
	where $g = \exp(\tfrac{1}{\hbar}[s, \,\cdot\,]_*)$. In particular, if $D^\gamma$ is a Fedosov connection on $W$, then $g \cdot D^\gamma$ is also a Fedosov connection on $W$ and $R^{s \bullet \gamma} = R^\gamma$.
\end{proposition}
\begin{proof}
	By definition, $g \cdot D^\gamma = g \circ D^\gamma \circ g^{-1}$. Let $X = \tfrac{1}{\hbar}[s, \,\cdot\,]_*$. We have the standard identity:
	\[
	[D^\gamma, g^{-1}] = -g^{-1} \circ \frac{e^{\operatorname{ad}_X} - \operatorname{Id}}{\operatorname{ad}_X} ([D^\gamma, X]).
	\]
	Here, $[\,\cdot\,,\,\cdot\,]$ denotes the commutator of operators and $\operatorname{ad}_X = [X, \,\cdot\,]$ denotes the adjoint action of $X$. Since $D^\gamma$ and $\tfrac{1}{\hbar}[s, \,\cdot\,]_*$ are $*$-derivations,
	\[
	\frac{e^{\operatorname{ad}_X} - \operatorname{Id}}{\operatorname{ad}_X} ([D^\gamma, X]) = \frac{1}{\hbar} \left[\frac{\exp(\frac{1}{\hbar} [s,\,\cdot\,]_*) - \operatorname{Id}}{\frac{1}{\hbar} [s,\,\cdot\,]_*} (D^\gamma s), \,\cdot\,\right]_*.
	\]
	Consequently,
	\[
	g \cdot D^\gamma = g \circ g^{-1} \circ D^\gamma + g \circ [D^\gamma, g^{-1}] = D^\gamma - \frac{1}{\hbar} \left[ \frac{\exp(\frac{1}{\hbar} [s,\,\cdot\,]_*) - \operatorname{Id}}{\frac{1}{\hbar} [s,\,\cdot\,]_*} D^\gamma s, \,\cdot\,\right]_* = D^{s \bullet \gamma}.
	\]
	
	To show that $R^{s \bullet \gamma} = g(R^\gamma)$, we construct the one-parameter family $\{\gamma(t)\}_{t \in [0, 1]}$ by defining $\gamma(t) := (ts) \bullet \gamma$ for each $t \in [0, 1]$. Note that $\gamma(0) = \gamma$ and $\gamma(1) = s \bullet \gamma$. Using the formula \eqref{Equation: lifted gauge transform} for $\gamma(t)$ and the fact that $\exp(\tfrac{1}{\hbar}[s, \,\cdot\,]_*) (s) = s$, we compute that
	\[
	\frac{d\gamma(t)}{dt} = -D^{\gamma(t)} (s).
	\]
	Using this formula, we compute the derivative of $R^{\gamma(t)} = R^\nabla + \nabla^W \gamma(t) + \tfrac{1}{2\hbar}[\gamma(t), \gamma(t)]_*$ and obtain
	\[
	\frac{d}{dt} R^{\gamma(t)} = -(D^{\gamma(t)})^2 s = \frac{1}{\hbar}[s, R^{\gamma(t)}]_*.
	\]
	On solving, $R^{\gamma(t)} = \exp \left( \tfrac{1}{\hbar}[ts, \,\cdot\,]_*\right) (R^\gamma)$ for all $0 \leq t \leq 1$. Therefore, $R^{s \bullet \gamma} = R^{\gamma(1)} = g(R^\gamma)$.
	
	Finally, suppose $D^\gamma$ is a Fedosov connection on $W$. Then $(g \cdot D^\gamma)^2 = g \circ (D^\gamma)^2 \circ g^{-1} = 0$. Therefore, $g \cdot D^\gamma = D^{s \bullet \gamma}$ is a Fedosov connection. Since $[R^\gamma, \,\cdot\,]_* = 0$, $R^{s \bullet \gamma} = g(R^\gamma) = R^\gamma$.
\end{proof}

\begin{remark}
	\label{Remark: leading term of gauge transformation}
	If $s \in \Gamma(M, W_{(r)})$ for some $r \geq 3$, then 
	\begin{equation*}
		\frac{\exp (\frac{1}{\hbar} [s,\,\cdot\,]_*) - \operatorname{Id}}{\frac{1}{\hbar} [s,\,\cdot\,]_*} (D^\gamma s) \in \operatorname{F}_r\Omega^1(M, W).
	\end{equation*}
	Moreover, its image in $\operatorname{Gr}_r\Omega^1(M, W)$ is precisely $(d_A^{\operatorname{Bott}}  - \delta)(s)$.
\end{remark}

We conclude this subsection by proving a key lemma about Fedosov connections.

\begin{lemma}\label{samecurvature}
	Let $D^\gamma$ and $D^{\gamma'}$ be Fedosov connections on $W$. If $D^\gamma$ and $D^{\gamma'}$ are gauge equivalent, then $[R^\gamma] = [R^{\gamma'}]$. Conversely, if $[R^\gamma] = [R^{\gamma'}]$ and the map $\pi_{0, *}: H^1(M, \mathbb{C}) \to H_A^1(M)$ is surjective, then $D^\gamma$ and $D^{\gamma'}$ are gauge equivalent.
\end{lemma}
\begin{proof}
	For the first statement, there exists $s \in \Gamma(M, W_{(3)})$ such that $D^{\gamma'} = D^{s \bullet \gamma}$. Then we have $s \bullet \gamma = \gamma' + \alpha$ for some $d$-closed $1$-form $\alpha \in \Omega^1(M, \mathbb{C})[[\hbar]] \cap \operatorname{F}_3\Omega^1(M, W)$. By Proposition \ref{Proposition: gauge transformation of Fedosov connection},
	\[
	R^\gamma = R^{s \bullet \gamma} = R^{\gamma'+\alpha} = R^{\gamma'} + d\alpha.
	\]
	
	For the converse, after replacing $\gamma'$ by $\gamma'+\alpha'$ for a suitable $d$-closed $1$-form $\alpha'$ lying in $\Omega^1(M, \mathbb{C})[[\hbar]] \cap \operatorname{F}_3\Omega^1(M, W)$, we may assume that $R^\gamma = R^{\gamma'}$. The desired conclusion is proved by induction on the filtration of $\Omega^*(M, W)$.
	
	As a preliminary step, choose a (not necessarily involutive) complex subbundle $T_{\mathbb C}M$ complementary to $A$. Under the induced identification $\Omega^*(B, W) \cong \Omega_A^{*, *}(M, W)$, we have a bidegree decomposition of $\nabla^W$ as
	\[
	\nabla^W = d_A^{\operatorname{Bott}} + d_B + \iota_R \quad \text{on } \Omega_A^{*, *}(M, W),
	\]
	for a bidegree-$(1, 0)$ operator $d_B$ and $R \in \Omega_A^{2, 0}(M, A \otimes \operatorname{End} (W))$ so that $\iota_R$ is of bidegree $(2, -1)$. On the other hand, for $\alpha, \beta \in \Omega^*(M, W)$, we write
	\[
	\alpha \equiv_{(k)} \beta \quad \text{if and only if} \quad \alpha - \beta \in \operatorname{F}_{k+1}\Omega^*(M, W).
	\]
	Now we assume that for some $n \geq 2$, we have chosen
	\begin{itemize}
		\item a section $s_{\leq n} \in \Gamma(M, W_{(3)})$ having a decomposition $s_{\leq n} = \sum_{k=3}^n s_k$ with $s_k \in \Gamma(M, \operatorname{Gr}_k W)$ for each $3 \leq k \leq n$, and
		\item a $d$-closed $1$-form $\alpha_{\leq n} \in \Omega^1(M, \mathbb{C})[[\hbar]]$ having a decomposition $\alpha_{\leq n} := \sum_{k=3}^n (\alpha_k^B + \alpha_k^A)$, where $\alpha_k^B \in \Omega_A^{1, 0}(M)[[\hbar]]$ is of weight $k-1$ and $\alpha_k^A \in \Omega_A^{0, 1}(M)[[\hbar]]$ is of weight $k$ for each $3 \leq k \leq n$,
	\end{itemize}
	such that
	\begin{equation}
		\label{Equation: gauge equivalence up to an order}
		s_{\leq n} \bullet \gamma \equiv_{(n)} \gamma' + \alpha_{\leq n},
	\end{equation}
	where $\bullet$ is defined in \eqref{Equation: lifted gauge transform}. By definition of a Fedosov connection, $I := s_{\leq n} \bullet \gamma + \tilde{\Omega} \in \operatorname{F}_3\Omega^1(M, W)$. Write $I = \sum_{k = 3}^\infty I_k$, where $I_k \in \operatorname{Gr}_k \Omega^1(M, W)$. Similarly, write $\gamma' + \tilde{\Omega} =: I' = \sum_{k=3}^\infty I'_k$. Then \eqref{Equation: gauge equivalence up to an order} is equivalent to
	\begin{align*}
		I_k = I'_k + \alpha_k^B + \alpha_k^A \quad \text{for all } k = 3, \ldots, n.
	\end{align*}
	The induction hypothesis is vacuous for $n = 2$, upon taking $s_{\leq 2} = 0$ and $\alpha_{\leq 2} = 0$.
	
	Note that $R^{s_{\leq n} \bullet \gamma} = R^\nabla - \Omega + (\nabla^W - \delta) I + \tfrac{1}{2\hbar} \left[ I, I \right]_*$ and similarly for $R^{\gamma'}$. It then follows from $R^{s_{\leq n} \bullet \gamma} = R^\gamma = R^{\gamma'}$ that
	\[
	(\nabla^W - \delta)(I - I') + \frac{1}{2\hbar} [ I, I ]_* - \frac{1}{2\hbar} [I', I']_* = 0.
	\]
	By the induction hypothesis and the fact that $\alpha_{\leq n}$ is central with respect to the commutator $[\,\cdot\,,\,\cdot\,]_*$, the $\operatorname{Gr}_{n+1} \Omega^2(M, W)$-component of the above equation is
	\begin{align*}
		(d_A^{\operatorname{Bott}} - \delta) (I_{n+1} - I'_{n+1}) + d_B (\alpha_n^B + \alpha_n^A) + \iota_R (\alpha_{n-1}^B + \alpha_{n-1}^A) = 0.
	\end{align*}
	Since $d\alpha_{\leq n} = 0$, comparison of the $\Omega_A^{2, 0}(M, \operatorname{Gr}_{n-1} W)$-component yields $d_B\alpha_n^B + \iota_R \alpha_{n-1}^A = 0$, while the $\Omega_A^{1, 1}(M, \operatorname{Gr}_n W)$-component gives $d_B\alpha_n^A = 0$. Moreover, $\iota_R$ annihilates $\Omega_A^{*, 0}(M, W)$, so $\iota_R \alpha_{n-1}^B = 0$. Consequently,
	\begin{align*}
		(d_A^{\operatorname{Bott}} - \delta) (I_{n+1} - I'_{n+1}) = 0.
	\end{align*}
	By Remark \ref{Remark: application of lemma}, it follows that
	\begin{equation}
		\label{Equation: difference in I}
		I_{n+1} - I'_{n+1} = (d_A^{\operatorname{Bott}} - \delta)(x_{n+1}) + \tau,
	\end{equation}
	for some $x_{n+1} \in \Gamma(M, \operatorname{Gr}_{n+1} W)$ and some $d_A$-closed form $\tau \in \Omega_A^{0, 1}(M, \mathbb{C})[[\hbar]]$ of weight $n+1$.
	
	Using the surjectivity of $\pi_{0, *}: H^1(M, \mathbb{C}) \to H_A^1(M)$, choose a form $\alpha_{n+1}^B \in \Omega_A^{1, 0}(M, \mathbb{C})[[\hbar]]$ of weight $n$ and a form $\alpha_{n+1}^A \in \Omega_A^{0, 1}(M, \mathbb{C})[[\hbar]]$ of weight $n+1$ such that $d(\alpha_{n+1}^B + \alpha_{n+1}^A) = 0$ and $\tau - \alpha_{n+1}^A$ is $d_A$-exact. Thus there exists $y_{n+1} \in \Omega^0(M, \mathbb{C})[[\hbar]]$ of weight $n+1$ satisfying $\tau - \alpha_{n+1}^A = d_Ay_{n+1}$.
	
	Define $\alpha_{\leq (n+1)} = \alpha_{\leq n} + \alpha_{n+1}^B + \alpha_{n+1}^A$. Because $\delta y_{n+1} = 0$, \eqref{Equation: difference in I} implies that
	\begin{align*}
		s_{\leq n} \bullet \gamma
		- \gamma' \equiv_{(n+1)} (d_A^{\operatorname{Bott}} - \delta)(s_{n+1}) + \alpha_{\leq (n+1)},
	\end{align*}
	where $s_{n+1}:=x_{n+1}+y_{n+1} \in \Gamma(M, \operatorname{Gr}_{n+1}W)$. By Remark \ref{Remark: leading term of gauge transformation}, we obtain
	\begin{align*}
		s_{n+1} \bullet (s_{\leq n} \bullet \gamma) & \equiv_{(n+1)} \gamma' + \alpha_{\leq (n+1)}.
	\end{align*}
	Because $s_{n+1}$ has weight $n+1$ and $s_{\leq n} \in \Gamma(M, W_{(3)})$, we have $\tfrac{1}{\hbar} [s_{n+1}, s_{\leq n}]_* \in \Gamma(M, W_{(n+2)})$. Therefore, using the Baker–Campbell–Hausdorff formula, we can deduce that the successive gauge actions combine additively modulo $\operatorname{F}_{n+2}\Omega^1(M, W)$:
	\begin{equation*}
		s_{\leq (n+1)} \bullet \gamma \equiv_{(n+1)} \gamma' + \alpha_{\leq (n+1)},
	\end{equation*}
	where $s_{\leq (n+1)} := s_{\leq n} + s_{n+1}$. This completes the induction argument.
	
	Eventually, the sequences $\{ s_{\leq n}\}_{n=2}^\infty$ and $\{\alpha_{\leq n}\}_{n=2}^\infty$ converge to a section $s \in \Gamma(M, W_{(3)})$ and a $d$-closed $1$-form $\alpha \in \Omega^1(M, \mathbb{C})[[\hbar]] \cap \operatorname{F}_3\Omega^1(M, W)$ in the filtration topology, respectively. In the limit, the relation becomes exact:
	\begin{equation*}
		s \bullet \gamma = \gamma' + \alpha.
	\end{equation*}
	Therefore, $D^\gamma$ is gauge equivalent to $D^{\gamma'+\alpha} = D^{\gamma'}$ via the gauge transformation $\exp (\tfrac{1}{\hbar}[s, \,\cdot\,]_*)$.
\end{proof}

\subsection{Fedosov quantization}
\label{4.4}
\quad\par
In this subsection, we elaborate how to transform the classification of deformation quantization of $(M, \Omega)$, up to equivalence, into the classification of Fedosov connections, up to gauge transformations of the Weyl bundle.

Associated to each Fedosov connection $D^\gamma$ on $W$ is a pair
\begin{equation}
	\label{Equation: DQ from Fedosov}
	(\mathcal{W}_{D^\gamma}, \psi_\gamma)
\end{equation}
constructed as follows. We define \(\mathcal{W}_{D^\gamma}\) as the sheaf of $D^\gamma$-flat sections of $W$. Since $D^\gamma$ is a $*$-derivation, the fibrewise Moyal product $*$ descends to a product on $\mathcal{W}_{D^\gamma}$ so that $\mathcal{W}_{D^\gamma}$ forms a sheaf of $\mathbb{C}[[\hbar]]$-algebras. On the other hand, recall that we introduced a decreasing filtration on $\Omega^*(M, W)$ in Section \ref{4.1}. Let $\sigma: \Omega^*(M, W) \to \operatorname{Gr}_0 \Omega^*(M, W) = \Omega_A^{0, *}(M)$ be the canonical projection. The leading term of the associated graded operator of $D^\gamma = \nabla^W + \tfrac{1}{\hbar}[\gamma, \,\cdot\,]_*$ restricted to $\Omega_A^{0, *}(M)$ is precisely $d_A$. Consequently,
$$\sigma: (\Omega^*(M, W), D^\gamma) \to (\Omega_A^{0, *}(M), d_A)$$
is a morphism of differential graded algebras. Restricting $\sigma$ to $D^\gamma$-flat sections therefore yields a morphism of $\mathbb{C}$-algberas $\psi_\gamma: \mathcal{W}_{D^\gamma} \to \mathcal{O}$.
We now state the main theorems of this subsection.

\begin{theorem}
	\label{Theorem: Fedosov quantization is DQ}
	Let $D^\gamma$ be a Fedosov connection on $W$. Then the pair $(\mathcal{W}_{D^\gamma}, \psi_\gamma)$ in \eqref{Equation: DQ from Fedosov} is a deformation quantization of the transverse holomorphic symplectic manifold $(M, \Omega)$.
\end{theorem}

To be more precise, we should specify the differential structure with which $(\mathcal{W}_{D^{\gamma}}, \psi_\gamma)$ is equipped. This structure will be constructed in the proof.

\begin{theorem}
	\label{Theorem: gauge equivalence}
	Two Fedosov connections $D^\gamma$ and $D^{\gamma'}$ are gauge equivalent if and only if the deformation quantizations $(\mathcal{W}_{D^\gamma}, \psi_\gamma)$ and $(\mathcal{W}_{D^{\gamma'}}, \psi_{\gamma'})$ are equivalent.
\end{theorem}

\begin{theorem}\label{dqtofedosov}
	Every deformation quantization of $(M, \Omega)$ is equivalent to one of the form $(\mathcal{W}_{D^\gamma}, \psi_\gamma)$ as in \eqref{Equation: DQ from Fedosov} for some Fedosov connection $D^\gamma$ on the Weyl bundle $W$.
\end{theorem}

By these theorems and Lemma \ref{samecurvature}, we can define, without additional assumptions, the map
\begin{equation}
	\label{Equation: period map}
	\operatorname{per}: \operatorname{DQ}(M, \Omega) \to H^2(M, \mathbb{C})[[\hbar]]
\end{equation}
by assigning to an equivalence class of deformation quantizations of $(M,\Omega)$ the cohomology class $[R^\gamma]$, where $D^\gamma$ is any Fedosov connection for which $(\mathcal{W}_{D^\gamma}, \psi_\gamma)$ represents the given equivalence class. Following the terminology of \cite{BezKal2004}, we call $\operatorname{per}$ the \emph{non-commutative period map}.

\begin{remark}
	Indeed, the image of $\operatorname{per}$ lies in $-[\Omega] + \hbar H^2(M, \mathbb{C})[[\hbar]]$. To see this, note that for any $I \in \operatorname{F}_3\Omega^1(M, W)$,
	\[
	R^{-\tilde{\Omega} + I} = R^\nabla - \Omega + (\nabla^W - \delta)I + \tfrac{1}{2\hbar} [I, I]_*.
	\]
	Here, $R^\nabla$ has no component in $\Omega^2(M, \mathbb{C})[[\hbar]]$, while $(\nabla^W - \delta)I + \tfrac{1}{2\hbar} [I, I]_* \in \operatorname{F}_3\Omega^2(M, W)$. Consequently, the $\Omega^2(M, \mathbb{C})$-component of $R^{-\tilde{\Omega} + I}$ is simply $-\Omega$.
\end{remark}

The proofs of the above theorems require the transverse holomorphic jet bundle as a bridge between local star products and Fedosov constructions. Hence, we first explain how star products on $\mathcal O[[\hbar]]$ and gauge equivalences prolong to the transverse holomorphic jet bundle. We also record the basic properties of these prolongations that will be used in the proof of the above theorem.

Let $\star$ be a star product on $\mathcal O[[\hbar]]$. Write $f \star g = \sum_{r=0}^\infty \hbar^r C_r(f, g)$, where $C_r$'s are transverse holomorphic bidifferential operators. This product induces a $\mathcal{C}^\infty(M)[[\hbar]]$-bilinear map
\begin{equation}
	\label{Equation: prolongation of star product}
	\widetilde{\star}: \Gamma(M, J_AM)[[\hbar]] \times \Gamma(M, J_AM)[[\hbar]] \to \Gamma(M, J_AM)[[\hbar]]
\end{equation}
defined by
$$a \,\widetilde{\star}\, b
:=
\sum_{r=0}^\infty \hbar^r
j^\infty(C_r) (a, b)$$
for all $a, b \in \Gamma(M, J_AM)[[\hbar]]$, where each $j^\infty(C_r)$ is the prolongation of $C_r$ (see Subsection \ref{Subsection: transverse jet bundles}). We call $\widetilde{\star}$ the \emph{prolongation} of $\star$. It has the following properties:
\begin{enumerate}
	\item $(\Gamma(M, J_AM)[[\hbar]],\widetilde\star)$ forms a $\mathcal{C}^\infty(M)[[\hbar]]$-algebra;
	\item for each $v\in\mathcal B(M)$, $\nabla_v^{\mathrm G}$ is a derivation of $\widetilde\star$, namely for all $s,t\in \Gamma(M, J_AM)[[\hbar]]$,
	\[
	\nabla^{\mathrm G}_v(s\,\widetilde\star\, t)
	=
	(\nabla^{\mathrm G}_v s)\,\widetilde\star\, t
	+
	s\,\widetilde\star\,(\nabla^{\mathrm G}_v t);
	\]
	
	\item modulo $\hbar$, the product $\widetilde\star$ reduces to the ordinary fiberwise multiplication on $J_AM$;
	
	\item the prolongation map is an algebra homomorphism:
	\[
	j^\infty(f\star g)
	=
	j^\infty(f)\,\widetilde\star\, j^\infty(g),
	\qquad
	f,g\in\mathcal O(U)[[\hbar]].
	\]
\end{enumerate}

For a gauge equivalence $T = \mathrm{Id} + \sum_{r=1}^\infty \hbar^r T_r$ between two star products on $\mathcal{O}[[\hbar]]$, where each $T_r$ is a transverse holomorphic differential operator, we similarly define the \emph{prolongation} of $T$ as the $\mathcal{C}^\infty(M)$-linear map
\begin{equation}
	\label{Equation: prolongation of gauge equivalence}
	\widetilde{T} := \mathrm{Id} + \sum_{r=1}^\infty \hbar^r j^\infty(T_r): \Gamma(M, J_AM)[[\hbar]] \to \Gamma(M, J_AM)[[\hbar]]
\end{equation}
This operator commutes with $\nabla^{\operatorname{G}}$, i.e. $\widetilde{T} \circ \nabla^{\operatorname{G}} = \nabla^{\operatorname{G}} \circ \widetilde{T}$.

Now, we make the following crucial observations. By Theorem \ref{Theorem: Darboux}, we can choose an open cover $\{U_i\}_{i \in I}$ of $M$ such that each $U_i$ is equipped with transverse holomorphic Darboux coordinates, which defines the Moyal product $\star_i$ on $\mathcal{O} \vert_{U_i} [[\hbar]]$ (see Example \ref{Example: Moyal product}).

Fix $i \in I$. Let $\widetilde{\star}_i$ be the prolongation of $\star_i$, given as in \eqref{Equation: prolongation of star product}. The fixed transverse holomorphic Darboux coordinates $(x^1, ..., x^m, z^1, ..., z^{2n})$ on $U_i$ induce a $\mathcal{C}^\infty(U_i)[[\hbar]]$-algebra isomorphism
\[
\phi_i: (\Gamma(U_i, J_AM)[[\hbar]], \widetilde{\star}_i) \to (\Gamma(U_i, W), * \vert_{U_i}),
\]
obtained by extending \eqref{Equation: local trivialization of infinite jet bundle} $\mathcal{C}^\infty(U_i)[[\hbar]]$-linearly. By \eqref{Equation: local formula of Grothendieck connection}, we have $\phi_i \circ \nabla^{\operatorname{G}} \circ \phi_i^{-1} = d_i - \delta$, where 
\[
d_i := dx^j \wedge \frac{\partial}{\partial x^j} + d\overline{z}^k \wedge \frac{\partial}{\partial \overline{z}^k} + dz^l \wedge \frac{\partial}{\partial z^l} \quad \text{on } \Gamma(U_i, W).
\]
Recall that we have fixed a symplectic connection $\nabla$ on $(M, \Omega)$, which extends the Bott connection on $B$. The induced connection $\nabla^W$ on the Weyl bundle $W$ takes the local form:
\[
\nabla^W \vert_{U_i} = d_i - \frac{1}{\hbar}[\Gamma_i, \,\cdot\,]_*,
\]
where $\Gamma_i \in \Omega_A^{1, 0}(U_i, \operatorname{Sym}^2 (B^*)) \subset \operatorname{F}_3\Omega^1(U_i, W)$. Consequently,
\[
D_i := \phi_i \circ \nabla^{\operatorname{G}} \circ \phi_i^{-1} = \nabla^W \vert_{U_i} - \delta + \frac{1}{\hbar}[\Gamma_i, \,\cdot\,]_* = \nabla^W \vert_{U_i} + \frac{1}{\hbar}[-\tilde{\Omega} + \Gamma_i, \,\cdot\,]_*
\]
is a Fedosov connection. Furthermore, $\sigma \circ \phi_i$ coincides with the composition
\[
\Gamma(U_i, J_AM)[[\hbar]] \to \Gamma(U_i, J_AM) \to \mathcal{C}^\infty(U_i).
\]
Restricting to $D_i$-flat sections yields a sheaf map $\psi_i: \mathcal{W}_{D_i} \to \mathcal{O} \vert_{U_i}$. Via the Moyal product $\star_i$ and the above constructions, we can obviously see that $\mathcal{W}_{D_i}$ is a sheaf of flat, $\hbar$-adically complete $\mathbb{C}[[\hbar]]$-algebras on $U_i$, that there is a short exact sequence:
\begin{equation*}
	\begin{tikzcd}
		0 \ar[r] & \mathcal{W}_{D_i} \ar[r, "\hbar (\,\cdot\,)"] & \mathcal{W}_{D_i} \ar[r, "\psi_i"] & \mathcal{O} \ar[r] & 0,
	\end{tikzcd}
\end{equation*}
and that for all $\widetilde{f}, \widetilde{g} \in \mathcal{W}_{D_i}(U_i)$, $\{\psi_i(\widetilde{f}), \psi_i(\widetilde{g})\} = \psi_i\left(\tfrac{1}{\hbar} [\widetilde{f}, \widetilde{g}]_* \right)$. In other words, we have shown a special case of Theorem \ref{Theorem: Fedosov quantization is DQ} that $(\mathcal{W}_{D_i}, \psi_i)$ is a deformation quantization of $(U_i, \Omega \vert_{U_i})$.

We make a further observation on each overlap $U_{ij} := U_i \cap U_j$, where $i, j \in I$. Suppose that
\[
T_{ij}: (\mathcal O[[\hbar]]|_{U_{ij}},\star_i)
\longrightarrow
(\mathcal O[[\hbar]]|_{U_{ij}},\star_j)
\]
is a gauge equivalence in the sense of Definition \ref{Definition: gauge equivalence}. Its prolongation
\[
\widetilde{T}_{ij}
:
(J_AU_{ij} [[\hbar]], \widetilde{\star}_i) \to (J_AU_{ij} [[\hbar]], \widetilde{\star}_j),
\]
given as in \eqref{Equation: prolongation of gauge equivalence}, commutes with $\nabla^{\operatorname{G}}$. Define the composition
\[
g_{ij} := \phi_j \circ \widetilde{T}_{ij} \circ \phi_i^{-1}: (\Gamma(U_{ij}, W), * \vert_{U_{ij}}) \to (\Gamma(U_{ij}, W), * \vert_{U_{ij}}).
\]

Note that $g_{ij} - \mathrm{Id} = \phi_j \circ (\widetilde{T}_{ij} - \operatorname{Id}) \circ \phi_i^{-1} + ( \phi_j \circ \phi_i^{-1} - \operatorname{Id} )$ strictly increases the weight of $W \vert_{U_{ij}}$. The first term does so because $T_{ij}$ is a gauge equivalence, while the second does so because $\phi_i$ and $\phi_j$ induce the same canonical identification on the associated graded of the transverse holomorphic jet bundle. Hence, $g_{ij}$ is a gauge transformation of $W \vert_{U_{ij}}$. Also, $\widetilde{T}_{ij} \circ \nabla^{\operatorname{G}} = \nabla^{\operatorname{G}} \circ \widetilde{T}_{ij}$ implies $g_{ij} \circ D_i = D_j \circ g_{ij}$.

Conversely, suppose that $g_{ij}$ is a gauge transformation of $W \vert_{U_{ij}}$ intertwining $D_i$ and $D_j$. Clearly, $\widetilde{T}_{ij} := \phi_j^{-1} \circ g_{ij} \circ \phi_i$  restricts to an automorphism of the sheaf of $\nabla^{\operatorname{G}}$-flat sections of $J_AU_{ij}[[\hbar]]$. We obtain an isomorphism of sheaves of $\mathbb{C}[[\hbar]]$-algebras
\[
T_{ij}: (\mathcal O[[\hbar]]|_{U_{ij}},\star_i)
\longrightarrow
(\mathcal O[[\hbar]]|_{U_{ij}},\star_j)
\]
by the defining relation $j^\infty \circ T_{ij} = \widetilde{T}_{ij} \circ j^\infty$ on $\mathcal{O} \vert_{U_{ij}}[[\hbar]]$. More explicitly,
\[
T_{ij} = \pi_{0, *} \circ g_{ij} \circ \phi_i \circ j^\infty \quad \text{on } \mathcal{O} \vert_{U_{ij}} [[\hbar]],
\]
where $\pi_{0, *}$ is given in \eqref{Equatoin: pi zero star}, noting that the composition $\pi_{0, *} \circ \phi_j$ coincides with the canonical projection $\Gamma(U_j, J_AM)[[\hbar]] \to C^\infty(U_j)[[\hbar]]$.

We now check that $T_{ij}$ admits an expansion of the form $T_{ij} = \operatorname{Id} + \sum_{r=1}^\infty \hbar^r T_{ij, r}$, where for all $r > 0$, $T_{ij, r}$ is a transverse holomorphic bidifferential operator. Recall that we have fixed transverse holomorphic Darboux coordinates on $U_i$. Let $u^k$ denote the local section $dz^k$ of $B^*$. The condition that $g_{ij} - \mathrm{Id}$ strictly increases the weight implies that, for all multi-index $\alpha$ on $\{1, \ldots, 2n\}$, $g_{ij}(u^\alpha)$ is of the form
\[
g_{ij}(u^\alpha) = \sum_{\lvert \beta \rvert + 2r \geq \lvert \alpha \rvert} \hbar^r a_{\beta, r}^\alpha u^\beta,
\]
where $a_{\beta, r}^\alpha\in \mathbb{C}$ and $a_{\alpha, 0}^\alpha = 1$. Let $U \subset U_{ij}$ be an open subset and $f \in \mathcal{O}(U)$. Recalling that
\[
\phi_i (j^\infty f) = \sum_\alpha \frac{1}{\alpha!} \frac{\partial^{\lvert \alpha \rvert} f}{\partial z^\alpha} u^\alpha,
\]
we obtain
\[
T_{ij}(f) = f + \sum_{r=1}^\infty \hbar^r \sum_{\lvert \alpha \rvert \leq 2r} \frac{a_{\emptyset, r}^\alpha}{\alpha!} \frac{\partial^{\lvert \alpha \rvert} f}{\partial z^\alpha},
\]
where $\emptyset$ denotes the empty index. Therefore, $T_{ij}$ is of the desired form. All in all, the above arguments prove a special case of Theorem \ref{Theorem: gauge equivalence} that $D_i \vert_{U_{ij}} $ and $D_j \vert_{U_{ij}}$ are gauge equivalent if and only if the deformation quantizations $(\mathcal{W}_{D_i} \vert_{U_{ij}}, \psi_i \vert_{U_{ij}})$ and $(\mathcal{W}_{D_j} \vert_{U_{ij}}, \psi_j \vert_{U_{ij}})$ are equivalent.

Now, we proceed to the general proofs of the above theorems.

\begin{proof}[\myproof{Theorem}{\ref{Theorem: Fedosov quantization is DQ}}]
	Choose an open cover $\{U_i\}_{i \in I}$ of $M$ in the form as in the above discussion and follow the notations for the relevant local data. Refine this open cover so that for each $i \in I$, $H^1(U_i, \mathbb{C}) = 0$ and $H_A^1(U_i) = 0$, and thus the assumption of Lemma \ref{samecurvature} is fulfilled. Then by that lemma, there exists a gauge transformation $g_i$ of $W \vert_{U_i}$ such that $D^\gamma \vert_{U_i} = g_i \circ D_i \circ g_i^{-1}$, therefore $g_i$ maps $\mathcal{W}_{D^\gamma} \vert_{U_I}$ to $\mathcal{W}_{D_i}$. Since $g_i - \operatorname{Id}$ strictly increases the weight of elements in the Weyl bundle $W$, we obtain $\psi_i \circ g_i = \psi_\gamma \vert_{U_i}$. The conditions that the sheaf of $\mathbb{C}[[\hbar]]$-algebras $\mathcal{W}_{D^\gamma}$ is flat and $\hbar$-complete and that $\mathcal{W}_{D^\gamma} \xrightarrow{\hbar (\,\cdot\,)} \mathcal{W}_{D^\gamma} \xrightarrow{\psi_\gamma} \mathcal{O}$ is a short exact sequence are local. It follows that $(\mathcal{W}_{D^\gamma}, \psi_\gamma)$ inherits these properties from the local models $(\mathcal{O} \vert_{U_i}[[\hbar]], \star_i)$.
	
	Consequently, $(U_i, \star_i, \tau_i)$ is a local differential trivialization of $(\mathcal{W}_{D^\gamma}, \psi_\gamma)$, where $\tau_i$ is determined by the condition that for all local $D^\gamma$-flat section $a$ of $W \vert_{U_i}$,
	\[
	(j^\infty \circ \tau_i)(a) = (\phi_i^{-1} \circ g_i)(a).
	\]
	Fix $i, j \in I$. The composition $g_j \circ g_i^{-1}$ is a gauge transformation of $W \vert_{U_{ij}}$. The discussion before this proof shows that
	\[
	\tau_j \circ \tau_i^{-1}: (\mathcal{O} \vert_{U_{ij}} [[\hbar]], \star_i) \to (\mathcal{O} \vert_{U_{ij}} [[\hbar]], \star_j)
	\]
	is a gauge equivalence between $\star_i$ and $\star_j$. Thus, $\{ (U_i, \star_i, \tau_i) \}_{i \in I}$ is a differential atlas of $(\mathcal{W}_{D^\gamma}, \psi_\gamma)$.
	
	Again, the condition \eqref{Equation: semiclassical to Poisson} is local, and is hence implied by the local models. To conclude, equipped with the differential structure containing the above local differential trivializations, $(\mathcal{W}_{D^\gamma}, \psi_\gamma)$ is a deformation quantization of $(M, \Omega)$.
\end{proof}

\begin{proof}[\myproof{Theorem}{\ref{Theorem: gauge equivalence}}]
	Let $U_i, \star_i, \tau_i, g_i$ be the local data given in the proof of Theorem \ref{Theorem: Fedosov quantization is DQ}. Define $g'_i, \tau'_i$ as the $(\mathcal{W}_{D^{\gamma'}}, \psi_{\gamma'})$-counterparts of $g_i, \tau_i$. Then $(U_i, \star_i, \tau_i)$ and $(U_i, \star_i, \tau'_i)$ are adapted local differential trivializations of $(\mathcal{W}_{D^\gamma}, \psi_\gamma)$ and $(\mathcal{W}_{D^{\gamma'}}, \psi_{\gamma'})$, respectively.
	
	Suppose that there exists an equivalence $\Phi: (\mathcal{W}_{D^\gamma}, \psi_\gamma) \to (\mathcal{W}_{D^{\gamma'}}, \psi_{\gamma'})$ of deformation quantizations of $(M, \Omega)$. Extend $\Phi$ to a bundle isomorphism $g: W \to W$ by $\mathcal{C}^\infty(M)[[\hbar]]$-linearity. For each $i \in I$, by the definition of $\Phi$, we obtain a gauge equivalence:
	\[
	T_i := \tau_i' \circ \Phi \circ \tau_i^{-1} : (\mathcal{O} \vert_{U_i}[[\hbar]], \star_i) \to (\mathcal{O} \vert_{U_i}[[\hbar]], \star_i).
	\]
	It follows from the discussion before the proof of Theorem \ref{Theorem: Fedosov quantization is DQ} that $\phi_i \circ \widetilde{T}_i \circ \phi_i^{-1}$ is a gauge transformation of $W \vert_{U_i}$ which commutes with $D_i$, where $\widetilde{T}_i$ denotes the prolongation of $T_i$. On the other hand, by construction, the equality
	\[
	\Phi (a) = ( (g'_i)^{-1} \circ \phi_i \circ \widetilde{T}_i \circ \phi_i^{-1} \circ g_i ) (a)
	\]
	holds for all local $D^\gamma$-flat section $a$ of $W \vert_{U_i}$. Taking $\mathcal{C}^\infty(U_{ij})[[\hbar]]$-linear extensions on both sides, we conclude that $g$ is locally the composition of the following gauge transformations:
	\[
	g \vert_{U_i} = (g'_i)^{-1} \circ \phi_i \circ \widetilde{T}_i \circ \phi_i^{-1} \circ g_i.
	\]
	Therefore, $g$ is a gauge transformation of $W$ intertwining $D^\gamma$ and $D^{\gamma
		'}$.
	
	It remains to prove the converse. Suppose that there exists a gauge transformation $g$ of $W$ such that $D^{\gamma'} = g \cdot D^\gamma$. Then it restricts to an isomorphism of $\mathbb{C}[[\hbar]]$-algebras $\Phi: \mathcal{W}_{D^\gamma} \to \mathcal{W}_{D^{\gamma'}}$ such that $\psi_{\gamma'} \circ g = \psi_\gamma$.  Fix $i, j \in I$. The composition $g'_i \circ g \circ g_i^{-1}$ is a local gauge transformation of $W$. By the discussion before the proof of Theorem \ref{Theorem: Fedosov quantization is DQ} again, $\tau'_i \circ \Phi \circ \tau_i^{-1}$ is a gauge equivalence between $\star_i$ and itself. Therefore, $(\mathcal{W}_{D^\gamma}, \psi_\gamma)$ and $(\mathcal{W}_{D^{\gamma'}}, \psi_{\gamma'})$ are equivalent.
\end{proof}

\begin{proof}[\myproof{Theorem}{\ref{dqtofedosov}}]
	Consider any deformation quantization $(\mathcal{A}_M^\hbar, \psi)$ of $(M, \Omega)$. By Proposition \ref{localstar}, we can choose a differential atlas $\{(U_i ,\star_i, \tau_i)\}_{i \in I}$ on $(\mathcal{A}_M^\hbar, \psi)$ such that each $\star_i$ is the Moyal product with respect to a transverse holomorphic Darboux chart on $U_i$. We then follow the notations for the relevant local data in the discussion before the proof of Theorem \ref{Theorem: Fedosov quantization is DQ}.
	
	\textbf{Step 1:} We first associate to $(\mathcal{A}_M^\hbar, \psi)$ a \v{C}ech $1$-cocycle with values in the sheaf $\mathcal{G}$ of groups of gauge transformations of $W$.
	
	On each overlap $U_{ij}$, the transition map
	\[
	T_{ij}:=\tau_j\circ\tau_i^{-1}
	:
	(\mathcal O[[\hbar]]|_{U_{ij}},\star_i)
	\longrightarrow
	(\mathcal O[[\hbar]]|_{U_{ij}},\star_j)
	\]
	is a gauge equivalence. It follows from the discussion before the proof of Theorem \ref{Theorem: Fedosov quantization is DQ} that $g_{ij} := \phi_j \circ \widetilde{T}_{ij} \circ \phi_i^{-1}$ is a gauge transformation of $W \vert_{U_{ij}}$, where $\widetilde{T}_{ij}$ is the prolongation of $T_{ij}$. Because the maps $\widetilde{T}_{ij}$ satisfy the cocycle condition, $\{g_{ij}\}$ forms a \v{C}ech $1$-cocycle with values in $\mathcal{G}$, i.e. $g_{jk} \circ g_{ij} = g_{ik}$; see Subsection \ref{4.1}.
	
	\textbf{Step 2:} Recall from Subsection \ref{4.1} that $\mathcal{G}$ admits a decreasing filtration. We show, by induction on this filtration, that the cocycle $\{g_{ij}\}$ is the coboundary of a $0$-cochain with values in $\mathcal{G}$.
	
	Assume we have constructed a gauge transformation $\phi_i^{(r)}$ of $W \vert_{U_i}$ for each $i$ such that
	\[
	g_{ij}^{(r)} := (\phi_j^{(r)})^{-1} \circ g_{ij} \circ \phi_i^{(r)} \in \mathcal{G}^{(r)}(U_{ij}), \quad \text{for all } i, j \in I.
	\]
	The base case $r=1$ holds trivially by setting $\phi_i^{(1)} = \mathrm{Id}$.
	
	Note that the quotient $\mathcal{G}^{(r)} / \mathcal{G}^{(r+1)}$ is canonically isomorphic to the Abelian sheaf $\mathcal{A}^{(r)}$ given in \eqref{Equation: graded piece of derivations}. Clearly, $\mathcal{A}^{(r)}$ is a locally free $\mathcal{C}_M^\infty$-module, hence a fine sheaf. Consequently, $\check{H}^1(\mathcal{U}, \mathcal{A}^{(r)}) = 0$, where $\mathcal{U} := \{U_i\}_{i \in I}$. The image of $\{g_{ij}^{(r)}\}$ in $\mathcal{G}^{(r)} / \mathcal{G}^{(r+1)} \cong \mathcal{A}^{(r)}$ is a \v{C}ech $1$-cocycle. Thus, there exists a $0$-cochain $\{s_i^{(r)}\}$ with values in $\mathcal{A}^{(r)}$ such that
	\begin{equation*}
		\log (g_{ij}^{(r)}) - \tfrac{1}{\hbar} [s_j^{(r)} - s_i^{(r)}, \,\cdot\,]_\star \in \operatorname{Der}_{\geq (r+1)}(U_{ij}), \quad \text{for all } i, j \in I.
	\end{equation*}
	We define the corrected gauge transformation for the subsequent inductive step by setting
	\begin{equation*}
		\phi_i^{(r+1)} := \phi_i^{(r)} \circ \exp(\tfrac{1}{\hbar}[s_i^{(r)}, \,\cdot\,]_*).
	\end{equation*}
	Using the Baker--Campbell--Hausdorff formula, one verifies that
	$$g_{ij}^{(r+1)} := (\phi_j^{(r+1)})^{-1} \circ g_{ij} \circ \phi_i^{(r+1)} \in \mathcal{G}^{(r+1)}(U_{ij}).$$
	
	Since each correction belongs to $\mathcal{G}^{(r)}(U_i)$, the sequence $\{\phi_i^{(r)}\}$ is Cauchy with respect to the filtration topology on $\mathcal{G}(U_i)$. Because $\mathcal{G}(U_i)$ is complete, it converges to a gauge transformation $\phi_i^{(\infty)} \in \mathcal{G}(U_i)$. Since $\bigcap_{r = 1}^\infty \mathcal{G}^{(r)}(U_{ij}) = \{\mathrm{Id}\}$, we obtain $(\phi_j^{(\infty)})^{-1} \circ g_{ij} \circ \phi_i^{(\infty)} = id$, or equivalently,
	\[
	g_{ij} = \phi_j^{(\infty)} \circ (\phi_i^{(\infty)})^{-1}.
	\]
	
	\textbf{Step 3:} Because $\widetilde{T}_{ij} \circ \nabla^{\operatorname{G}} \circ \widetilde{T}_{ij}^{-1} = \nabla^{\operatorname{G}}$, we can globally define a connection $D$ on $W$ by setting
	\[
	D \vert_{U_i} := (\phi_i^{(\infty)})^{-1} \circ D_i \circ \phi_i^{(\infty)}, \quad \text{for each } i \in I.
	\]
	This connection inherits the flatness from $\nabla^{\operatorname{G}}$. We claim that $D$ is a Fedosov connection on $W$.
	
	Choose a section $s_i \in \Gamma(U_i, W_{(3)})$ such that $\phi_i^{(\infty)} = \exp\left(-\frac{1}{\hbar}[s_i, \,\cdot\,]_* \right)$. By Proposition \ref{Proposition: gauge transformation of Fedosov connection},
	\[
	D \vert_{U_i} = \nabla^W + \tfrac{1}{\hbar} [-\tilde{\Omega} + \Gamma_i + \tilde{\gamma}_i, \,\cdot\,]_*.
	\]
	where
	\begin{align*}
		\tilde{\gamma}_i := -\frac{\exp(\frac{1}{\hbar} [s_i,\,\cdot\,]_*) - \operatorname{Id}}{\frac{1}{\hbar} [s_i,\,\cdot\,]_*} \circ D_i s_i \in \operatorname{F}_3\Omega^1(U_i, W).
	\end{align*}
	
	Furthermore, note that $\tfrac{1}{\hbar} [\Gamma_i + \tilde{\gamma}_i, \,\cdot\,]_* = \tfrac{1}{\hbar} [\Gamma_j + \tilde{\gamma}_j, \,\cdot\,]_*$ on the overlap $U_{ij}$. This implies that
	\[
	(\Gamma_i + \tilde{\gamma}_i) - (\Gamma_j + \tilde{\gamma}_j) \in \Omega^1(U_{ij}, \mathbb{C})[[\hbar]] \cap \operatorname{F}_3\Omega^1(U_{ij}, W).
	\]
	Since the sheaf of smooth complex $1$-forms on $M$ is fine, we can adjust the local $W$-valued $1$-forms $\Gamma_i + \tilde{\gamma}_i$ by a coboundary to glue them into a global section $\tilde{\gamma} \in \operatorname{F}_3\Omega^1(M, W)$. This confirms that $D = D^\gamma$, where $\gamma := -\tilde{\Omega} + \tilde{\gamma}$, is a Fedosov connection, completing the proof of our claim.
	
	\textbf{Step 4:} We claim that we can well define an isomorphism of sheaves of $\mathbb{C}[[\hbar]]$-algebras
	\[
	\Phi: \mathcal{A}_M^\hbar \to \mathcal{W}_D,
	\]
	where $\mathcal{W}_D$ denotes the sheaf of $D$-flat sections of $W$, by $\Phi \vert_{U_i} = (\phi_i^{(\infty)})^{-1} \circ \phi_i \circ j^\infty \circ \tau_i$ for all $i \in I$. To prove this claim, fix $i, j \in I$. On $\mathcal{A}_M^\hbar \vert_{U_{ij}}$, we compute:
	\begin{align*}
		(\phi_i^{(\infty)})^{-1} \circ \phi_i \circ j^\infty \circ \tau_i = & (\phi_j^{(\infty)})^{-1} \circ g_{ij} \circ \phi_i \circ j^\infty \circ \tau_i\\
		= & (\phi_j^{(\infty)})^{-1} \circ \phi_j  \circ \widetilde{T}_{ij} \circ j^\infty \circ \tau_i\\
		= & (\phi_j^{(\infty)})^{-1} \circ \phi_j  \circ j^\infty \circ T_{ij} \circ \tau_i\\
		= & (\phi_j^{(\infty)})^{-1} \circ \phi_j  \circ j^\infty \circ \tau_j.
	\end{align*}
	In the first line, we used the result in Step 2. In the second and the fourth lines, we used the definitions of $g_{ij}$ and $T_{ij}$, respectively. In the third line, we used a basic property of the prolongation $\widetilde{T}_{ij}$ of $T_{ij}$. Finally, straightforward arguments show that $\psi_\gamma \circ \Phi = \psi$ and, via adapted local differential trivializations of $(\mathcal{A}_M^\hbar, \psi)$ and of $(\mathcal{W}_D, \psi_\gamma)$, $\Phi$ is identified with a gauge equivalence between star products, thereby completing the proof.
\end{proof}

\subsection{Classification theorem for deformation quantization}\label{4.5}
\quad\par
In this subsection, we prove our main result on classification of deformation quantizations of $(M, \Omega)$ via the non-commutative period map given in \eqref{Equation: period map}. We begin with proving the existence of Fedosov connections under the assumption that the map
\[
\pi_{0, *}: H^2(M, \mathbb{C}) \to H_A^2(M)
\]
is surjective. Whenever this assumption holds, we fix a $\mathbb{C}$-linear right inverse $H_A^2(M) \to H^2(M, \mathbb{C})$ of it. In addition, we fix a (not necessarily $\mathbb{C}$-linear) lifting map associated with this right inverse:
$$\Omega_{A, d_A\text{-closed}}^{0, 2}(M) \to \Omega_{d\text{-closed}}^2(M, \mathbb{C}), \quad \tau \mapsto \Lambda_\tau.$$
In particular, $[\pi_{0, *}(\Lambda_\tau)] = [\tau] \in H_A^2(M)$. We further require the normalization condition $\Lambda_0 = 0$. We extend the lifting map to a map
$$\Omega_{A, d_A\text{-closed}}^{0, 2}(M)[[\hbar]] \to \Omega_{d\text{-closed}}^2(M, \mathbb{C})[[\hbar]], \quad \tau \mapsto \Lambda_\tau,$$
by defining, for $\tau = \sum_{r=0}^\infty \hbar^r \tau_r$ with $\tau_r \in \Omega_{A, d_A\text{-closed}}^{0, 2}(M)$, $\Lambda_\tau := \sum_{r=0}^\infty \hbar^r \Lambda_{\tau_r}$.

\begin{theorem}
	\label{Theorem 3.18}
	Let $(M, \Omega)$ be a transverse holomorphic symplectic manifold. Assume that the map $\pi_{0, *}: H^2(M, \mathbb{C}) \to H_A^2(M)$ is surjective. Then for all $d$-closed form
	$$\Lambda \in \hbar \operatorname{F}_1 \Omega^2(M, \mathbb{C})[[\hbar]],$$
	there exists a Fedosov connection $D^\gamma$ on $W$ and a $d_A$-closed element $\tau \in \hbar^2 \, \Omega_A^{0, 2}(M)[[\hbar]]$ such that
	\begin{equation}
		R^\gamma = -\Omega + \Lambda + \Lambda_\tau.
	\end{equation}
	If furthermore $D^{\gamma'}$ is a Fedosov connection on $W$ and $\tau' \in \hbar^2 \, \Omega_A^{0, 2}(M)[[\hbar]]$ is a $d_A$-closed element such that $R^{\gamma'} = -\Omega + \Lambda + \Lambda_{\tau'}$, then $[\tau] = [\tau']$.
\end{theorem}
\begin{proof}
	We begin with setting $\gamma_2 = I_2 = -\tilde{\Omega} \in \operatorname{F}_2\Omega^1(M, W)$ and noting that
	\[
	R^{\gamma_2} + \Omega - \Lambda = R^\nabla - \Lambda \in \operatorname{F}_3 \Omega^2(M, W).
	\]
	Next, assume that, for some $r \geq 2$, we can choose $I_2 \in \operatorname{F}_2\Omega^1(M, W)$, ..., $I_r \in \operatorname{F}_r\Omega^1(M, W)$ and a $d_A$-closed element $\tau \in \hbar^2 \Omega_A^{0, 2}(M) \oplus \hbar^3 \Omega_A^{0, 2}(M) \oplus \cdots \oplus \hbar^{\lfloor \frac{r}{2} \rfloor} \Omega_A^{0, 2}(M)$, such that
	\[
	R^{\gamma_r} + \Omega - \Lambda - \Lambda_\tau \in \operatorname{F}_{r+1} \Omega^2(M, W),
	\]
	where $\gamma_r := I_2 + \cdots + I_r$. Note that the lift $\Lambda_\tau$ of $\tau$ is guaranteed by the surjectivity of $\pi_{0, *}: H^2(M, \mathbb{C}) \to H_A^2(M)$. This assumption holds for $r = 2$ by setting $\tau = 0$.
	
	By the Bianchi identity for $D^{\gamma_r}$ and the fact that $\omega - \Lambda - \Lambda_\tau$ is central and $d$-closed, we obtain
	\[
	D^{\gamma_r} (R^{\gamma_r} + \Omega - \Lambda - \Lambda_\tau) = D^{\gamma_r} (R^{\gamma_r}) + d(\Omega - \Lambda - \Lambda_\tau) = 0.
	\]
	Projecting this equation onto the associated graded quotient $\operatorname{Gr}_{r+1} \Omega^*(M, W)$, the leading-order differential operator reduces to exactly $d_A^{\operatorname{Bott}} - \delta$. This implies
	$$(d_A^{\operatorname{Bott}} - \delta) [R^{\gamma_r} + \Omega - \Lambda - \Lambda_\tau]_{r+1} = 0.$$
	Here, for any $\alpha \in \operatorname{F}_{r+1}\Omega^*(M, W)$, we define $[\alpha]_{r+1} := \alpha + \operatorname{F}_{r+2} \Omega^*(M, W) \in \operatorname{Gr}_{r+1} \Omega^*(M, W)$.
	
	By Remark \ref{Remark: application of lemma}, there exist a $1$-form $I_{r+1} \in \operatorname{F}_{r+1} \Omega^1(M, W)$ and an element $\tau_{r+1} \in \Omega_A^{0, 2}(M)[[\hbar]]$ such that
	\[
	[R^{\gamma_r} + \Omega - \Lambda - \Lambda_\tau]_{r+1} + (d_A^{\operatorname{Bott}} - \delta)[I_{r+1}]_{r+1} = \tau_{r+1}.
	\]
	Moreover, $\tau_{r+1}=0$ whenever $r+1$ is odd. If $r+1$ is even, then $\tau_{r+1} = \hbar^{\frac{r+1}{2}} \widetilde{\tau}$, where $\widetilde{\tau} \in \Omega_A^{0, 2}(M)$ is $d_A$-closed. Clearly,
	$$\tau' := \tau + \tau_{r+1} \in \hbar^2 \Omega_A^{0, 2}(M) \oplus \hbar^3 \Omega_A^{0, 2}(M) \oplus \cdots \oplus \hbar^{\lfloor \frac{r+1}{2} \rfloor} \Omega_A^{0, 2}(M)$$
	is $d_A$-closed and $\Lambda_{\tau'} = \Lambda_\tau + \Lambda_{\tau_{r+1}}$. It further follows that
	\[
	R^{\gamma_r} + D^{\gamma_r}I_{r+1} + \omega - \Lambda - \Lambda_{\tau'} \in \operatorname{F}_{r+2} \Omega^2(M, W).
	\]
	Note that, setting $\gamma_{r+1} := \gamma_r + I_{r+1}$, we have
	\[
	R^{\gamma_{r+1}} = R^{\gamma_r} + D^{\gamma_r} I_{r+1} + \frac{1}{2\hbar} [I_{r+1}, I_{r+1}]_*.
	\]
	Since $I_{r+1} \in \operatorname{F}_{r+1}\Omega^1(M, W)$ and $r \geq 2$, the commutator $\frac{1}{2\hbar} [I_{r+1}, I_{r+1}]_*$ lies in $\operatorname{F}_{r+2}\Omega^2(M, W)$. It forces that
	\[
	R^{\gamma_{r+1}} + \Omega - \Lambda - \Lambda_{\tau'} \in \operatorname{F}_{r+2} \Omega^2(M, W).
	\]
	
	Proceeding inductively degree by degree, we construct a formal sequence of $1$-forms $\{\gamma_r\}_{r\geq 2}$. By the completeness of the filtration, this sequence converges to a $1$-form $\gamma \in \Omega^1(M, W)$ satisfying the conditions that $\gamma + \tilde{\Omega} \in \operatorname{F}_3\Omega^1(M, W)$ and that $R^\gamma = -\Omega + \Lambda + \Lambda_\tau$ for some $d_A$-closed element $\tau \in \hbar^2 \, \Omega_A^{0, 2}(M)[[\hbar]]$.
	
	Suppose $D^{\gamma'}$ is another Fedosov connection on $W$ and $\tau' \in \hbar^2\Omega_A^{0,2}(M)[[\hbar]]$ is a $d_A$-closed form such that $R^{\gamma'} = -\Omega + \Lambda + \Lambda_{\tau'}$. Setting $\eta := \gamma' - \gamma \in F_3\Omega^1(M,W)$, the curvature variation formula yields
	\begin{equation}\label{Curvaturevariation}
		R^{\gamma'} - R^\gamma = D^\gamma\eta + \frac{1}{2\hbar}[\eta,\eta]_* = \Lambda_{\tau'} - \Lambda_\tau.
	\end{equation}
	
	Write $\tau - \tau' = \sum_{r \ge 2} \hbar^r \mu_r$, where each $\mu_r \in \Omega_A^{0,2}(M)$ is $d_A$-closed, and decompose $\eta = \sum_{N \ge 3} \eta_{[N]}$ into its weight components. Let $d_0 := d_A^{\mathrm{Bott}} - \delta$, and let $\pi_{0,*}$ denote the projection from the associated graded complex to the $A$-complex. By Remark~\ref{Remark: application of lemma}, we have the chain-level identity $\pi_{0,*} \circ d_0 = d_A \circ \pi_{0,*}$.
	
	Let $N_0$ be the lowest Fedosov weight for which $\eta_{[N_0]} \neq 0$. Since $N_0 \ge 3$, the quadratic term $\frac{1}{2\hbar}[\eta,\eta]_*$ has weight at least $2N_0 - 2 > N_0$, and thus does not contribute to the weight-$N_0$ component. Moreover, the leading-order term of $D^\gamma$ on the associated graded complex is precisely $d_0$. Extracting the weight-$N_0$ component of \eqref{Curvaturevariation} therefore gives $(\Lambda_{\tau'} - \Lambda_\tau)_{[N_0]} = d_0\eta_{[N_0]}$. Applying $\pi_{0,*}$, the corresponding $A$-complex component becomes $d_A \pi_{0,*}(\eta_{[N_0]})$. Thus, the first non-vanishing coefficient of $\tau' - \tau$ is necessarily $d_A$-exact.
	
	We now proceed by induction on the power of $\hbar$. Assume that $[\mu_2] = \cdots = [\mu_{k-1}] = 0$ in $H_A^2(M)$. Evaluating the identity \eqref{Curvaturevariation} at weight $2k$ yields
	\[
	(\Lambda_{\tau'} - \Lambda_\tau)_{[2k]} = d_0\eta_{[2k]} + \Psi_{2k},
	\]
	where $\Psi_{2k}$ depends only on the lower-weight components of $\eta$ and of the connection. By the induction hypothesis and the fact that $d_0$ acts as a derivation of the graded commutator, $\pi_{0,*}(\Psi_{2k})$ is $d_A$-exact. Applying $\pi_{0,*}$ to the weight-$2k$ equation then shows that $-\mu_k$ has same cohomology as $d_A \pi_{0,*}(\eta_{[2k]})$. Hence, $\mu_k$ is also $d_A$-exact.
	
	By induction, every coefficient $\mu_r$ is $d_A$-exact. Therefore, $\tau - \tau'$ is $d_A$-exact in $\hbar^2\Omega_A^{0,2}(M)[[\hbar]]$, which implies $[\tau] = [\tau']$ in $\hbar^2 H_A^2(M)[[\hbar]]$. We conclude that the cohomology class $[\tau] \in \hbar^2 H_A^2(M)[[\hbar]]$ is uniquely determined by $[\Lambda] \in \hbar H^2(F_1\Omega^*(M,\mathbb C), d)[[\hbar]]$.
\end{proof}

In general, the non-commutative period map $\operatorname{per}: \operatorname{DQ}(M, \omega) \to -[\Omega] + \hbar H^2(M, \mathbb{C})[[\hbar]]$ might be non-injective. As we shall see, a sufficient condition for its injectivity is the surjectivity of the map $\pi_{0, *}: H^1(M, \mathbb{C}) \to H_A^1(M)$. Motivated by this observation and the assumptions of Theorem \ref{Theorem 3.18}, we make the following definition.

\begin{definition}\label{admissible}
	A transverse holomorphic symplectic manifold $(M, \Omega)$ is said to be \emph{admissible} if the map $\pi_{0, *}: H^i(M, \mathbb{C}) \to H_A^i(M)$ is surjective for $i \in \{1, 2\}$.
\end{definition}

The admissibility condition ensures the following short exact sequence of vector spaces:
$$ 0 \longrightarrow H^2(F_1 \Omega^*(M, \mathbb{C}), d) \longrightarrow H^2(M, \mathbb{C}) \longrightarrow H_A^2(M) \longrightarrow 0,$$
where $H^2(M, \mathbb{C}) \to H_A^2(M)$ is induced by the restriction map $\pi_{0, *}: \Omega^*(M, \mathbb{C})\rightarrow \Omega_A^{0, *}(M)$. When $(M, \Omega)$ is admissible, the fixed $\mathbb{C}$-linear right inverse $H_A^2(M) \to H^2(M, \mathbb{C})$ at the beginning of this subsection becomes a splitting of the above short exact sequence:
\begin{equation}
	\label{Equation: splitting}
	H^2(M,\mathbb{C}) \cong H^2(F_1\Omega^*(M,\mathbb{C}),d) \oplus H_A^2(M).
\end{equation}

\begin{theorem}
	\label{Theorem: classification}
	Let $(M, \Omega)$ be an admissible transverse holomorphic symplectic manifold. Then the equivalence classes of deformation quantizations of $(M, \Omega)$ are in bijective correspondence with the elements in $\hbar \, H^2(\operatorname{F}_1 \Omega^*(M, \mathbb{C}), d)[[\hbar]]$.
\end{theorem}
\begin{proof}
	Recall from Theorem \ref{dqtofedosov} that every equivalence class of deformation quantizations of $(M, \Omega)$ admits a representative of the form $(\mathcal{W}_{D^\gamma}, \psi_\gamma)$ for some Fedosov connection $D^\gamma$ on $W$. For each Fedosov connection $D^\gamma$, denote by $[(\mathcal{W}_{D^\gamma}, \psi_\gamma)]$ the equivalence class of $(\mathcal{W}_{D^\gamma}, \psi_\gamma)$ and define $\widetilde{\operatorname{per}}([(\mathcal{W}_{D^\gamma}, \psi_\gamma)])$ as the projection of $[R^\gamma + \Omega]$ onto $\hbar H^2(\operatorname{F}_1\Omega^*(M,\mathbb C),d)[[\hbar]]$ with respect to the fixed splitting \eqref{Equation: splitting}. It suffices to show that the map
	\[
	\widetilde{\operatorname{per}}: \operatorname{DQ}(M, \Omega) \to \hbar H^2(\operatorname{F}_1\Omega^*(M,\mathbb C),d)[[\hbar]]
	\]
	is bijective. The surjectivity of $\widetilde{\operatorname{per}}$ directly follows from Theorem \ref{Theorem 3.18}.
	
	For the injectivity, suppose that $D^\gamma, D^{\gamma'}$ are Fedosov connections on $W$ such that
	\begin{equation}
		[\Lambda] := \widetilde{\operatorname{per}}([(\mathcal{W}_{D^\gamma}, \psi_\gamma)]) = \widetilde{\operatorname{per}}([(\mathcal{W}_{D^{\gamma'}}, \psi_{\gamma'})]).
	\end{equation}
	Due to the fixed splitting \eqref{Equation: splitting}, we can choose classes $[\tau], [\tau'] \in \hbar H_A^2(M)[[\hbar]]$ such that
	\[
	R^\gamma = -\Omega + \Lambda + \Lambda_\tau \quad \text{and} \quad R^{\gamma'} = -\Omega + \Lambda + \Lambda_{\tau'}.
	\]
	Recall that $R^\nabla \in \Omega^2(M, \operatorname{Sym}^2B^*)$ satisfies $\pi_{0, *}(R^\nabla) = 0$, whence $R^\nabla \in \operatorname{F}_3\Omega^2(M, W)$. It then follows from $I := \gamma + \tilde{\Omega} \in \operatorname{F}_3\Omega^1(M, W)$ that $R^\gamma + \Omega = R^\nabla + (\nabla^W - \delta)I + \tfrac{1}{2\hbar}[I, I]_* \in \operatorname{F}_3 \Omega^2(M, W)$. Thus,
	\[
	\tau = \pi_{0, *}(R^\gamma) \in \hbar \Omega_A^{0, 2}(M)[[\hbar]] \cap \operatorname{F}_3\Omega^2(M, W) = \hbar^2 \, \Omega_A^{0, 2}(M)[[\hbar]].
	\]
	Similarly, $\tau' \in \hbar^2 \, \Omega_A^{0, 2}(M)[[\hbar]]$. By Theorem \ref{Theorem 3.18} again, we obtain $[\tau] = [\tau']$ and thus $[R^\gamma] = [R^{\gamma'}]$. By the admissibility assumption together with Lemma \ref{samecurvature} and Theorem \ref{Theorem: gauge equivalence}, $(\mathcal{W}_{D^\gamma}, \psi_\gamma)$ and $(\mathcal{W}_{D^{\gamma'}}, \psi_{\gamma'})$ are equivalent. This shows the injectivity of $\widetilde{\operatorname{per}}$, concluding the proof.    
\end{proof}

\section{Deformation quantization and coisotropic A-branes}
\label{Section: Examples}
In this section, we relate our deformation-quantization framework
to coisotropic A-branes and examine its applicability through
examples. In Subsection \ref{Subsection: A-branes and THS}, we recall how coisotropic A-branes
carry transverse holomorphic symplectic structures and establish
a converse realization. Subsections \ref{Subsection: eg1}--\ref{Subsection: eg4} study examples arising
from mapping tori, $T$-duality, and coisotropic lifts in cotangent bundles, focusing on admissibility and the transverse Atiyah class. These examples illustrate the distinction between deformation quantizations and global star products, as well as the separate hypotheses of our existence and classification theorems. Finally, building on Subsection \ref{Subsection: eg4}, Subsection \ref{section5.6} constructs, under suitable compatibility assumptions, a deformation quantization together with a compatible module structure. This construction is expected to model, at the perturbative level, the A-model composition
\[
\operatorname{Hom}(\mathcal B,\mathcal B)
\times
\operatorname{Hom}(\mathcal L,\mathcal B)
\longrightarrow
\operatorname{Hom}(\mathcal L,\mathcal B)
\]
for a coisotropic A-brane $\mathcal{B}$ and a Lagrangian A-brane $\mathcal{L}$ whose support is contained in that of $\mathcal{B}$.

\subsection{Coisotropic A-branes and transverse holomorphic symplectic manifolds}
\quad\par
\label{Subsection: A-branes and THS}

We recall the definition of (rank-$1$) coisotropic A-branes, including the presence of a B-field.

\begin{definition}
	Let $(X, \omega)$ be a symplectic manifold and $B \in \Omega^2(X)$ be a closed real-valued $2$-form. A \emph{coisotropic A-brane} on $(X, \omega, B)$ consists of a coisotropic submanifold $M \subset X$ and a Hermitian line bundle $L$ over $M$ with a unitary connection $\nabla^L$ satisfying the following conditions:
	\begin{itemize}
		\item Denote the curvature of $\nabla^L$ by $-2\pi\sqrt{-1} F_L$, and let $A_\mathbb{R} = \ker (\omega \vert_M)$ be the characteristic distribution of $M$. Then $F_L - B \vert_M$ descends to a non-degenerate $2$-form $\underline{F}$ on $TM/ A_\mathbb{R}$.
		\item The endomorphism $\underline{\omega}^{-1}\circ \underline{F}: TM/ A_\mathbb{R} \to TM/ A_\mathbb{R}$ defines a complex structure on $TM/ A_\mathbb{R}$, where $\underline{\omega}$ is the non-degenerate $2$-form on $TM / A_\mathbb{R}$ descended from $\omega \vert_M$.
	\end{itemize}
\end{definition}

Consequently, for every coisotropic A-brane $\mathcal{B} = (M, L, \nabla^L)$ on $(X, \omega, B)$, $M$ carries a transverse holomorphic structure for which $\Omega = (F_L - B \vert_M) + \sqrt{-1}\omega$ defines a transverse holomorphic symplectic form (cf. \cite{KapOrl2003}). We now record a converse.

\begin{proposition}
	\label{Proposition: THS as coisotropic A-branes}
	Let $(M, \Omega)$ be a transverse holomorphic symplectic manifold, and let $L$ be a Hermitian line bundle over $M$ equipped with a unitary connection $\nabla^L$. Then there exist a symplectic manifold $(X, \omega)$, a closed real $2$-form $B \in \Omega^2(X)$, and a closed embedding $M \hookrightarrow X$ such that
	\[
	\mathcal{B} = (M, L, \nabla^L)
	\]
	is a coisotropic A-brane on $(X, \omega, B)$. Moreover, if the curvature of $\nabla^L$ is $-2\pi\sqrt{-1}\operatorname{Re} \Omega$, then one may take $B = 0$.
\end{proposition}
\begin{proof}
	The imaginary part $\operatorname{Im} \Omega$ of $\Omega$ defines a presymplectic form on $M$. Applying Gotay's coisotropic embedding theorem \cite{coisoimbedd} yields a vector bundle over $M$ together with a symplectic form $\omega$ on its total space $X$ such that $\omega \vert_M = \operatorname{Im}\Omega$, where $M$ is embedded into $X$ via the zero section. Then $B = \pi^*(F_L - \operatorname{Re}\Omega)$ defines a closed real $2$-form on $X$. Consequently, $\mathcal{B}$ is a coisotropic A-brane on $(X, \omega, B)$.
\end{proof}

\subsection{Example I: transverse holomorphic symplectic mapping tori}
\quad\par
\label{Subsection: eg1}

Mapping-torus constructions link symplectic dynamics with the categorical structures of homological mirror symmetry. Seidel \cite{Sei2014} uses Fukaya categories to study symplectic mapping tori and their flux-type
invariants, while Kartal \cite{Kar2021} describes wrapped Fukaya categories of open symplectic mapping tori by categorical constructions incorporating the monodromy. In a complementary direction, Qin--Wang \cite{QinWan2018} use a mapping-torus construction based on Kummer surface automorphisms to construct compact non-K\"ahler manifolds carrying
both complex and symplectic Calabi--Yau structures.

We use this construction to obtain admissible coisotropic A-branes on compact symplectic Calabi--Yau manifolds. In the
$\mathbb T^4$ example below, every characteristic leaf is nonclosed and dense in the brane support; in the Kummer example, this holds for almost every leaf. The latter also has nonvanishing transverse
Atiyah class, and its ambient manifold admits no K\"ahler structure.

Let $Y$ be a compact K\"ahler holomorphic symplectic manifold and let $\phi:Y\to Y$ be a biholomorphism preserving the holomorphic symplectic form $\Omega_Y$. Its mapping torus is
\begin{equation}
	M:=(\mathbb R\times Y)/\bigl((t+1,y)\sim(t,\phi(y))\bigr).
\end{equation}
Write $\pi_Y: \mathbb{R} \times Y \to Y$ for the projection. The $2$-form $\pi_Y^*\Omega_Y$ descends to a transverse holomorphic symplectic form $\Omega$ on $M$. We denote by $A$ the underlying transverse holomorphic structure. The characteristic distribution is generated by the suspension flow of $\phi$. A leaf is dense in $M$ if and only if its corresponding $\phi$-orbit is dense in $Y$.

The standard symplectic mapping torus construction (cf. \cite[Section 4.1]{Sei2014}) equips
\[
X:=\mathbb S^1 \times M,
\quad
\omega :=dt\wedge ds + \pi_M^*\operatorname{Im}\Omega,
\]
with a compact symplectic structure. Here $\pi_M: X \to M$ is the projection, $dt$ comes from the base
circle of $M$, and $ds$ from the additional circle. The inclusion $M\cong M\times\{1\}\subset X$ is coisotropic, with $\omega|_M=\operatorname{Im}\Omega$.

Since $\phi$ is biholomorphic, the product complex structure determined by the complex structure on $Y$ and the coordinate $t+\sqrt{-1}s$ descends to a complex structure $J$ on $X$. The arguments of \cite[Propositions 3.3 and 3.4]{QinWan2018} apply in this setting, showing that the canonical bundle of
$(X,J)$ is holomorphically trivial and that $c_1(X,\omega)=0$. Hence $X$ carries both complex and symplectic Calabi--Yau structures in the sense of \cite{QinWan2018}.

\begin{proposition}\label{prop:mapping-torus}
	For the above mapping torus:
	\begin{enumerate}
		\item The transverse holomorphic symplectic manifold $(M, \Omega)$ is admissible.
		\item If $Y$ has nonzero Atiyah class, then so does $(M, A)$.
	\end{enumerate}
\end{proposition}
\begin{proof}
	We prove that the restriction map
	$H^i(M,\mathbb C)\to H_A^i(M)$ is surjective for every $i$. Both ordinary forms and $A$-forms on $M$ identify with $\mathbb Z$-invariant forms on $\mathbb R\times Y$. We work upstairs, using $A$ also to denote the lifted distribution
	$(d\pi_Y)^{-1}(T^{0,1}Y)$.
	
	Since $Y$ is compact K\"ahler, each Dolbeault class in $H_{\bar\partial}^{0,q}(Y)$ has a unique antiholomorphic representative. These representatives are $d$-closed, and their choice commutes with $\phi^*$.
	
	Let $\gamma \in \Omega_A^{0, i}(\mathbb{R} \times Y)$ be $d_A$-closed and invariant. Write
	\[
	\gamma=\alpha(t)+dt\wedge\beta(t), \quad \alpha(t) \in \Omega^{0, i}(Y), \quad \beta(t) \in \Omega^{0, i-1}(Y).
	\]
	The cocycle equations are $\overline{\partial} \alpha = 0$ and $\partial_t\alpha = \overline{\partial} \beta$. Hence $[\alpha(t)]_{\bar\partial}$ is independent of $t$, and invariance of $\gamma$ implies that its antiholomorphic representative $a$ is $\phi^*$-invariant.
	
	Choose a smooth invariant family $\eta(t)$ satisfying $\alpha(t) = a + \bar\partial\eta(t)$. Such primitives exist locally
	over the base circle by Dolbeault Hodge theory and can be patched
	by a partition of unity on that circle. Subtracting $d_A\eta$ gives
	\[
	\gamma-d_A\eta
	=a+dt\wedge\beta'(t),
	\qquad
	\beta':=\beta-\partial_t\eta,
	\qquad
	\bar\partial\beta'=0.
	\]
	Let $b(t)$ be the antiholomorphic representative of $[\beta'(t)]$.
	It is smooth and invariant. Choosing invariant primitives
	$\beta'-b=\bar\partial\xi$ in the same way, we obtain
	\[
	\gamma-d_A(\eta-dt\wedge\xi)
	=a+dt\wedge b(t).
	\]
	The right-hand side is an invariant ordinary closed form. It descends to a closed form on $M$ whose restriction represents $[\gamma]$, proving surjectivity and hence (1).
	
	For (2), take the embedding $\iota:Y\hookrightarrow M$, $y\mapsto[0,y]$. The bundle $B := T_{\mathbb C}M / A$ restricts along $\iota$ to
	$T^{1,0}Y$. Since $d\iota(T^{0,1}Y)\subset A$, pullback by $\iota$ sends $A$-compatible connections on $B$ to holomorphic connections on $T^{1,0}Y$. This proves the assertion.
\end{proof}

Let $E := \mathbb{C}^2/(\mathbb{Z} \oplus \sqrt{-1}\mathbb{Z})^2$, and choose $a_1,\ldots,a_4\in\mathbb R$ such that
$1,a_1,\ldots,a_4$ are linearly independent over $\mathbb Q$. Take $Y=E$ with its standard holomorphic symplectic form $\Omega_Y=dz^1\wedge dz^2$, and let $\phi$ be the translation $(z^1, z^2) \mapsto (z^1 + a_1 + \sqrt{-1}a_2, z^2 + a_3 + \sqrt{-1}a_4)$. The resulting mapping torus $M$ is diffeomorphic to $\mathbb{T}^5 = \mathbb R^5/\mathbb{Z}^5$ and is admissible by Proposition \ref{prop:mapping-torus}. Every characteristic leaf is dense, and the translation-invariant flat connection descends to an $A$-compatible connection, so the Atiyah class of $(M,A)$ vanishes.

For an example with nonvanishing Atiyah class, we follow \cite[Section 2]{QinWan2018}.

\begin{theorem}
	There exists a compact complex and symplectic Calabi--Yau manifold $(X, \omega)$ together with a rank-$1$ coisotropic A-brane $\mathcal{B}$ on it such that
	\begin{enumerate}
		\item $X$ admits no K\"ahler structure;
		\item The underlying transverse holomorphic symplectic manifold $(M, \Omega)$ of $\mathcal{B}$ is admissible and has nonzero Atiyah class;
		\item Almost every leaf of the characteristic distribution of $M$ is dense.
	\end{enumerate}
\end{theorem}
\begin{proof}
	Let $Y = \operatorname{Km}(E)$ be the Kummer surface associated with $E$ defined above, and let $\phi: Y \to Y$ be the biholomorphism induced by
	$\begin{pmatrix}
		2 & 1\\
		1 & 1
	\end{pmatrix} \in \operatorname{SL}(2, \mathbb{Z}[\sqrt{-1}])$. Writing $\operatorname{Bl}: \widetilde{E} \to E$ for the blow-up at the sixteen fixed points of the involution $z \mapsto -z$, the lifted involution defines a branched double cover $q:\widetilde E\to Y$, and $Y$ carries a $\phi$-invariant holomorphic symplectic form $\Omega_Y$ determined by $q^*\Omega_Y=\operatorname{Bl}^*(2 dz^1\wedge dz^2)$. Its pullback to $\mathbb R\times Y$ descends to the transverse holomorphic symplectic form $\Omega$ on the mapping torus $M$. Since the Atiyah class of $Y$ is nonzero, Proposition \ref{prop:mapping-torus} shows that $(M,\Omega)$ is admissible and has nonvanishing Atiyah class.
	
	Since no eigenvalue of $\begin{pmatrix}
		2 & 1\\
		1 & 1
	\end{pmatrix}$ is a root of unity, the induced toral automorphism of $E$ is ergodic with respect to normalized Haar measure. Ergodicity passes to the quotient and its resolution, since the exceptional locus has measure zero. Thus $\phi$ is
	ergodic with respect to the normalized holomorphic symplectic volume measure on $Y$. This measure has full support, so almost
	every $\phi$-orbit, and hence almost every characteristic leaf
	of the mapping torus, is dense.

	Since $[\operatorname{Re}(dz^1\wedge dz^2)]$ is integral, the identity $q_!q^*=2\operatorname{Id}$ gives
	\[
	[\operatorname{Re}\Omega_Y]=q_!\operatorname{Bl}^*[\operatorname{Re}(dz^1\wedge dz^2)]
	\in H^2(Y,\mathbb Z),
	\]
	where $q_!$ denotes the integral cohomological pushforward. Thus, there exists a Hermitian line bundle $L_Y$ with a unitary connection of curvature
	$-2\pi\sqrt{-1}\operatorname{Re}\Omega_Y$. The induced connection on
	$\phi^*L_Y\otimes L_Y^{-1}$ is flat and, since $Y$ is simply connected, admits a global parallel
	unit section. Hence, $\phi$ admits a unitary,
	connection-preserving lift, whose suspension defines a Hermitian line bundle on $M$ with curvature $-2\pi\sqrt{-1}\operatorname{Re}\Omega$.
	
	Consequently, $M$ supports a rank-$1$ coisotropic A-brane on the compact symplectic Calabi--Yau six-manifold $X=M\times\mathbb S^1$ with zero
	$B$-field. This ambient manifold belongs to the family of \cite[Theorem 1.1]{QinWan2018} and admits no K\"ahler structure.
\end{proof}

\subsection{Example II: coisotropic A-branes and $T$-duality}
\quad\par
\label{Subsection: eg2}

The Strominger--Yau--Zaslow proposal \cite{StrYauZas1996} interprets mirror symmetry as fibrewise $T$-duality for special Lagrangian torus fibrations. In the semi-flat setting, Fourier--Mukai-type transforms relate Lagrangian A-branes to B-branes \cite{LeuYauZas2000}. Beyond the Lagrangian case, Aldi--Zaslow \cite{AldZas2005} studied explicit coisotropic A-branes on flat four-tori, exhibiting noncommutative algebras acting on spaces of open-string ground states.

The first author, with K.-L. Chan and Zhang \cite{ChaLeuZha2018}, extended this correspondence to \emph{semi-affine coisotropic A-branes}. The first two authors, with K. Chan, Li and Suen \cite{ChaLeuLiSueYau2026}, subsequently studied rank-$1$ space-filling semi-affine coisotropic A-branes on $X=T^*Q/T_{\mathbb Z}^*Q$, for an integral affine manifold $Q$, under a fibrewise non-degeneracy assumption equivalent to the mirror B-brane being space-filling. The induced holomorphic symplectic structure admits Darboux charts with affine symplectic transition maps. Thus the standard flat connections glue, giving vanishing of the Atiyah class of $X$, and the local Moyal products define a global star product on $\mathcal O_X[[\hbar]]$. Neither vanishing of the transverse Atiyah class nor our admissibility condition should be expected to hold automatically for general semi-affine coisotropic A-branes.

Further $T$-duality examples on nilmanifolds appear in Katzarkov--Soriani \cite{KatSor2018}. See also Sis\'ak \cite{Sis2024} for $T$-duality of hyperK\"ahler branes, and Qin's thesis \cite{Qin2020} for a twisted doubling approach to coisotropic A-branes on symplectic tori.

\subsection{Example III: coisotropic lifts of complex submanifolds to $T^*Q$}
\quad\par
\label{Subsection: eg4}

Cotangent bundles connect A-branes with sheaf theory and quantization: the Nadler--Zaslow correspondence \cite{NadZas2009} relates Lagrangian A-branes to constructible sheaves, while Kuwagaki's $\hbar$-Riemann--Hilbert correspondence \cite{Kuw2022} links deformation quantization and sheaf quantization on holomorphic cotangent bundles.

For coisotropic A-branes on cotangent bundles, Bischoff--Gualtieri \cite{BisGua2022} use twisted cotangent bundles to construct space-filling examples. Combining this with the conormal construction of Lagrangian branes, the first two authors with Xie \cite{qcr} construct coisotropic lifts of complex submanifolds in their study of quantization and reduction. We show that Atiyah-class vanishing and admissibility for these lifts reduce to the corresponding conditions on the associated twisted cotangent bundles. We begin by briefly recalling the construction.

Let $Q$ be a smooth manifold and $Z\subset Q$ a smooth submanifold
equipped with a complex structure and a Hermitian holomorphic line bundle
$L_Z$ with Chern connection. The associated twisted cotangent bundle
$(Y,\Omega_Y)$ is smoothly identified with the real cotangent bundle $T^*Z$.
Restriction of covectors defines a projection
\[
\pi:M:=T^*Q|_Z\longrightarrow Y.
\]
The construction in \cite{qcr} equips $M$ with a rank-$1$ coisotropic
A-brane on $X = T^*Q$, inducing a transverse holomorphic symplectic structure $(\Omega, A)= (\pi^*\Omega_Y, (d\pi)^{-1}(T^{0,1}Y))$ on $M$. The short exact sequence of real vector bundles
\[
0\longrightarrow N^*_{Z/Q}
\longrightarrow T^*Q|_Z
\longrightarrow T^*Z
\longrightarrow0
\]
shows that $\pi$ is a smooth real affine bundle modelled on the conormal bundle $N^*_{Z/Q}$ of $Z \subset Q$. The following propositions make precise the reduction of the conditions of interest to $Y$.

\begin{proposition}\label{prop:Atiyah-descent}
	Let $\pi:M\to Y$ be a smooth surjective submersion with connected fibres over a complex manifold $Y$. Equip $M$ with the transverse holomorphic structure $A:=(d\pi)^{-1}(T^{0,1}Y)$. Then the Atiyah class of $(M,A)$ vanishes if and only if the Atiyah class of $Y$ vanishes.
\end{proposition}
\begin{proof}
	Identify
	$B :=T_{\mathbb C}M/A$ with $\pi^*T^{1,0}Y$. Pullback by $\pi$, after extension by $\bar\partial$, sends holomorphic
	connections on $T^{1,0}Y$ to $A$-compatible connections on $B$. 
	
	Conversely, let $\nabla$ be an $A$-compatible connection on $B$. The Bott condition and connectedness of the fibres imply that for any open subset $U \subset Y$, a smooth section $b \in \Gamma(\pi^{-1}(U), B)$ is pulled back from a unique
	smooth section of $T^{1,0}Y$ over $U$ if and only if $\nabla_V b = 0$ for all vertical vector field $V$ on $\pi^{-1}(U)$. For a local smooth section $s$ of $T^{1,0}Y$ and a projectable lift $\widetilde v$ of a local complex vector field $v$ on $Y$, $A$-compatibility gives
	\[
	\nabla_V\nabla_{\widetilde v}\pi^*s=\nabla_{\tilde{v}} \nabla_V \pi^*s + \nabla_{[V, \tilde{v}]} \pi^*s = 0,
	\]
	for every vertical vector field $V$, since $[V,\widetilde v]$ is vertical. Hence,
	\[
	\pi^*(D_vs)=\nabla_{\widetilde v}\pi^*s
	\]
	defines a connection $D$ on $T^{1,0}Y$, independently of the chosen lift, with $\nabla=\pi^*D$. $A$-compatibility implies that $D^{1,0}$ is a holomorphic connection, proving the assertion.
\end{proof}

\begin{proposition}
	Let $\pi:M\to Y$ be a smooth real affine bundle over a holomorphic symplectic manifold $(Y,\Omega_Y)$. Equip the total space $M$ with the transverse holomorphic symplectic structure
	\[
	\Omega=\pi^*\Omega_Y, \quad A=(d\pi)^{-1}(T^{0,1}Y).
	\]
	Then $(M,\Omega)$ is admissible if and only if $(Y,\Omega_Y)$ is admissible.
\end{proposition}
\begin{proof}
	The inclusion $d\pi(A) \subset T^{0, 1}Y$ induces a cochain map
	\begin{equation*}
		\pi_A^*: (\Omega^{0, *}(Y), \overline{\partial}) \to (\Omega_A^{0, *}(M), d_A)
	\end{equation*}
	compatible with ordinary pullback $\pi^*: \Omega^*(Y, \mathbb{C}) \to \Omega^*(M, \mathbb{C})$ through the restriction maps. Since ordinary pullback induces an isomorphism on de Rham cohomology, the equivalence of admissibility follows once $\pi_A^*$ is shown to be a quasi-isomorphism.
	
	Choose a smooth section $s:Y\to M$, which exists by a partition of unity. Since $ds(T^{0,1}Y)\subset A$, restriction along $s$ induces a cochain map $s_A^*:(\Omega_A^{0, *}(M), d_A)\to(\Omega^{0,*}(Y), \overline{\partial})$ satisfying $s_A^*\pi_A^*=\operatorname{Id}$. Consider the fibrewise contraction
	\[
	H_t(m) := s(\pi(m))+t\bigl(m-s(\pi(m))\bigr), \quad m \in M \text{ and } 0 \leq t \leq 1.
	\]
	We have $dH_t(A)\subset A$ and
	$\partial_tH_t\in \Gamma(M, A)$. Thus, for $\alpha\in\Omega_A^{0,i}(M)$ with $i\geq1$, the formula
	\[
	K(\alpha)(a_1, \dots, a_{i-1})=\int_0^1 \alpha (\partial_tH_t(m), dH_t(a_1), \dots, dH_t(a_{i-1})) dt, \quad m \in M, a_1, \dots, a_{i-1} \in A_m,
	\]
	defines an operator of degree $-1$, with $K=0$ in degree zero. The ordinary homotopy formula, restricted to $A$, gives $d_AK+Kd_A=\operatorname{Id}-\pi_A^*s_A^*$, implying that $\pi_A^*$ is a homotopy equivalence.
\end{proof}

\subsection{Example IV: coisotropic lifts of transverse complex manifolds to $T^*Q$}
\quad\par
\label{Subsection: eg3}
Another type of coisotropic lift construction in cotangent bundles appears in the work of Aganagic--Costello--McNamara--Vafa \cite{AgaCosMcnVaf2017}. In their three-dimensional setting, strings between Lagrangian and coisotropic branes supply the matter fields of Chern--Simons/matter theory, while transverse deformation quantization enters the proposed five-dimensional noncommutative Chern--Simons theory on the coisotropic branes. See also Aharony--Feldman--Honda \cite{AhaFelHon2019} for further developments on $T^*\mathbb S^3$, and Saulina \cite{Sau2014} for related
coisotropic branes on toric Calabi--Yau threefolds.

We recall the construction of \cite[Section 3.3]{AgaCosMcnVaf2017} in arbitrary dimension. Let $(Q,A')$ be a transverse complex manifold
with characteristic distribution $A'_{\mathbb R}$. The dual bundle
\[
M:=(TQ/A'_{\mathbb R})^*
\]
is naturally identified with the annihilator of $A'_{\mathbb R}$ in $T^*Q$, and hence is a coisotropic submanifold of $T^*Q$. Write $p:M\to Q$ for the projection and $\lambda$ for
the real tautological form on $T^*Q$, with $\omega=-d\lambda$. The transverse holomorphic structure on $Q$ induces the decomposition
\[
M \otimes \mathbb{C} = (B')^* \oplus \overline{(B')^*}, \quad \text{where } B' = (TQ \otimes \mathbb{C})/A'.
\]
For $\alpha\in M$,
write $\alpha^{1,0}$ for its component in $(B')^*$. Define
$\Theta\in\Omega^1(M,\mathbb C)$ by
\begin{equation}
	\label{Equation: primitive 1-form}
	\Theta(v)
	:=-2\sqrt{-1}\,\alpha^{1,0}(dp(v)),
	\quad \alpha \in M, v\in T_\alpha M\otimes\mathbb C.
\end{equation}
Then $\Omega := d\Theta$ is a transverse holomorphic symplectic form on $M$ with $\operatorname{Im}\Omega=\omega|_M$. Thus the trivial Hermitian line bundle with connection
$d-2\pi\sqrt{-1}\operatorname{Re}\Theta$ defines a rank-$1$ coisotropic
A-brane on $T^*Q$ with support $M$.

Now, we specialize to an example studied in \cite{AgaCosMcnVaf2017}: the Hopf fibration $h: Q = \mathbb S^3\to\mathbb P^1$ with transverse holomorphic structure $A'=(dh)^{-1}(T^{0,1}\mathbb P^1)$. The lift then has a
particularly concrete transverse geometry: writing $Y$
for the holomorphic cotangent bundle of $\mathbb{P}^1$ with canonical holomorphic symplectic form $\Omega_Y$, the preceding identification gives
\[
M \cong Y \times_{\mathbb P^1} \mathbb S^3
\]
under which $\Omega = \pi^*\Omega_Y$, where $\pi:M\to Y$ is the projection. Thus $M$ is a principal circle bundle over $Y$. The total space $Y$ of the holomorphic cotangent bundle of $\mathbb{P}^1$ has nonvanishing Atiyah class. Thus, Proposition \ref{prop:Atiyah-descent} implies that the Atiyah class of $(M,A)$ is also nonzero, where $A$ is the underlying transverse holomorphic structure of $(M, \Omega)$.

We conclude by examining admissibility, beginning with degree two.
Let $\theta\in\Omega^1(M)$ be the pulled-back Hopf connection and
put $\theta_A:=\theta|_A\in\Omega_A^{0,1}(M)$.
Its curvature is pulled back from a $(1,1)$-form on $Y$, so
$d_A\theta_A=0$. Cartan's homotopy formula shows that circle
averaging is chain homotopic to the identity on the $A$-complex; hence circle-invariant forms compute its cohomology. Every invariant $A$-form of degree $k$ is uniquely expressible as $\pi_A^*\alpha+\theta_A\wedge\pi_A^*\beta$, where
$\alpha\in\Omega^{0,i}(Y)$, $\beta\in\Omega^{0,i-1}(Y)$, and $\pi_A^*$ denotes the pullback induced by
$d\pi \vert_A:A\to T^{0,1}Y$. Moreover,
\[
d_A(\pi_A^*\alpha+\theta_A\wedge\pi_A^*\beta)
=
\pi_A^*\bar\partial\alpha
-\theta_A\wedge\pi_A^*\bar\partial\beta.
\]
The standard vanishings
$H_{\bar\partial}^{0,1}(Y)=H_{\bar\partial}^{0,2}(Y)=0$
therefore give
\[
H_A^2(M)\cong
H_{\bar\partial}^{0,2}(Y)
\oplus H_{\bar\partial}^{0,1}(Y)[\theta_A]=0.
\]
Thus the restriction map $H^2(M, \mathbb{C}) \to H_A^2(M)$ is surjective.

In degree one, since $M$ is a vector bundle over $\mathbb S^3$, $H^1(M,\mathbb C) \cong H^1(\mathbb{S}^3,\mathbb C)=0$. However, $[\theta_A]\ne0$ in $H_A^1(M)$: the circle fibres are tangent to $A$, and the integral of $\theta_A$ along each fibre is nonzero, whereas that of any $d_A$-exact form vanishes. Hence the restriction map $H^1(M, \mathbb{C}) \to H_A^1(M)$ is not surjective. Admissibility therefore fails only in degree one. Our existence theorem ensures that $(M,\Omega)$ admits deformation quantizations, although our classification theorem via the noncommutative period map does not apply.

\subsection{\textbf{Perturbative brane morphism from transverse complex geometry}\label{section5.6}}
\quad\par
We continue with a transverse complex manifold $(Q,A')$ and the associated quintuple
\begin{equation}
	\label{Equation: transverse cotangent bundle}
	(M, A, \Omega, p, j).
\end{equation}

Here $M=(TQ/A'_{\mathbb R})^*$, where $A'_{\mathbb R}$ is the characteristic distribution of $(Q, A')$; $A$ and $\Omega$ are the transverse holomorphic structure and transverse holomorphic symplectic form on $M$ constructed in Subsection \ref{Subsection: eg3}; and $p: M \to Q$ and $j: Q \to M$ denote the bundle projection and zero section.

For the pair of A-branes $\mathcal{B}$ and $\mathcal{L}$ supported on $M$ and on $Q$, respectively, the transverse holomorphic structure $A'$ on $Q = Q \cap M$ is adapted in the sense \cite[Definition 5.2]{qcr}. From the perspective of \cite{qcr}, and in the absence of instanton corrections, $H_{A'}^*(Q)$ is therefore a candidate model for the A-model morphism space $\operatorname{Hom}(\mathcal{L},\mathcal{B})$. Under composition, this space is expected to carry a natural left module structure over the A-model endomorphism algebra $\operatorname{Hom}(\mathcal{B},\mathcal{B})$.

In the present work, we formulate the perturbative counterpart of this picture at the sheaf level. Rather than working directly with the cohomological model $H_{A'}^*(Q) = H^*(Q, \mathcal{O}_{A'})$, we construct a deformation quantization \((\mathcal A_M^\hbar,\psi)\) of $(M, \Omega)$, together with a \(\mathcal A_M^\hbar\)-module structure on \(j_*\mathcal O_{A'}[[\hbar]]\) whose reduction modulo \(\hbar\) recovers the left \(\mathcal O_A\)-module structure on \(j_*\mathcal O_{A'}\) given by \(s\cdot f=j^*(s)f\). We expect this construction to model, at the perturbative level, the module structure arising from the A-model composition discussed above.

\begin{theorem}
	\label{Theorem: lift of flat connections}
	Let $(Q,A')$ be a transverse complex manifold and let $(M, A, \Omega, p, j)$ be given as in \eqref{Equation: transverse cotangent bundle}. Then there exist a deformation quantization \((\mathcal A_M^\hbar,\psi)\) of \((M,\Omega)\), satisfying
	\[
	\operatorname{per}([(\mathcal A_M^\hbar,\psi)]) = \hbar \pi \sqrt{-1} p^*(c_1(B')),
	\]
	and a left \(\mathcal A_M^\hbar\)-module structure $\rho_\hbar$ on \(j_*\mathcal O_{A'}[[\hbar]]\) such that the following diagram commutes:
	\begin{center}
		\begin{tikzcd}
			\mathcal A_M^\hbar
			\otimes_{\mathbb C[[\hbar]]}
			j_*\mathcal O_{A'}[[\hbar]]
			\ar[r, "\rho_\hbar"]
			\ar[d,"\psi\otimes j_*\pi_0"']
			&
			j_*\mathcal O_{A'}[[\hbar]]
			\ar[d,"j_*\pi_0"]
			\\
			\mathcal O_A
			\otimes_{\mathbb C}
			j_*\mathcal O_{A'}
			\ar[r, "\rho_{\operatorname{cl}}"']
			&
			j_*\mathcal O_{A'}
		\end{tikzcd}
	\end{center}
	where $\pi_0: \mathcal O_{A'}[[\hbar]] \to \mathcal O_{A'}$ is reduction modulo $\hbar$ and $\rho_{\operatorname{cl}}$ is given by $s \otimes f \mapsto j^*(s)f$.
\end{theorem}

We now develop the constructions needed to prove the theorem. For several local computations, it will be convenient to fix adapted coordinates. Henceforth, let
\[
(x^1,\ldots,x^m,z^1,\ldots,z^n)
\]
denote local coordinates on $Q$ adapted to the transverse holomorphic structure $A'$. They induce local transverse holomorphic Darboux coordinates
\[
(x^1, \dots, x^m, z^1, \dots, z^n, \xi_1, \dots, \xi_n)
\]
on $(M,\Omega)$, in which $\Theta=-\xi_idz^i$. We write $u^i$ and $v_i$ for the corresponding fiberwise symmetric variables associated with $dz^i$ and $d\xi_i$, respectively.

Recall that $\Omega$ descends to a non-degenerate $2$-form $\underline{\Omega}$ on \(B=T_{\mathbb C}M/A\). Let \(B':=T_{\mathbb C}Q/A'\). There is a natural short
exact sequence
\(0\to p^*(B')^*\to B\to p^*B'\to0\) with $p^*(B')^*$ embedded as a Lagrangian subbundle of $(B, \underline{\Omega})$. We now choose
\begin{equation}
	\label{Equation: choice of data}
	\begin{cases}
		\text{a torsion-free $T_{\mathbb C}Q$-connection \(\nabla'\) on \(B'\);}\\
		\text{a splitting $B' \to T_\mathbb{C}Q$ of the short exact sequence $0 \to A' \to T_\mathbb{C}Q \to B' \to 0$;}\\
		\text{a Lagrangian splitting $p^*B' \to B$,}
	\end{cases}
\end{equation}
and henceforth identify $p^*B'$ with its image in $B$. These choices induce compatible structures on the bundles
\begin{equation*}
	\mathcal F:=\widehat{\operatorname{Sym}}((B')^*)[[\hbar]] \quad \text{and} \quad W := \widehat{\operatorname{Sym}} (B^*)[[\hbar]]
\end{equation*}
as follows. On one hand, $\nabla'$ induces a connection \(\nabla^{\mathcal F}\) on $\mathcal F$. By \cite[Theorem 3.10 and (3.9)]{WanYau2025} (see also \cite[Lemmas 5.14 and 5.17]{PBWXuping}), the first two data in \eqref{Equation: choice of data} together with the trivial connection $d$ on the trivial line bundle $Q \times \mathbb{C}$ yield a flat connection on $\mathcal{F}$ of the form
\begin{equation*}
	\nabla^{K} =\nabla^{\mathcal F} + \Xi, \quad \Xi\in -dz^i \otimes \tfrac{\partial}{\partial z^i} + \Omega^1(Q,\widehat{\operatorname{Sym}}^{\geq 2}((B')^*) \otimes B').
\end{equation*}
Here, by abuse of notation, we use the same symbol \(\Xi\) for the corresponding \(1\)-form with values in derivations of $\widehat{\operatorname{Sym}}((B')^*)$. Furthermore, \cite[Proposition 3.15]{WanYau2025} gives a sheaf isomorphism $T$ between $\mathcal{O}_{A'}[[\hbar]]$ and the sheaf of $\nabla^K$-flat sections of $\mathcal{F}$.

On the other hand, the pullback connection \(p^*\nabla'\), together with the chosen Lagrangian splitting, induces a symplectic connection \(\nabla\) on \((M,\Omega)\), which in turn induces a connection \(\nabla^W\) on \(W\). In this case,
\begin{itemize}
	\item the $1$-form $\tilde{\Omega} \in \Omega_A^{1, 0}(M, B^*)$ satisfying $\delta=\frac{1}{\hbar}[\tilde{\Omega},\,\cdot\,]_*$ is given by $\tilde\Omega=v_i\,dz^i-u^i\,d\xi_i$;
	\item the $2$-form \(R^\nabla\in\Omega^2(M, \operatorname{Sym}^2 B^*)\) satisfying $\pi_{0, *}(R^\nabla) = 0$ and \((\nabla^W)^2=\tfrac{1}{\hbar} [R^\nabla,\,\cdot\,]_*\) is contained in $\Omega^2(M, p^*(B')^* \otimes p^*B')$. Under the identification $p^*(\operatorname{End} B') = p^*(B')^* \otimes p^*B'$, $R^\nabla$ coincides with $-p^*F^{\nabla'}$, where
	\begin{equation}
		\label{Equation: curvature of transverse tangent bundle}
		F^{\nabla'} \in \Omega^2(Q, \operatorname{End} B)
	\end{equation}
	is the curvature of the connection $\nabla'$ on $B'$.
\end{itemize}

The chosen Lagrangian splitting defines a faithful fibrewise Fock action
\[
\Gamma(M,W)\times\Gamma(M,p^*\mathcal F)
\to \Gamma(M,p^*\mathcal F),
\quad
(a,s)\mapsto\rho(a)s.
\]
On the two Lagrangian summands, it is characterized by
\[
\rho(u)s=u s,
\quad
\rho(v)s=\hbar\,\partial_v s,
\]
for \(u\in\Gamma(M,p^*(B')^*)\), \(v\in\Gamma(M,p^*B')\), and \(s\in\Gamma(M,p^*\mathcal F)\), where \(\partial_v\) denotes contraction with \(v\), extended as a derivation of \(p^*\mathcal F\). These operators satisfy the canonical commutation relations and hence extend uniquely, fibrewise, to an action of \((W,*)\) on \(p^*\mathcal F\).

The pullback $1$-form $p^*(\Xi)$ defines a $1$-form with values in derivations of $p^*\mathcal{F}$. It yields a $W$-valued $1$-form $\gamma$ satisfying
\[
\rho(\gamma)=\hbar \, p^*\Xi.
\]
Such a $1$-form $\gamma$ is unique by the faithfulness of the Fock representation $\rho$. Indeed,
\begin{equation}
	\gamma=-v_i\,dz^i+\gamma_{\geq3}, \quad \gamma_{\geq3} := (p^*\Xi_{\geq 2}^i)*v_i \in \operatorname{F}_3\Omega^1(M, W),
\end{equation}
where we write $\Xi + dz^i \otimes \tfrac{\partial}{\partial z^i}  = \Xi_{\geq 2}^i \otimes \tfrac{\partial}{\partial z^i}$ with local $\widehat{\operatorname{Sym}}^{\geq 2} (B')^*$-valued $1$-forms $\Xi_{\geq 2}^i$. This is due to the identity $\rho(a*v) = \rho(a) \rho(v) = \hbar a \partial_v$ for all $a \in \Gamma(M, p^*\mathcal{F})$ and $v \in \Gamma(M, p^*B')$.

\begin{lemma}
	\label{Lemma: intermediate gamma}
	The $1$-form $\gamma$ satisfies
	\[
	R^{\gamma} := R^\nabla+\nabla^W\gamma+\tfrac{1}{2\hbar}[\gamma,\gamma]_* = -\tfrac{\hbar}{2}p^*\operatorname{tr}F^{\nabla'},
	\]
	and the connection $D^{\gamma} := \nabla^W+\tfrac{1}{\hbar}[\gamma,\,\cdot\,]_*$ on $W$ is flat. Moreover, for all \(a\in \Gamma(M, W)\),
	\begin{equation}
		\label{Equation: intertwining relation}
		[p^*\nabla^K, \rho(a)] = \rho (D^{\gamma} a).
	\end{equation}
	
\end{lemma}
\begin{proof}
	We first compute $R^{\gamma}$ term by term. For \(u\in\Gamma(M,p^*(B')^*)\) and \(v\in\Gamma(M,p^*B')\), the Weyl-ordered Fock action of $uv$ is given by $\rho(uv) = \hbar ( u\partial_v + \tfrac{1}{2} (\partial_vu) \operatorname{Id} )$. More generally, for $a \in \Omega^*(M, p^*(B')^* \otimes p^*B') \cong \Omega^*(M, \operatorname{End} (p^*B'))$, one has $\rho(a) = \hbar (\widetilde{a} + \tfrac{1}{2} \operatorname{tr} a \operatorname{Id})$, where $\widetilde{a}$ is the extension of $a$ by derivations. Applying this identity to $R^\nabla = -p^*F^{\nabla'}$ yields
	\[
	\rho(R^\nabla)
	=
	\hbar(p^*\nabla^{\mathcal F})^2
	-
	\tfrac{\hbar}{2}
	p^*\operatorname{tr}F^{\nabla'}\,\operatorname{Id}.
	\]
	Using \(\rho(\gamma)=\hbar p^*\Xi\), we further obtain
	\begin{align*}
		\rho(\nabla^W\gamma)=&[p^*\nabla^{\mathcal F}, \rho(\gamma)] = \hbar [p^*\nabla^{\mathcal F}, p^*\Xi],\\
		\rho(\tfrac{1}{2\hbar} [\gamma,\gamma]_*)=&\tfrac{1}{2\hbar}[\rho(\gamma), \rho(\gamma)]=\tfrac{\hbar}{2}\,p^*[\Xi, \Xi].
	\end{align*}
	Combining the identities, together with using $(\nabla^{\mathcal F})^2+[\nabla^{\mathcal F}, \Xi]+\tfrac{1}{2} [\Xi, \Xi] = (\nabla^K)^2 = 0$ and faithfulness of $\rho$, gives \(R^{\gamma}=-\frac{\hbar}{2}p^*\operatorname{tr}F^{\nabla'}\). In particular, $(D^{\gamma})^2 = \tfrac{1}{\hbar} [R^{\gamma}, \,\cdot\,]_* = 0$.
	
	Since \(\nabla^W\) and \(\nabla^{\mathcal F}\) are induced by the same connection $\nabla'$ and \(\rho\) is a fiberwise Weyl algebra representation, the induced connections satisfy \([p^*\nabla^{\mathcal F}, \rho(a)]=\rho(\nabla^W a)\) for
	\(a\in \Gamma(M, W)\). The identities
	\(\rho(\gamma)=\hbar\,p^*\Xi\) and \( [\rho(\gamma), \rho(a)] = \rho([\gamma, a]_*) \) give $[p^*\nabla^K, \rho(a)] = \rho (D^{\gamma} a)$.
\end{proof}

In order to obtain a $1$-form $\gamma' \in -\tilde\Omega + \operatorname{F}_3\Omega^1(M, W)$ such that $\nabla^W + \tfrac{1}{\hbar}[\gamma', \,\cdot\,]_*$ is flat, we add a correction term to $\gamma \in v_i dz^i + \operatorname{F}_3\Omega^1(M, W)$ obtained by the following construction. The composition $T_\mathbb{C}M \overset{dp}{\longrightarrow} p^*T_\mathbb{C}Q \longrightarrow p^*B'$ is a surjective vector bundle map. Its kernel $A'' \subset T_\mathbb{C}M$ defines another transverse holomorphic structure on $M$, and is annihilated by the primitive $\Theta$ of $\Omega$ given in \eqref{Equation: primitive 1-form}. Thus, the connection $d + \Theta$ on the trivial line bundle $M \times \mathbb{C}$ restricts to a flat $A''$-connection. The splittings in the chosen data \eqref{Equation: choice of data} induces a splitting of $A'' \to T_\mathbb{C}M \to p^*B'$. With this splitting, the torsion-free connection $p^*\nabla'$ on $p^*B'$ and the connection $d + \Theta$, \cite[Theorem 3.10]{WanYau2025} gives a flat connection on $p^*\mathcal{F}$ of the form
\[
\nabla^{K, \Theta} = p^*\nabla^\mathcal{F} + \Xi^M + \Theta + F^\Theta,
\]
with $\Xi^M \in \Omega^1(M, \widehat{\operatorname{Sym}} \, p^*(B')^* \otimes p^*B' )$ and $F^\Theta \in \Omega^1(M, \widehat{\operatorname{Sym}}^{\geq 1} \, p^*(B')^*)$. Here, $\Xi^M$ and $\Theta + F^\Theta$ act on $p^*\mathcal{F}$ by derivations and multiplication, respectively. Due to the compatibility of the splittings,
\[
\Xi^M = p^*\Xi.
\]

We are now ready to prove our theorem.

\begin{proof}[\myproof{Theorem}{\ref{Theorem: lift of flat connections}}]
	We regard $F^\Theta$ as a $W$-valued $1$-form and set $\gamma' := \gamma + F^\Theta$,
	\[
	D^{\gamma'}:=\nabla^W+\tfrac{1}{\hbar}[\gamma',\,\cdot\,]_* \quad \text{and} \quad R^{\gamma'}:=R^\nabla+\nabla^W\gamma'+\tfrac{1}{2\hbar}[\gamma',\gamma']_*.
	\]
	We claim that
	\[ \gamma' \in -\tilde\Omega + \operatorname{F}_3\Omega^1(M, W) \quad \text{and} \quad R^{\gamma'} = -\Omega - \frac{\hbar}{2} p^* \operatorname{tr} F^{\nabla'}.
	\]
	Since \((D^{\gamma'})^2=\tfrac{1}{\hbar}[R^{\gamma'},\,\cdot\,]_*\), this claim implies that $D^{\gamma'
	}$ is a Fedosov connection on $W$ and the deformation quantization $(\mathcal{W}_{D^{\gamma'}}, \psi_{\gamma'})$ associated to $D^{\gamma'}$, given as in \eqref{Equation: DQ from Fedosov}, satisfies
	\[
	\operatorname{per}([(\mathcal{W}_{D^{\gamma'}}, \psi_{\gamma'})]) = -[\Omega] - \frac{\hbar}{2} p^* [\operatorname{tr} F^{\nabla'}] = \hbar \pi \sqrt{-1} p^*(c_1(B')).
	\]
	We first check the filtration condition. Write $F^\Theta = \sum_{r=1}^\infty F_{(r)}^\Theta$ with $F_{(r)}^\Theta \in \Omega^1(M, \widehat{\operatorname{Sym}}^r (p^*(B')^*))$. Let $\delta_{p^*B'}^{-1}$ denote the counterpart of the homotopy operator $\delta^{-1}$ with respect to $p^*B'$. Using the splitting of $A'' \to T_\mathbb{C}M \to p^*B'$, we identify $\delta_{p^*B'}^{-1}$ as an operator on $\Omega^*(M, p^*\mathcal{F})$. Since $\delta_{p^*B'}^{-1}$ and $\Xi^M = p^*\Xi$ annihilate $d\xi_i$, the recursive relations in \cite[Proposition 3.11]{WanYau2025} are simplified to
	\begin{align*}
		F_{(1)}^\Theta = \delta_{p^*B'}^{-1}\Omega = u^i d\xi_i, \quad F_{(r+1)}^\Theta = \delta_{p^*B'}^{-1}(p^*\nabla')F_{(r)}^\Theta, \quad r \geq 1,
	\end{align*}
	Moreover, for all $r \geq 0$, $F_{(r)}^\Theta = (F_{(r)}^\Theta)^i d\xi_i$ with local sections $(F_{(r)}^\Theta)^i$ of $\operatorname{Sym}^{r} (p^*(B')^*)$. Therefore, $F^\Theta \in u^id\xi_i + \operatorname{F}_3\Omega^1(M, W)$. By $\gamma \in -v_i\,dz^i+\operatorname{F}_3\Omega^1(M, W)$ and $\tilde\Omega=v_i\,dz^i-u^i\,d\xi_i$, we have
	\[
	\gamma' = \gamma + F^\Theta \in -\tilde{\Omega} + \operatorname{F}_3\Omega^1(M,W).
	\]
	Next, we compute $R^{\gamma'}=R^\gamma+D^\gamma F^\Theta +\frac1{2\hbar}[F^\Theta,F^\Theta]_*$. By Lemma \ref{Lemma: intermediate gamma}, $R^{\gamma}
	=-\frac{\hbar}{2}p^*\operatorname{tr}F^{\nabla'}$. The flatness of $p^*\nabla^{K}$ and $\nabla^{K,\Theta}$, together with $[\Theta + F^\Theta, \Theta + F^\Theta] = 0$, implies
	\[
	\Omega + D^{\gamma}F^\Theta
	=d\Theta+(p^*\nabla^K)F^\Theta = (p^*\nabla^K)(\Theta + F^\Theta) = 0.
	\]
	Note that $p^*\mathcal F$ is a commutative subalgebra bundle of $(W, *)$, thus $[F^\Theta,F^\Theta]_*=0$. Consequently, $R^{\gamma'} =
	-\Omega - \tfrac{\hbar}{2} p^* \operatorname{tr} F^{\nabla'}$. Therefore, our claim holds.
	
	We proceed to construct the module structure. For all open $U \subset M$, $s \in (\mathcal{W}_{D^{\gamma'}})(U)$ and $f \in \mathcal{O}_{A'}(j^{-}(U))[[\hbar]]$, we define \(s\cdot f \in \mathcal{O}_{A'}(j^{-1}(U))[[\hbar]]\) uniquely by
	\[
	T(s\cdot f) = \rho(j^*s)\,T(f),
	\]
	where $T$ is the map that identifies $\mathcal O_{A'}[[\hbar]]$ with the sheaf of $\nabla^K$-flat sections of $\mathcal F$ given in \cite[Proposition 3.15]{WanYau2025}. To see that $s\cdot f$ is well defined, we show that \(\rho(j^*s)\,T(f)\) is \(\nabla^K\)-flat. Indeed, because $F_{(r)}^\Theta = (F_{(r)}^\Theta)^i d\xi_i$ for all $r \geq 0$, $j^*F^\Theta = 0$. Restricting the intertwining relation \eqref{Equation: intertwining relation} along \(j\) gives $[\nabla^K, \rho(j^*s)] = \rho(j^*(D^\gamma s)) = \rho(j^*(D^{\gamma'} s)) = 0$. Therefore,
	\[
	\nabla^K ( \rho(j^*s)\, T(f) ) = \rho(j^*s) \, (\nabla^K T(f)) = 0.
	\]
	If \(s,s'\in \mathcal{W}_{D^{\gamma'}}(U)\), then \((s* s')\cdot f=s\cdot(s'\cdot f)\) because
	\[
	T((s* s')\cdot f)
	=\rho(j^*(s * s')) T(f)
	=\rho(j^*s)\rho(j^*s')T(f)
	=T(s\cdot(s'\cdot f)).
	\]
	Therefore,
	\begin{equation}
		\label{Equation: module structure}
		\mathcal{W}_{D^{\gamma'}} \otimes_{\mathbb{C}[[\hbar]]} j_*\mathcal O_{A'}[[\hbar]] \to j_*\mathcal O_{A'}[[\hbar]], \quad (s, f) \mapsto s \cdot f,
	\end{equation}
	defines a $\mathcal{W}_{D^{\gamma'}}$-module structure on $j_*\mathcal{O}_{A'}[[\hbar]]$.
	
	We also denote by \(\pi_0:\mathcal F\to\widehat{\operatorname{Sym}}((B')^*)\) the reduction modulo \(\hbar\). Modulo \(\hbar\), the operator \(\rho(j^*s)\) reduces to multiplication by \(T(j^*(\psi_{\gamma'}(s)))\). Since \(T\) commutes with reduction modulo \(\hbar\), we obtain
	\[
	T\bigl(\pi_0(s\cdot f)\bigr)=
	\pi_0\bigl(\rho(j^*s)T(f)\bigr)=
	T\bigl( j^*(\psi_{\gamma'}(s))\bigr)
	\,\pi_0\bigl(T(f)\bigr)=
	T \bigl(j^*(\psi_{\gamma'}(s))\,\pi_0(f)\bigr).
	\]
	By the injectivity of \(T\), it follows that $\pi_0(s\cdot f) = j^* \bigl(\psi_{\gamma'}(s)\bigr)\,\pi_0(f)$, which is exactly the commutativity of the required diagram.
\end{proof}

To extend the above construction to a more general \(A\)-model setting, we introduce two notions.

\begin{definition}\label{transverse holomorohic lagrangian}
	Let $(C, A_C, \Omega)$ be a transverse holomorphic symplectic manifold with characteristic distribution $A_{C, \mathbb{R}}$. A submanifold $Q \subset C$ is called a \emph{transverse holomorphic Lagrangian submanifold} if $A_C \cap T_\mathbb{C} Q$ defines a transverse holomorphic structure on $Q$ and $(TQ)^{\perp \Omega_C} = TQ + A_{C, \mathbb{R}} \vert_Q$.
\end{definition}

\begin{remark}
	Unlike in the holomorphic symplectic setting, the dimension of a transverse holomorphic Lagrangian submanifold is not determined by that of the ambient manifold. More precisely,
	\begin{equation}
		\label{Equation: dimension counting}
		\dim Q = \operatorname{rank} ((A_{C, \mathbb{R}} \cap TQ) \otimes \mathbb{C}) + \operatorname{rank} (T_\mathbb{C}C / A_C),
	\end{equation}
	where $\dim Q$ denotes the real dimension of $Q$, while $\operatorname{rank}$ denotes complex rank.
\end{remark}

\begin{definition}\label{Definition5.11}
	A \emph{transverse holomorphic Lagrangian fibration} is a smooth fibration $\pi:C\to Q$ from a transverse holomorphic symplectic manifold $(C, A_C, \Omega_C)$ to a transverse complex manifold $(Q, A')$ such that $d\pi(A_C) \subset A'$ and every fibre of $\pi$ is a transverse holomorphic Lagrangian submanifold of $(C, \Omega_C)$.
\end{definition}

Let \((X,\omega)\) be a symplectic manifold, \((C, A_C, \Omega_C)\) be the underlying transverse holomorphic symplectic manifold of a coisotropic \(A\)-brane on $(X, \omega)$, and \(Q\) be a Lagrangian submanifold of \(X\) contained in \(C\). 
Let $\iota: Q \to C$ be the inclusion. Since $Q\subset C$, one has
\[
(TC)^{\perp_\omega} \vert_Q \subset (TQ)^{\perp_\omega} \vert_Q =TQ.
\]

\begin{remark}
	Suppose $A' := A_C \cap T_\mathbb{C}Q$ defines a transverse holomorphic structure on $Q$. Then it automatically satisfies the adaptedness condition in the sense of \cite[Definition 5.2]{qcr}. In this case, simple linear algebra (see \cite[Proposition A.1]{qcr}) shows that
	\[
	(TQ)^{\perp_{\Omega_C}} = (TQ)^{\perp_\omega} = TQ = TQ + (TC)^{\perp_\omega} \vert_Q,
	\]
	hence $Q$ is a transverse holomorphic Lagrangian submanifold of $(C, \Omega_C)$.
\end{remark}

\begin{proposition}\label{Prop5.9}
	Suppose that $Q$ is a transverse holomorphic Lagrangian submanifold of $(C, \Omega_C)$ and $C$ admits a transverse holomorphic Lagrangian fibration $\pi:C\rightarrow Q$ such
	that $\pi \circ \iota =\operatorname{Id}_Q$. Let $(M, \Omega, p, j)$ be given as in \eqref{Equation: transverse cotangent bundle}. Then there exist neighborhoods $U$ of $j(Q)$ in $M$ and $V$ of $Q$ in $C$, together with a diffeomorphism $\Phi: U \to V$, such that
	\[
	\Phi \circ j = \iota, \quad p = \pi \circ \Phi \quad \text{and} \quad \Phi^*(\Omega_C \vert_V) = \Omega \vert_U.
	\]
	
\end{proposition}
\begin{proof}
	Since every fibre of $\pi$ is transverse holomorphic Lagrangian, $\ker d\pi$ is isotropic with respect to $\Omega_C$. We therefore obtain a well-defined real vector bundle morphism
	\[
	J: \ker d\pi \to \pi^*(T_\mathbb{C}^*Q)
	\]
	as follows. For $m \in C$, $u\in\ker d\pi_m$ and $v \in T_{\mathbb{C}, \pi(m)}Q$, set
	\begin{equation}
		\label{Equation: definition of J}
		(J(u))(v):=\Omega_C(u, \widetilde v),
	\end{equation}
	where $\widetilde v\in T_{\mathbb C,m}C$ is any lift of $v$ under $d\pi$.
	
	Let $A' = A_C \cap T_\mathbb{C}Q$ and $B' = T_\mathbb{C}Q/A'$. We claim that $J$ is a real vector bundle isomorphism onto $\pi^*(B')^*$. First, $J(u)$ annihilates $A'$. Indeed, by Lemma \ref{Lemma: transverse holomorphic Lagrangian fibration}, for all $v \in A'_{\pi(m)}$, we may choose $\widetilde{v} \in A_{C, m}$ such that $v = d\pi(\widetilde{v})$, yielding $J(u)(v) = \Omega_C(u, \widetilde{v}) = 0$. Therefore $J(u) \in (B')_{\pi(m)}^*$. Next, suppose $J(u) = 0$. Then $\Omega_C(u, \,\cdot\,)$ annihilates $T_mC$ and thus $u \in (A_{C, \mathbb{R}})_m$. Using $A_{C, \mathbb{R}} \cap \ker d\pi = 0$ by Lemma \ref{Lemma: transverse holomorphic Lagrangian fibration}, we obtain $u = 0$. Hence, $J$ is injective. Finally, $\ker d\pi$ and $\pi^*B'$ have the same real rank by the dimension computation in Lemma \ref{Lemma: transverse holomorphic Lagrangian fibration}; hence $J: \ker d\pi \to \pi^*B'$ is an isomorphism.
	
	For $x \in Q$, write $C_x := \pi^{-1}(\{x\})$. Recall that $M = (TQ/A'_\mathbb{R})^*$, where $A'_\mathbb{R}$ denote the characteristic distribution of $(Q, A')$, and that $M \to (B')^*$ given by $\xi \mapsto \xi^{1, 0}$ is a vector bundle isomorphism. Each $\xi \in M$ defines a vector field $X^\xi$ on $C_{p(\xi)}$ determined by
	\begin{equation}
		\label{Equation: vector field on fibre}
		J(X^\xi(m)) = -2\sqrt{-1} \xi^{1, 0}, \quad \text{for all } m \in C_{p(\xi)}.
	\end{equation}
	
	The assignment $\xi\mapsto X^\xi(x)$ identifies $M_x$ with the normal space of $Q$ in $C$ at $x$. Hence, by the standard tubular-neighborhood argument, after shrinking around the zero section there exist neighborhoods $U\subset M$ of $j(Q)$ and $V\subset C$ of $Q$ such that
	\[
	\Phi: U \to V, \quad \xi \mapsto \varphi_1^\xi(p(\xi)),
	\]
	defines a diffeomorphism satisfying $\Phi\circ j=\iota$ and $\pi\circ\Phi=p$, where $\varphi_t^\xi$ denotes the local flow of the vector field $X^\xi$ on $C_{p(\xi)}$.
	
	It remains to prove the claim that the $2$-form $\Lambda:=\Phi^*(\Omega_C \vert_V) -\Omega \vert_U$ is horizontal for the projection $p \vert_U: U \to Q$. If this claim holds, then $d$-closedness of $\Lambda$ implies that $\Lambda$ is basic and hence, shrinking $U$ and $V$ if necessary, $\Lambda = p^*\beta$ for some complex $2$-form $\beta$ on $Q$. We obtain $\beta=j^*\Lambda = 0$, because $j^*\Omega = 0$ by construction and $j^*\Phi^*(\Omega_C \vert_V) = \Omega_C \vert_Q = 0$ by the hypothesis that $Q$ is a transverse holomorphic Lagrangian submanifold of $(C, \Omega_C)$. Consequently,
	\[
	\Phi^*(\Omega_C \vert_V)=\Omega \vert_U.
	\]
	
	Now, we prove this last claim. Fix a vertical tangent vector $u$ at a point $\xi \in U$ and a tangent vector $v \in T_{\xi}U$. Set $x = p(\xi)$. Then there exists unique $\eta \in M_x$ such that $u = \left.\frac{d}{dt}\right|_{t=0}(\xi+t\eta)$. We compute $\Omega(u, v)$ using the primitive $\Theta$ of $\Omega$ given in \eqref{Equation: primitive 1-form}. Choose commuting local vector fields $Y$ and $Z$ extending $u$ and $v$, with $Y$ being vertical.
	Then $\Theta(Y) = \Theta([Y, Z]) = 0$, yielding
	\[
	\Omega(u, v) = \mathcal{L}_u (\Theta(Z)) = -2\sqrt{-1} \tfrac{d}{dt} (\xi + t \eta)^{1, 0}(dp(Z)) \vert_{t=0} = -2\sqrt{-1} \eta^{1, 0}(dp(v)).
	\]
	We next compute $(\Phi^*\Omega_C)(u, v)$. Choose local $d$-closed sections $\alpha$ and $\beta$ of $(B')^*$ near $x$ extending $-2\sqrt{-1}\xi^{1, 0}$ and $-2\sqrt{-1} \eta^{1, 0}$. Then $X^\xi$ and $X^\eta$ extend to local vertical vector fields on $C$ satisfying $\iota_{\widetilde{X}^\xi} \Omega_C = \pi^*\alpha$ and $\iota_{\widetilde{X}^\eta} \Omega_C = \pi^*\beta$. Since $\alpha, \beta$ and $\Omega_C$ are all $d$-closed, a direct computation shows
	\[
	\iota_{[\widetilde{X}^\xi, \widetilde{X}^\eta]} \Omega_C = \mathcal{L}_{\widetilde{X}^\xi} \iota_{\widetilde{X}^\eta} \Omega_C - \iota_{\widetilde{X}^\eta} \mathcal{L}_{\widetilde{X}^\xi} \Omega_C = d(\iota_{\widetilde{X}^\xi} \pi^*\alpha).
	\]
	As $\widetilde{X}^\xi$ is a vertical vector field, one has $\iota_{\widetilde{X}^\xi} \pi^*\alpha = 0$ and hence $\iota_{[\widetilde{X}^\xi, \widetilde{X}^\eta]} \Omega_C = 0$. Thus, the vector fields $X^\xi$ and $X^\eta$ on $C_x$ commute. It implies that $\varphi_1^{\xi + t\eta}(x) = \varphi_t^\eta(\varphi_1^\xi(x)) = \varphi_t^\eta(m)$, where $m:=\Phi(\xi) \in C_x$, yielding $d\Phi(u) = X^\eta(m)$. Observe that $d\pi(d\Phi(v)) = d(\pi \circ \Phi)(v) = dp(v)$. Using \eqref{Equation: definition of J} and \eqref{Equation: vector field on fibre}, we obtain
	\[
	(\Phi^*\Omega_C)(u, v)
	=\Omega_C(X^\eta(m), d\Phi(v)) = J(X^\eta(m))(dp(v)) = -2\sqrt{-1}\eta^{1, 0}(dp(v)).
	\]
	Consequently, $\Omega(u, v) = (\Phi^*\Omega_C)(u, v)$, implying that $\Lambda$ is horizontal.
\end{proof}

\begin{lemma}
	\label{Lemma: transverse holomorphic Lagrangian fibration}
	Under the hypothesis of Proposition \ref{Prop5.9}, $A_{C, \mathbb{R}} \cap \ker d\pi = 0$ and $d\pi: T_\mathbb{C}C \to \pi^*T_\mathbb{C}Q$ restricts to a surjective vector bundle morphism $A_C \to \pi^* A'$.
\end{lemma}
\begin{proof}
	The characteristic distribution $A_{C, \mathbb{R}}$ of $(C, A_C)$ is precisely $(TC)^{\perp_\omega}$. Since $A_{C, \mathbb{R}} \vert_Q \subset TQ$, one has $A'_\mathbb{R} = A_{C, \mathbb{R}} \vert_Q$ and thus $\dim Q = \operatorname{rank} A_C$ by \eqref{Equation: dimension counting}. For each $x \in Q$, $C_x := \pi^{-1}(\{x\})$ is a transverse holomorphic Lagrangian submanifold of $(C, \Omega_C)$ and
	\[
	\dim C_x = \dim C - \dim Q = \operatorname{rank} (T_\mathbb{C}C / A_C).
	\]
	By \eqref{Equation: dimension counting} again, $A_{C, \mathbb{R}} \cap T C_x = 0$. It implies that $A_{C, \mathbb{R}} \cap \ker d\pi = 0$.
	
	Since $d\pi(A_C) \subset \pi^*A'$, we have $d\pi(A_{C, \mathbb{R}}) \subset A'_\mathbb{R}$. Using $A_{C, \mathbb{R}} \cap \ker d\pi = 0$ and that $A_{C, \mathbb{R}}$ and $A'_\mathbb{R}$ have the same real rank, we obtain $d\pi(A_{C, \mathbb{R}}) = \pi^*A'_\mathbb{R}$. Now, suppose $m \in C$ and $u \in A'_{\pi(m)}$. Pick $v \in A_{C, m}$ and $w \in \overline{A_{C, m}}$ such that $d\pi(v+w) = u$. Then $d\pi(w) = u - d\pi(v)$ lies in $A'_{\pi(m)} \cap \overline{A'_{\pi(m)}} = (A'_{\mathbb{R}})_{\pi(m)} \otimes \mathbb{C}$, so there exists $v' \in A_{C, \mathbb{R}, m} \otimes \mathbb{C}$ such that $d\pi(v') = d\pi(w)$. We obtain $v + v' \in A_{C, m}$ and $u = d\pi(v + v')$. This yields $d\pi(A_C)=\pi^*A'$.
\end{proof}

Under the above assumptions, the local identification of $(C,\Omega_C)$ with $(M,\Omega)$ near $Q$ allows us to apply the preceding construction to obtain the corresponding perturbative module action.

\appendix

\subsection*{Acknowledgement}
\quad\par

The work of N. C. Leung was substantially supported by grants of the Hong Kong Research Grants Council (Project No. 14305923, 14302224, 14302225).
The work of the second author was supported by World Premier International Research Center Initiative (WPI), MEXT, Japan, and Yau Mathematical Sciences Center, Tsinghua University, China. He is grateful to Mathieu Sti\'enon and Ping Xu for their valuable discussions on formality theorems for Lie pairs.
The third author thanks Ki Fung Chan and Eddie Lam for many useful discussions. He also thanks the first and second authors for introducing him to this field and for their guidance and support throughout this work.

\bibliographystyle{amsplain}
\bibliography{References}

\end{document}